\documentclass[10pt]{article}
\usepackage{amsfonts,color}
 \usepackage{amsmath}
 \usepackage{epsfig}
\usepackage{lscape}
\usepackage{amssymb}
\usepackage{graphicx}
 \usepackage{footnote}
\usepackage[utf8]{inputenc}
\usepackage{mathtools}
\usepackage{amsthm,wasysym}
\usepackage[english]{babel}
\usepackage{tikz}
\usepackage{bbm}
\usepackage{colonequals}
 \usepackage{cancel}

\allowdisplaybreaks[1]

\makeindex

\def\12{{\frac{1}{2}}}                                                             

\newcommand{\f}[1]{{{#1}}}                                      

\newcommand{\e}[1]{{\emph{#1}}}
\newcommand{\fg}[1]{{\mbox{$#1$}}}                             

\newcommand{\R}{\mathbb{R}}

\newcommand{\norm}[1]{\|\, #1 \, \|}                      
\newcommand{\Cof}{\rm Cof \,}              
\newcommand{\m}{\mathfrak{m}}                      
\newcommand{\n}{\mathfrak{n}}

\DeclareMathOperator{\sym}{sym}
\DeclareMathOperator{\tr}{tr}

\DeclareMathOperator{\axl}{axl}
\DeclareMathOperator{\anti}{anti}
\DeclareMathOperator{\dev}{dev}
\DeclareMathOperator{\sL}{\mathfrak{sl}}
\DeclareMathOperator{\so}{\mathfrak{so}}
\DeclareMathOperator{\gl}{\mathfrak{gl}}

\DeclareMathOperator{\Lin}{Lin}

\DeclareMathOperator{\Curl}{Curl\,}
\DeclareMathOperator{\curl}{curl}
\newcommand{\Sym}{ {\rm{Sym}} }

\newcommand{\id}{{{\mathbbm{1}}}}

\def\Div{\textrm{Div\,}}
\def\DIV{\textrm{DIV\,}}
\newcommand{\gradx}[1]{{\rm grad}_{\f x}[{#1}]}
\newcommand{\gradxi}[1]{{\rm grad}_{\fg \xi}[{#1}]}
\newcommand{\Gradx}[1]{{\rm Grad}_{\f x}[{#1}]}
\newcommand{\Gradxi}[1]{{\rm Grad}_{\fg \xi}[{#1}]}
\newcommand{\GRADx}[1]{{\rm GRAD}_{\f x}[{#1}]}
\newcommand{\GRADxi}[1]{{\rm GRAD}_{\fg \xi}[{#1}]}
\newcommand{\divx}{{\rm div}_{\f x}}
\newcommand{\divxi}{{\rm div}_{\fg \xi}}
\newcommand{\Divx}{{\rm Div}_{\f x}}
\newcommand{\Divxi}{{\rm Div}_{\fg \xi}}
\newcommand{\DIVx}{{\rm DIV}_{\f x}}
\newcommand{\DIVxi}{{\rm DIV}_{\fg \xi}}
\newcommand{\curlx}{{\rm curl}_{\f x}}
\newcommand{\curlxi}{{\rm curl}_{\fg \xi}}
\newcommand{\Curlx}{{\rm Curl}_{\f x}}
\newcommand{\Curlxi}{{\rm Curl}_{\fg \xi}}
\newcommand{\inc}{{\rm inc}}
\newcommand{\incx}{{\rm inc}_{\f x}}
\newcommand{\incxi}{{\rm inc}_{\fg \xi}}

\newcommand{\ux}{\f u (\f x)}
\newcommand{\uxi}{\f u^\sharp (\fg \xi)}
\newcommand{\ks}{\widetilde{\f k}}
\newcommand{\ksx}{\widetilde{\f k}(\f x)}
\newcommand{\ksxi}{\widetilde{\f k}^{\sharp}(\fg \xi)}
\newcommand{\nablaxi}{\fg \nabla_{\! \fg \xi}}
\newcommand{\nablax}{\fg \nabla_{\! \f x}}
\newcommand{\sigmax}{\fg \sigma(\f x)}
\newcommand{\sigmaxi}{\fg \sigma^{\sharp}(\fg \xi)}
\newcommand{\phix}{\fg \varphi(\f x)}
\newcommand{\phixi}{\fg \varphi^{\sharp}(\fg \xi)}
\newcommand{\Cx}{\f C(\f x)}
\newcommand{\Cxi}{\f C^{\sharp}(\fg \xi)}
\newcommand{\varepsilonx}{\fg \varepsilon(\f x)}
\newcommand{\varepsilonxi}{\fg \varepsilon^{\sharp}(\fg \xi)}
\newcommand{\Px}{\f P(\f x)}
\newcommand{\Pxi}{\f P^{\sharp}(\fg \xi)}
\newcommand{\Fx}{\f F(\f x)}
\newcommand{\Fxi}{\f F^{\sharp}(\fg \xi)}

\newcommand{\oQ}{\overline{Q}}

\newcommand{\Laplace}{\mathop{}\!\mathbin\bigtriangleup}

\newcommand{\SO}{{\rm SO}}

\def\skew{\, \text{skew} \, }

\newcommand{\di}{\,{\rm d}}

\newcommand{\Grad}[1]{{\rm Grad}[{#1}]}
\newcommand{\GRAD}[1]{{\rm GRAD}[{#1}]}
\newcommand{\scal}[1]{\langle#1\rangle}            

\DeclareMathOperator{\Euklid}{\mathbb{R}^3}        
\DeclareMathOperator{\Skalar}{\mathbb{R}}           
\newcommand{\B}{{\cal B}}                          

\def\C{\mathbb{C}}

\def\L{\mathbb{L}}
\makeatletter
\let\@fnsymbol\@arabic

\begin{document}

\title{Rotational invariance conditions in elasticity, gradient elasticity and its connection to isotropy
}
\author{\normalsize{Ingo M\"unch\thanks{Corresponding author: Ingo M\"unch, Institute for Structural Analysis, Karlsruhe Institute of Technology, Kaiserstr. 12, 76131 Karlsruhe, Germany, email: ingo.muench@kit.edu}
\quad and \quad
Patrizio Neff\thanks{Patrizio Neff,  \ \ Head of Lehrstuhl f\"{u}r Nichtlineare
Analysis und Modellierung, Fakult\"{a}t f\"{u}r Mathematik, Universit\"{a}t Duisburg-Essen,
Thea-Leymann Str. 9, 45127 Essen, Germany, email: patrizio.neff@uni-due.de}
%
%
}
} 
\maketitle

\begin{abstract}
For homogeneous higher gradient elasticity models we discuss frame-indifference and isotropy requirements. To this end, we introduce the notions of local versus global SO(3)-invariance and identify frame-indifference (traditionally) with global left SO(3)-invariance and isotropy with global right SO(3)-invariance. For specific restricted representations, the energy may also be local left SO(3)-invariant as well as local right SO(3)-invariant. Then we turn to linear models and consider a consequence of frame-indifference together with isotropy in nonlinear elasticity and apply this joint invariance condition to some specific linear models. The interesting point is the appearance of finite rotations in transformations of a geometrically linear model. It is shown that when starting with a linear model defined already in the infinitesimal symmetric strain $\fg \varepsilon = \sym \Grad{\f u}$, the new invariance condition is equivalent to isotropy of the linear formulation. Therefore, it may be used also in higher gradient elasticity models for a simple check of isotropy and for extensions to anisotropy. In this respect we consider in more detail variational formulations of the linear indeterminate couple stress model, a new variant of it with symmetric force stresses and general linear gradient elasticity.
\\
\vspace*{0.25cm}

\noindent AMS Math 74A10 (Stress), 74A35 (polar materials), 74B05 (classical linear elasticity), 74A30 (nonsimple materials), 74A20 (theory of constitutive functions), 74B20 (nonlinear elasticity)
\end{abstract}

\newpage

\tableofcontents

%
\section{Introduction}\label{KapIntro}
This paper is motivated by our endeavour to better understand isotropy conditions in higher gradient elasticity theories. Therefore,  we make use of the old idea expressed by Truesdell \cite[Lecture 6]{Truesdell66} albeit including higher gradient continua.\footnote{Truesdell writes \cite[Lecture 6]{Truesdell66}: "In the older literature a material is said to be "isotropic" if it is "unaffected by rotations." This means that if we first rotate a specimen of material and then do an experiment upon it, the outcome is the same as if the specimen had not been rotated. In other words, within the class of effects considered by the theory, [referential] rotations cannot be detected by any experiment. The response of the material with respect to the reference configuration $\kappa$ is the same as that with respect to any other obtained from it by rotation."} To this end we reconsider frame-indifference and isotropy first in nonlinear elasticity and then in classical linear elasticity theories. These investigations are extended to linear higher gradient elasticity with a focus on the indeterminate couple stress theory. We do not discuss general material symmetry conditions a la Noll's group theoretic framework, as e.g. presented in \cite{ElzanowskiEpstein92,Steigmann07b}.

In nonlinear elasticity, isotropy requirements are fundamentally different from frame-indifference requirements. Indeed, isotropy rotates the referential coordinate system, while frame-indifference rotates the spatial coordinate system. Both spaces, referential and spatial, are clearly independent of each other and the corresponding rotations of frames are not connected to each other.

In this paper, however, we will later apply the {\bf same rotation} to the referential and spatial frame simultaneously. It is clear that this is not equivalent to frame-indifference or isotropy as we show with simple explicit examples. What, then, is the use of such a specific transformation (which we will call the $\sharp$-transformation in Section \ref{KapRotInvariance}). Indeed, once we consider reduced energy expressions which already encode frame-indifference, then applying the $\sharp$-transformation and requiring form-invariance of the energy under this transformation is equivalent to isotropy. Surprisingly, the $\sharp$-transformation remains operative in exactly the same way when applying it to the geometrically linear context. It is this insight that we follow when discussing isotropy conditions for higher gradient continua. We do not believe that our development leads to new results for the representation of isotropic formulations,\footnote{Indeed, in Murdoch \cite[eq.(16)-(20)]{Murdoch77} we find the same conclusions regarding isotropy of a material point.} but we believe that there is a conceptual gain of understanding when working directly with transformation connected to finite rotations. In this respect we remind the reader that for linearized frame-indifference it would be sufficient to work with the Lie-algebra of skew-symmetric rotations instead of the orthogonal group. This is, however, not possible with respect to isotropy considerations.\\

Further, we give a contribution to the discussion in \cite{Lederer15b}, where a fundamental aspect of strain gradient elasticity as proposed by Mindlin \cite{Mindlin63} is doubted. There we read: "However, the approach described by equation (1) includes a serious flaw:\footnote{This remark refers to eq.(1), Mindlin \cite{Mindlin63}, which reads, in the notation of the present paper,
$ \frac{\partial \, m_{31}}{\partial \, x_1} + \frac{\partial \, m_{32}}{\partial \, x_2} + \sigma_{21} - \sigma_{12} = 0$. Here, $\f m$ is the second-order couple stress tensor and $\sigma_{21} - \sigma_{12} = 2 \, \axl(\skew \widetilde{\sigma})$ represents the skew-symmetric part of the total force stress tensor $\widetilde{\sigma}$. The equation is given for the planar situation, see also eq.\eqref{CoupleStressBalanceMomentum}.
} Couple stresses are components of a tensor which has a rank higher than two. In conclusion, the torque arising from an asymmetric stress tensor cannot consequently be compensated by the torque associated with a tensor of higher rank. In fact, there are different transformation rules for tensors of different rank during rotation of a coordinate system. Therefore, couple stress theory does not consequently fulfill the requirement of invariance with respect to rotation of the coordinate system. In spite of this deficiency, strain gradient theory still contains important aspects which are corrected." Motivated by this statement we show to the contrary that all proper invariance requirements are satisfied.\\

The paper is now structured as follows. After notational agreements we recall the gradient of continuum rotation, which is the curvature measure in linear couple stress-elasticity. Identities from the scalar triple product built a basis to study transformation rules for this curvature measure later in the text. In Section \ref{KapObjNonlinearElas} we start our investigation in the context of nonlinear hyperelasticity by discussing frame-indifference and isotropy for first and second gradient continua. Then we consider the simultaneous transformation of the deformation to new spatial and referential coordinates. The consequence of this transformation in linearized elasticity follows in Section \ref{KapObjLinearElas}. Further, we prepare for higher gradient elasticity by giving details on transformation rules and apply them to the linear momentum balance equations of several models in Section \ref{KapRotInvariance}. We then define form-invariance in higher gradient elasticity in Section \ref{KapObjectIso} and specify our result for several couple stress theories in Section \ref{KapIndeterminate}. Finally, we conclude and give an outlook.
\subsection{Notational agreements and preliminary results}\label{KapNotations}
By $\R^{3\times 3}$ we denote the set of real $3\times 3$ second order tensors, written
with capital letters. Vectors in $\R^{3}$ are denoted by small letters. Additionally, tensors in $\R^{3\times 3 \times 3}$ are necessary for our discussion, where significant symbols like $\fg \epsilon$, $\m$, and $\n$ will be used. The components of vectors and tensors are given according to orthogonal unit vectors $\f e_1,
\, \f e_2, \, \f e_3$, which may be rotated by $\f Q \in \SO(3)$. Throughout this paper Latin subscripts specify the direction of components in index notation and take the values $1,2,3$. For repeating subscripts Einstein's
summation convention applies.

We adopt the usual abbreviations of Lie-algebra theory, i.e.,  $\so(3):=\{\f X\in\mathbb{R}^{3\times3}\;|\f X^T=- \f X\}$ is the Lie-algebra of  skew symmetric tensors
and $\sL(3):=\{\f X\in\mathbb{R}^{3\times3}\;| \tr({\f X})=0\}$ is the Lie-algebra of traceless
tensors. For all $\f X\in\mathbb{R}^{3\times3}$ we set $\sym \f X=\frac{1}{2}(\f X^T+\f X)\in\Sym(3)$,
$\skew \f X=\frac{1}{2}(\f X-\f X^T)\in \so(3)$ and the deviatoric part $\dev \f X=\f
X-\frac{1}{3}\;\tr(\f X) \,\id\in \sL(3)$ and we have the \emph{orthogonal Cartan-decomposition
of the Lie-algebra} $\gl(3)$
\begin{align}\label{Cartan}
\gl(3)&=\{\sL(3)\cap \Sym(3)\}\oplus\so(3) \oplus\mathbb{R}\!\cdot\! \id \, ,\qquad \f
X=\dev \sym \f X+ \skew \f X+\frac{1}{3}\tr(\f X) \, \id \, ,
\end{align}
simply allowing to split every second order tensor uniquely into its trace-free symmetric part, skew-symmetric part and spherical part, respectively. For
\begin{align}
\overline{\f A}  \colonequals \left(\begin{array}{ccc}
0 &-a_3&a_2\\
a_3&0& -a_1\\
-a_2& a_1&0
\end{array}\right)\in \so(3) \, , \qquad
\overline{\f A} \, \f v= \f a \times \f v \, , \qquad \forall \, \f v\in\mathbb{R}^3 \, , \qquad  \f a = \axl [\overline{\f A}] \,
, \qquad \overline{\f A} = \anti [\f a] \, ,
\end{align}
the functions $\axl:\so(3)\rightarrow\mathbb{R}^3$ and $\anti:\mathbb{R}^3\rightarrow
\so(3)$ are given by
\begin{align}\label{axlanti}
\axl [\overline{\f A}]  \colonequals -\frac{1}{2} \, \overline{A}_{ij} \, \epsilon_{ijk} \, \f e_k \, ,
\quad \quad \anti[\f a] := - \epsilon_{ijk} \, a_k \, \f e_i \otimes \f e_j = \overline{A}_{ij} \, \f e_i \otimes \f e_j \, ,
\end{align}
 with $\epsilon_{ijk}=+1$ for even permutation, $\epsilon_{ijk}=-1$ for odd permutation, and $\epsilon_{ijk}=0$ else.
 Note that the skew-symmetric part of a tensor $\f X \in \Skalar^{3 \times 3}$ can be written by  the combination of permutations via
\begin{align}\label{SkewTensor}
\skew \f X = \12 \, \epsilon_{nij} \, \epsilon_{nab} \, X_{ab} \, \f e_i \otimes \, \f e_j \,.
\end{align}
From  eq.\eqref{SkewTensor} and eq.\eqref{axlanti}$_1$ one obtains
\begin{align}\label{axlSkewTensor}
\axl [ \skew \f X ] = - \12 \, \underbrace{\epsilon_{ijk} \, \12 \, \epsilon_{nij}}_{\displaystyle \delta_{nk}} \, \epsilon_{nab} \, X_{ab} \, \f e_k = - \12 \, X_{ab} \, \epsilon_{abk} \, \f e_k
\,.
\end{align}
%
In index notation the typical conventions for differential operations are implied such as comma followed by a subscript to denote the partial derivative with respect to  the corresponding coordinate. The gradient of a scalar field $\phi \in \Skalar$ and the gradient of a vector field $\f b \in \Euklid$ are given by
\begin{align}\label{GradSkalDef}
\gradx{\phi}  \colonequals \frac{\partial \, \phi}{\partial \, \f x} = \frac{\partial \, \phi}{\partial \, x_i} \, \f e_i = \phi_{,i} \, \f e_i \quad \in \Euklid \,,
\end{align}
\begin{align}\label{GradVektDef}
\Gradx{\f b}  \colonequals \frac{\partial \, \f b}{\partial \, \f x} = \frac{\partial \, b_i}{\partial \, x_j} \, \f e_i \otimes \f e_j = b_{i,j} \, \f e_i \otimes \f e_j \quad \in \Skalar^{3 \times 3} \,.
\end{align}
Similarly, we also define the gradient of a tensor field $\f X\in\mathbb{R}^{3\times3}$
:
\begin{align}\label{GradTensor2Def}
\GRADx{\f X}  \colonequals \frac{\partial \, \f X}{\partial \, \f x} = \frac{\partial \, X_{ij}}{\partial \, x_k} \, \f e_i \otimes \f e_j \otimes \f e_k = X_{ij,k} \, \f e_i \otimes \f e_j \otimes \f e_k \quad \in \Skalar^{3 \times 3 \times 3} \,.
\end{align}
%
%
%
For vectors $\f a, \f b\in\R^3$ we let
\begin{align}\label{VektInnerProduct}
\scal{\f a,\f b}_{\R^3} = a_i \, b_i \, \in \Skalar \, , \quad \f a \times \f b = a_i \, b_j \, \epsilon_{ijk} \, \f e_k \quad \in \Euklid \,,
\end{align}
denote the inner and outer product on $\R^3$ with associated vector norm $\|{\f a}\|^2_{\R^3}=\langle {\f a},{\f a}\rangle_{\R^3}$.
The standard Euclidean inner product on $\R^{3\times 3}$ is given by
\begin{align}\label{TensorInnerProduct}
\scal{\f X , \f Y}_{\R^{3 \times 3}} = X_{ij} \, Y_{ij} \quad \in \Skalar \, ,
\end{align}
and thus the Frobenius tensor norm is $\|{\f X}\|^2=\langle{\f X},{\f X}\rangle_{\R^{3\times3}}$. Similarly, on $\R^{3 \times 3 \times 3}$ we consider the inner product
\begin{align}\label{Tensor3StufeInnerProduct}
\scal{\m , \n}_{\R^{3 \times 3 \times 3}} = \e m_{ijk} \, \e n_{ijk} \quad \in \Skalar \, ,
\end{align}
and the tensor norm $\|{\m}\|^2_{\R^{3\times3 \times 3}}=\langle{\m},{\m}\rangle_{\R^{3\times3 \times 3}}$. The identity tensor on $\R^{3\times3}$ will be denoted by $\id = \delta_{ij} \, \f e_i \otimes \f e_j$, so that $\tr({\f X})=\scal{\f X,\id}_{\R^{3\times 3}} = X_{ij} \, \delta_{ij}=X_{ii}$ and $\scal{\f X , \f Y}_{\R^{3\times 3}} =\tr({\f X \, \f Y^T})$. Further, we will make repeated use of the identity $\tr[\f X \, \f Y \, \f Z] = \scal{\f X \, \f Y \, \f Z , \id} = \scal{\f Y , \f X^T \, \f Z}$. The divergence of a vector field $\f v$ reads:
\begin{align}\label{VektDiv}
\divx \f v  \colonequals \tr ( \Gradx{\f v}) = \langle \Gradx{\f v} , \id \rangle = v_{i,j} \, \delta_{ij} =  v_{i,i} \quad \in \Skalar \, .
\end{align}
Note that we do not introduce symbolic notation for the double contraction or more complicated contractions for higher order tensors, since this may be confusing and also limited for certain cases. Instead, we use generally index notation for contractions being not defined by the scalar product. Thus, the divergence of a second order tensor is given by
\begin{align}\label{TensorDiv}
\Divx \f Y  \colonequals (\GRAD{\f Y})_{ijk} \, \delta_{jk} \, \f e_i = Y_{ij,j} \, \f e_i \quad \in \Euklid \, .
\end{align}
Similarly, the divergence of a third order tensor reads
\begin{align}\label{TensorDivThirdOrder}
\DIVx \m  \colonequals
\e m_{ijk,k} \, \f e_i \otimes \f e_j \quad \in \Skalar^{3 \times 3} \, .
\end{align}
The curl of a vector $\f v$ is given by
\begin{align}\label{CurlDefVektor}
\curlx \f v \colonequals  -  v_{a,b} \, \epsilon_{abi} \, \f e_i \quad \in \Euklid \, .
\end{align}
Similarly, the Curl  of a second order tensor field $\f X \in  \Skalar^{3 \times 3}$ is defined by
\begin{align}\label{CurlDefTensor}
\Curlx {\f X} \colonequals  - X_{ia,b} \, \epsilon_{abj} \, \f e_i \otimes \f e_j \quad \in \Skalar^{3 \times 3} \, .
\end{align}
Using eq.\eqref{axlSkewTensor}  the curl of a vector can be written in terms of a gradient via
\begin{align}\label{axlgradvcurl}
\axl ( \skew \Gradx{\f v}) & = - \12 \, (\Gradx{\f v})_{ab} \, \epsilon_{abk} \, \f e_k = \12 \,\curlx \f v \,.
\end{align}
In this work we consider a body which  occupies  a bounded open set $\Omega \subseteq \Euklid$ of the three-dimensional Euclidian space $\R^3$ and assume that its boundary $\partial \Omega$ is a piecewise smooth surface. An elastic material fills the domain $\Omega \subset \R^3$ and we refer the motion of the body to the displacement field $\ux : \Omega \subseteq \Euklid \mapsto \R^3$ shifting any point $\f x$ of the reference configuration to the actual configuration $\f x + \ux$.\\

\noindent
Since the gradient of a scalar $\phi$ is curl-free
\begin{align}\label{CurlGradSkal}
\curlx \gradx{\phi} & = \curlx [\phi_{,i} \, \f e_i] = - \phi_{,ij} \, \epsilon_{ijk} \, \f e_k \notag \\
&=  -(\phi_{,12} - \phi_{,21}) \, \f e_3 - (\phi_{,23} - \phi_{,32}) \, \f e_1 - (\phi_{,31} - \phi_{,13}) \, \f e_2 = 0 \,,
\end{align}
the definition in eq.\eqref{CurlDefTensor} is such that the gradient of a vector field $\f u$ is also curl-free:
\begin{align}\label{CurlGradVect}
\Curlx \Gradx{\f u} & =\Curlx [ u_{i,j} \, \f e_i \otimes \f e_j] = - u_{i,jb} \, \epsilon_{jbk} \, \f e_i \otimes \f e_k \notag \\
&= - (u_{i,12} - u_{i,21}) \, \f e_i \otimes \f e_3 - (u_{i,23} - u_{i,32}) \, \f e_i \otimes \f e_1 - (u_{i,31} - u_{i,13}) \, \f e_i \otimes \f e_2 = 0 \, .
\end{align}
However, the skew-symmetric part of the gradient of a vector field $\f u$ is generally {\bf not} curl-free:
\begin{align}\label{CurlskewGradVect}
\Curlx \skew \Gradx{\f u} &= \Curlx(\12 \, \epsilon_{nij} \, \epsilon_{nab}  \, u_{a,b} \, \f e_i \otimes \f e_j) =-\12 \, \epsilon_{nij} \, \epsilon_{nab}  \, u_{a,bk}  \, \epsilon_{jkp} \, \f e_i \otimes \f e_p  \notag \\
&= - \12 \, (\delta_{kn} \, \delta_{ip} - \delta_{ki} \, \delta_{pn}) \, \epsilon_{nab} \, u_{a,bk} \, \f e_i \otimes \f e_p
= - \12 (u_{a,bn} \, \epsilon_{abn} \, \delta_{ip} - u_{a,bi} \, \epsilon_{abp} ) \, \f e_i \otimes \f e_p \notag \\
& = - \12 \underbrace{(u_{1,23} - u_{1,32} + u_{2,31} - u_{2,13} + u_{3,12} - u_{3,21})}_{\displaystyle \tr \Curlx \Gradx{\f u} = 0} \, \id + \frac{\partial}{\partial \, x_i} ( \underbrace{\12 \, u_{a,b} \, \epsilon_{abp}}_{\displaystyle (- \axl \skew \Gradx{\f u})_p} ) \, \f e_i \otimes \f e_p \notag \\
&= - (\Gradx{\axl \skew \Gradx{\f u}})_{pi} \, \f e_i \otimes \f e_p \notag \\
& = -  (\Gradx{\axl \skew \Gradx{\f u}})^T = - \12 \, (\Gradx{\curlx \f u})^T  \, .
\end{align}
In the indeterminate couple stress model the gradient of the continuum rotation $\curlx \, \ux$ defines the curvature measure
\begin{align}\label{Curvature_V1}
\ksx \colonequals \12 \, \Gradx{ \curlx \ux}  = - \12 \,  u_{a,bj} \, \epsilon_{abi} \, \f e_i \otimes \f e_j \, .
\end{align}
From eq.\eqref{axlgradvcurl} and \eqref{CurlskewGradVect}  it follows that
\begin{align}\label{Curvature_V2}
\ksx  = \Gradx{\axl \skew \Gradx{\ux}} = - (\Curlx \skew \Gradx{\ux})^T \, .
\end{align}
Further, it has been shown in \cite{GhiNeMaMue15} (see also \cite[p.~27]{Sokolnikoff56}) that $\ksx$ can be also written in terms of
\begin{align}\label{Curvature_V3}
\ksx  = (\Curlx \sym \Gradx{\ux})^T \, .
\end{align}
The representations of $\ksx$ with the skew-symmetric gradient in eq.\eqref{Curvature_V2} on the one hand, and the symmetric gradient in eq.\eqref{Curvature_V3} on the other hand seems amazing. For the proof of relation \eqref{Curvature_V3} we use eq.\eqref{CurlGradVect}, \eqref{CurlskewGradVect},  and $\f X = \sym \f X + \skew \f X$ to obtain the above statement:
\begin{align}\label{CartanGradu}
\underbrace{\Curlx \Gradx{\ux}}_{\displaystyle = \, 0} & = \Curlx \{ \sym \Gradx{\ux} + \skew \Gradx{\ux} \} \notag \\
\Leftrightarrow \, 0 & = \Curlx \sym \Gradx{\ux} +  \underbrace{\Curlx \skew \Gradx{\ux}}_{\displaystyle = - (\ksx)^T} \notag \\
\Leftrightarrow \quad (\ksx)^T & =\Curlx \sym \Gradx{\ux}   \,.
\end{align}
Note that the trace of $\ksx$ is zero:
\begin{align}\label{traceCurvature}
2 \, \tr ( \ksx ) &=  \tr ( \Gradx{ \curlx \ux})  = \divx \curlx \ux = \scal{\Gradx{ \curlx \ux}, \id}  = - u_{a,bj} \, \epsilon_{abi} \, \delta_{ij}  \notag \\
&= - u_{a,bi} \, \epsilon_{abi} = - u_{1,23} + u_{1,32} - u_{2,31} + u_{2,13} - u_{3,12} + u_{3,21} = 0 \, .
\end{align}
\subsection{Identities from the Levi-Civita tensor}\label{KapLeviCivita}
In section \ref{KapNotations} we have used the baseless permutation $\epsilon_{ijk}$ to define some functions. Therefore, the base of the functions argument yields the base of the function automatically. On the one hand, this is convenient but on the other hand not sufficiently precise to investigate some transformation rules to be considered in section \ref{KapRotInvariance}. Thus, we make use of the third order Levi-Civita tensor
\begin{align}\label{LeviCivitaDef}
\fg \epsilon \colonequals \scal{\f e_i , \f e_j \times \f e_k} \, \f e_i \otimes \f e_j \otimes \f e_k = \epsilon_{ijk} \, \f e_i \otimes \f e_j \otimes \f e_k \quad \in \Skalar^{3 \times 3 \times 3} \, .
\end{align}
The components of the Levi-Civita tensor are the scalar triple product of an orthogonal unit base, which may be rotated. Therefore, let $\f Q$ be a constant rotation tensor with $\f Q^T \, \f Q=\f Q \, \f Q^T = \id$ and $\det\f Q = +1$, mapping the orthogonal referential system of Euclidean vectors $\f e_i$ to a rotated system $\f d_i$ via
\begin{align}\label{QmapsEuklid}
\f d_i = \f Q \, \f e_i \, , \qquad \f Q = \f d_i \otimes \f e_i \, , \qquad Q_{ia} = \scal{\f d_i , \f e_a}  = \scal{\f e_a , \f d_i} \, , \qquad i=1,2,3 \, , \qquad a=1,2,3 \,.
\end{align}
Note that the components $Q_{ia}$ are defined by the commutative inner product. However, the rotation tensor $\f Q$ is generally not symmetric: $\f Q \neq \f Q^T$. Since the referential basis vectors $\f e_i$ are considered to be orthogonal and of unit length the rotated basis vectors $\f d_i$ are orthogonal and of unit length as well:
\begin{align}
\scal{\f d_i , \f d_j} = \scal{Q_{ia} \, \f e_a , Q_{jb} \ \f e_b} = Q_{ia} \, Q_{jb} \, \delta_{ab}  = Q_{ia} \, Q_{ja}  = Q_{ia} \, Q^T_{aj} = \delta_{ij}\, .
\end{align}
Thus, the Levi-Civita tensor $\fg \epsilon = \epsilon^{\sharp}_{ijk} \, \f d_i \otimes \f d_j \otimes \f d_k$  can be written in terms of the rotated base with components
\begin{align}\label{LeviCivitaIdentity}
\epsilon^{\sharp}_{ijk} = \scal{\f d_i , \f d_j \times \f d_k } = \scal{\f Q_{ia} \, \f e_a , \f Q_{jb} \, \f e_b \times \f Q_{kc} \, \f e_c } = Q_{ia} \, Q_{jb} \, Q_{kc} \, \scal{\f e_i , \f e_j \times \f e_k} = Q_{ia} \, Q_{jb} \, Q_{kc} \, \epsilon_{abc}  \, .
\end{align}
Since the scalar triple product does not depend on the direction of the orthogonal unit base, the components of the Levi-Civita tensor are independent concerning the direction of the orthogonal unit base\footnote{This is similar to the identity tensor $\id \in \Skalar^{3 \times 3}$. Its components are independent concerning the direction of the orthogonal unit base.} in $\Skalar^{3 \times 3 \times 3}$. Thus, the Levi-Civita tensor is an isotropic tensor of order 3 with respect to proper orthogonal transformations, which is reading $\epsilon^{\sharp}_{ijk} = \epsilon_{ijk}$ in eq.\eqref{LeviCivitaIdentity}, and yielding the identity
\begin{align}\label{EpsilonQQQ}
\epsilon_{ijk} = Q_{ia} \, Q_{jb} \, Q_{kc} \, \epsilon_{abc}  \, .
\end{align}
Since the symbolic notation of tensor products is only established in $\Skalar^{3 \times 3}$ we prescind from defining a new kind of simple contraction. This would be necessary to switch from index to symbolic notation in eq.\eqref{EpsilonQQQ}, where each rotation tensor $\f Q$ has one simple contraction with the corresponding index of the permutation.\footnote{The simple contraction concerning the index $b$ in eq.\eqref{EpsilonQQQ} is not possible via standard symbolic notation.} Since the rotation $\f Q$ in eq.\eqref{EpsilonQQQ} is arbitrary, the identity must also hold for the transposed of $\f Q$, which is also in $\SO(3)$. Thus, a similar identity to eq.\eqref{EpsilonQQQ} is easily deduced:
\begin{align}\label{EpsilonQQQ2}
\epsilon_{ijk} = Q_{ai} \, Q_{bj} \, Q_{ck} \, \epsilon_{abc}  \, .
\end{align}
 Another identity can be found via simple contraction of $\f Q$ and the Levi-Civita tensor, which is considered in the basis $\f e_i$ on the left hand side and by the basis $\f d_i$ on the right hand side:
\begin{equation}\label{LeviCivitaQEpsQQ}
 {\begin{array}{lrcl}
  & \f Q \, \fg \epsilon & = & \f Q \, \fg \epsilon  \\
 \Leftrightarrow \, &\f Q \, (\epsilon_{ijk} \, \f e_i \otimes \f e_j \otimes \f e_k ) & = & \f Q \, (\epsilon_{ijk} \, \f d_i \otimes \f d_j \otimes \f d_k )  \vspace{0.1cm} \\
 \Leftrightarrow \, & \f Q \, (\epsilon_{ijk} \, \f e_i \otimes \f e_j \otimes \f e_k ) & = & \f Q \, (Q_{ia} \,Q_{jb} \,Q_{kc} \, \epsilon_{ijk} \, \f e_a \otimes \f e_b \otimes \f e_c )  \vspace{0.1cm} \\
\Leftrightarrow \,  & Q_{mn} \, \epsilon_{ijk} \, \f e_m \otimes \underbrace{\f e_n \cdot \f e_i}_{\displaystyle \delta_{ni}} \otimes \, \f e_j \otimes \f e_k & = & Q_{mn} \, Q_{ia} \,Q_{jb} \,Q_{kc} \, \epsilon_{ijk} \, \, \f e_m \otimes \underbrace{\f e_n \cdot  \f e_a}_{\displaystyle \delta_{na}} \otimes \, \f e_b \otimes \f e_c  \\
\Leftrightarrow \,  & Q_{mi} \, \epsilon_{ijk} \, \f e_m \otimes \f e_j \otimes \f e_k & = & \underbrace{Q_{ma} \, Q^T_{ai}}_{\displaystyle \delta_{mi}} \,Q_{jb} \,Q_{kc} \, \epsilon_{ijk} \, \, \f e_m \otimes  \f e_b \otimes \f e_c  \\
\Leftrightarrow \,  & Q_{mi} \, \epsilon_{ijk} \, \f e_m \otimes \f e_j \otimes \f e_k & = & Q_{jb} \,Q_{kc} \, \epsilon_{ijk} \, \, \f e_i \otimes  \f e_b \otimes \f e_c \vspace{0.1cm}  \\
\Leftrightarrow \,  & Q_{mi} \, \epsilon_{ibc} \, \f e_m \otimes \f e_b \otimes \f e_c & = & Q_{jb} \,Q_{kc} \, \epsilon_{mjk} \, \, \f e_m \otimes  \f e_b \otimes \f e_c \vspace{0.1cm}  \\
\Leftrightarrow \,  & Q_{mi} \, \epsilon_{ibc} & = & Q_{jb} \,Q_{kc} \, \epsilon_{jkm} \, .
 \end{array}}
\end{equation}
Note that eq.\eqref{LeviCivitaQEpsQQ} links a linear to a quadratic term of $\f Q$.
\section{Objectivity and isotropy in nonlinear elasticity}\label{KapObjNonlinearElas}
\setcounter{equation}{0}
\subsection{Objectivity and isotropy in nonlinear elasticity - the local case}\label{KapObjNonlinearElasLocal}
In geometrically nonlinear hyperelasticity we know that frame-indifference is left-invariance of the energy under SO(3)-action and isotropy is right-invariance under SO(3), i.e.\footnote{We do not discuss O(3)-invariance for simplicity.} for the energy density $W$ we have
\begin{align}
{\bf frame-indifference:} \qquad & W(\f Q_1 \, \f F) = W(\f F) \qquad \forall \, \f Q_1 \in \SO(3) \,, \label{FrameDiff} \\
{\bf isotropy:} \qquad & W(\f F \, \f Q_2) = W(\f F) \qquad \forall \, \f Q_2 \in \SO(3) \label{Isotropy} \,.
\end{align}

\noindent
Therefore, by specifying $\f Q_2 = \f Q_1^T$  we obtain also the necessary invariance condition
\begin{align}\label{wqfqt}
W\,(\f Q \, \f F \, \f Q^T) =W \, (\f F) \qquad \forall \, \f Q \in \SO(3) \,.
\end{align}
Note that condition \eqref{wqfqt} does not imply objectivity, as can be seen from considering the energy expressions
\begin{align}\label{IsoNonElas}
W(\f F) &= \norm{\f F - \id}^2 \neq \norm{\f Q \, \f F - \id}^2 = W(\f Q \, \f F)   \, , \notag \\
W(\f F) &= \norm{\sym(\f F - \id)}^2 \neq \norm{\sym(\f Q \, \f F - \id)}^2 = W(\f Q \, \f F)   \, ,
\end{align}
satisfying \eqref{wqfqt} but not being frame-indifferent \cite{FosdickSerrin79}. With the same example one sees that condition \eqref{wqfqt} does not imply isotropy. Therefore, applying the transformation $F \mapsto Q \, F \, Q^T$ (i.e. the $\sharp$-transformation defined in eq.\eqref{phisharpxi}) has {\bf no} intrinsic meaning in nonlinear elasticity theory as such.\\

\noindent
It is clear that every frame-indifferent elastic energy $W$ can be expressed in the right Cauchy-Green tensor
$\f C = \f F^T \, \f F$ in the sense that there is a function $\Psi: \Sym^+(3) \mapsto \Skalar$ such that $\forall \, \f F \in \Skalar^{3 \times 3}$
\begin{align}\label{WvonC}
W(\f F) = \Psi(\f F^T \f F) \,.
\end{align}
Of course, any $W$ of the form \eqref{WvonC} is automatically frame-indifferent. Applying the isotropy condition \eqref{Isotropy} to the representation in eq. \eqref{WvonC} we must have
\begin{align}\label{WvonFQ}
W(\f F \, \f Q) = W(\f F) \quad \Leftrightarrow \quad \Psi(\f Q^T  F^T \f F \, \f Q) = W(\f F \, \f Q) = W(\f F) = \Psi(\f F^T \f F) \qquad \forall \, \f Q \in \SO(3) \,.
\end{align}
Therefore, for reduced energy expressions \eqref{WvonC} the isotropy requirement can be equivalently stated as
\begin{align}\label{IsoRequirement}
\Psi(\f Q \, \f C \, \f Q^T) = \Psi(\f C) \qquad \forall \, \f Q \in \SO(3) \,.
\end{align}
In eq.\eqref{IsoRequirement} we see that isotropy is invariance of the function $\Psi$ under simultaneous spatial and referential rotation of the coordinate system with the same rotation $\f Q$ by interpreting the Cauchy-Green tensor $\f C$ as a linear mapping $\f C: \Euklid \rightarrow \Euklid$, which transforms under such a change of coordinate system as
\begin{align}\label{CtoQCQT}
\f C \rightarrow \f Q \, \f C \, \f Q^T \, .
\end{align}
We may now define a transformation of the deformation $\fg \varphi: \, \Omega \subseteq \Euklid \mapsto \R^3$ to new spatial and referential coordinates via\footnote{For isotropy considerations it is useful to think of the domain $\Omega \subseteq \Euklid$ to be a finite-sized ball $B \subseteq \Euklid$ which is invariant under rotations. Then $\varphi^{\sharp}: \, B \subseteq \Euklid \mapsto \Euklid$.}
\begin{align}\label{phisharpxi}
\phixi \colonequals \f Q \, \fg \varphi (\f Q^T \, \fg \xi) = \f Q \, \phix \qquad {\rm with} \quad \fg \xi = \f Q \, \f x  \, , \quad  \f x = \f Q^T \,\fg \xi \,.
\end{align}
Then\footnote{Antman \cite[p.~470]{Antman95} starts directly by considering $F^{\sharp}=Q \, F \, Q^T$ in his discussion of isotropy for frame-indifferent formulations, however, he remains on the level of stresses.}
\begin{align}\label{Gradxiphi}
\Gradxi{\phixi} = \Gradxi{ \f Q \, \fg \varphi (\f Q^T \, \fg \xi)} = \frac{\partial \, [ \f Q \, \fg \varphi (\f Q^T \, \fg \xi)]}{\partial \, \fg \xi}= \f Q \, \frac{\partial \, \fg \varphi (\f x)}{\partial \, \f x} \, \frac{\partial \, \f x}{\partial \, \fg \xi}= \f Q \, \Gradx{\fg \varphi (\f x)} \, \f Q^T
\end{align}
and
\begin{align}\label{Gradxiphi}
\Cxi = (\Gradxi{\phixi})^T \, \Gradxi{\phixi}=(\Gradx{\phix} \, \f Q^T)^T \, \Gradx{\phix} \f Q^T = \f Q \, \Cx \, \f Q^T  \,.
\end{align}
Gathering our findings so far we can state:
\begin{align}\label{WvonGradxi}
\fbox{
\parbox[][2cm][c]{14cm}{
The nonlinear hyperelastic formulation is frame-indifferent and isotropic {\bf if and only if} there exists $\Psi: \Sym^+(3) \rightarrow \Skalar$ such that $W(\f F) = \Psi(\f F^T \, \f F)$ and
$
W(\Gradxi{\phixi}) \overbrace{=}^{\rm form-invariance} W(\Gradx{\phix}) \quad \left[{\rm i.e.} \quad \Psi(\f Q \, \f C \, \f Q^T) = \Psi(\f C) \quad \forall \, \f Q \in \SO(3) \right] \, .
$
}}
\end{align}
In nonlinear elasticity, invariance of the formulation under eq.\eqref{phisharpxi} is a consequence of objectivity (left $\SO(3)$-invariance of the energy). Therefore, the transformation in eq.\eqref{phisharpxi} may be used to probe objectivity and isotropy of the formulation, while it is not equivalent to both. If we assume already the reduced representation in $\Psi$, then rotational invariance of the formulation under eq.\eqref{phisharpxi} is equivalent to isotropy, see also \cite[p.~220]{Marsden83} .\\

\noindent
The well-known representation theorems imply that any $W$ satisfying the isotropy condition and the classical format of material frame-indifference,
\begin{align}\label{WQFWSchlange}
W(QF) = W(F) = \Psi(C) \quad \forall Q \in \SO(3) \,,
\end{align}
must be expressible in terms of the principal invariants of $C=F^T \, F$, i.e.,
\begin{align}\label{WQFWSchlange2}
W(F) = \Psi(C) = \widehat{\Psi}(I_1(C), \, I_2(C), \, I_3(C)) \,,
\end{align}
where
\begin{align}\label{InvariantenvonC}
I_1(C) \colonequals \tr(C) \,, \quad I_2(C) \colonequals \tr(\Cof C) \, \quad I_3(C) \colonequals \det C \,.
\end{align}
Incidentally, $W(F) = \Psi(C) = \widehat{\Psi}(I_1(C),I_2(C),I_3(C))$ is not only invariant under compatible changes of the reference configuration with rigid rotations $\overline{Q}$ (constant rotations), but also under inhomogeneous rotation fields $Q(x)$, that is,
\begin{align}\label{WQFWSchlange3}
W(F \, Q(x)) = W(F) \qquad \forall Q(x) \in \SO(3) \,,
\end{align}
since
\begin{align}\label{InvariantenvonFQ}
I_k((F \, Q(x))^T \, F \, Q(x)) = I_k(Q(x)^T \, F^T \, F \, Q(x)) = I_k(F^T \, F) \,, \quad k=1,2,3 \,.
\end{align}
Therefore, the principal invariants $I_k$ are unaffected by inhomogeneous rotations. However, we must stress that isotropy per se is not defined as invariance under right multiplication with $Q = Q(x) \in \SO(3)$. The difference between form-invariance under compatible transformations with rigid rotations $\overline{Q}$ (isotropy) and right-invariance under inhomogeneous rotation fields $Q = Q(x) \in \SO(3)$ will only become visible in higher gradient elasticity treated in the next section. In any first gradient theory both requirements coincide. Note that likewise, \eqref{WQFWSchlange3} is invariant under left multiplication with inhomogeneous rotation fields. We define for further use
 \begin{align}\label{leftrightSO3invarianceLocal}
\fbox{
\parbox[][1.0cm][c]{14cm}{
$
\text{\bf left-local SO(3)-invariance} \,\,\,\,\, : \quad W(Q(x) \, F(x)) = W(F(x)) \quad \forall Q(x) \in \SO(3) \\
\\
\text{\bf right-local SO(3)-invariance} \, : \quad W(F(x) \, Q(x)) = W(F(x)) \quad \forall Q(x) \in \SO(3)  \,.
$}}
 \end{align}

Similarly, we define for {\bf constant} rotations $\overline{Q} \in \SO(3)$
 \begin{align}\label{leftrightSO3invarianceGlobal}
\fbox{
\parbox[][1.0cm][c]{14cm}{
$
\text{\bf left-global SO(3)-invariance} \,\,\,\,\, : \quad W(\overline{Q} \, F(x)) = W(F(x)) \quad \forall \overline{Q} \in \SO(3) \quad \text{\bf objectivity} \\
\\
\text{\bf right-global SO(3)-invariance} \, : \quad W(F(x) \, \overline{Q}) = W(F(x)) \quad \forall \overline{Q} \in \SO(3)   \quad \text{\bf isotropy} \,.
$}}
 \end{align}

\noindent
With these definitions we have shown that for $W(F) = \widehat{\Psi}(I_1(C), \, I_2(C), \, I_3(C))$ both, local and global, left and right SO(3) invariance are satisfied. Left-global SO(3)-invariance is identical to the Cosserat's invariance under "action euclidienne" \cite{Cosserat09}.
\subsection{Isotropy in second gradient nonlinear elasticity}\label{KapSecondGradElas}
In this subsection we would like to extend the "local" picture of the previous subsection to second gradient materials.
Using the Noether theorem and Lie-point symmetries, the explicit expressions of the isotropy condition in linear gradient elasticity of grade-2 and in linear gradient elasticity of grade-3 have been derived by Lazar and Anastassiadis
\cite[eq.(4.76)]{LazarAnastassiadis07} and Agiasofitou and Lazar \cite[eqs.(57) and (58)]{Agiasofitou2009}, respectively, as consequence of global rotational invariance.\\

\begin{figure}[h]
\centering
 \begin{picture}(130,110)
  \put(0,0){\epsfig{file=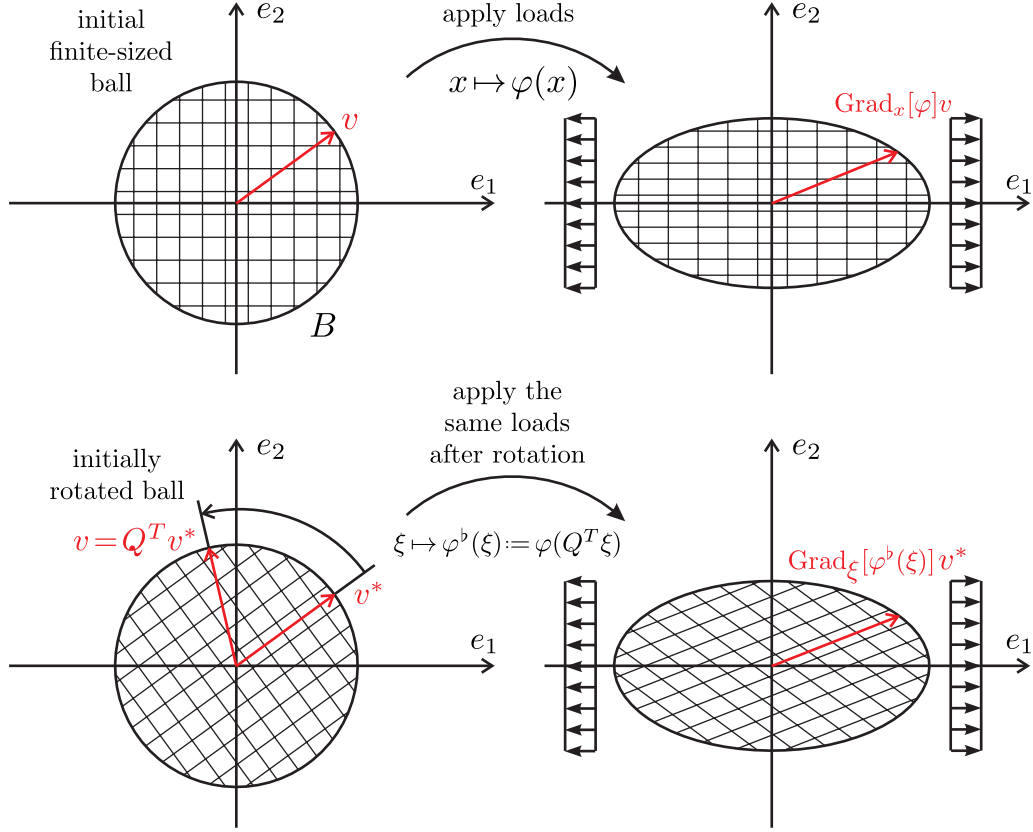, height=110mm, angle=0}}
 \end{picture}
\caption{Isotropy: We consider a finite-sized ball $B$ of homogeneous material in an initial undistorted configuration, which may be rotated by a constant rotation $Q \in SO(3)$ (left). Isotropy requires that the deformation response is independent of $Q$ if we apply the same loading to both configurations: there is no preferred direction (right).}
 \label{ObjectivityPic}
\end{figure}

For the sake of clarity we leave first objectivity aside and discuss only isotropy. Moreover, we make the following simplifying assumptions. We consider a {\bf homogeneous material} given as a {finite-sized ball} $B \subset \Euklid$ and we want to formalize the statement that first rotating the ball $B$ and then to apply the loads leaves invariant the response in the sense that this is indistinguishable from not rotating the ball, see Fig. \ref{ObjectivityPic}. This is in accordance with the statement in Truesdell \& Noll \cite[p.~78]{Trues65}: "For an isotropic material in an undistorted state, a physical test cannot detect whether or not the material has been rotated arbitrarily before the test ist made."

Let  the elastic energy of the body $B \subset \Euklid$ depend also on second gradients of the deformation, i.e. we consider therefore,
\begin{align}\label{IntWGradGrad}
\int \limits_{\displaystyle x \in B} \!\!\!\! W \left( F(x) , \, \GRADx{F(x)} \right) \, \di x =\int \limits_{\displaystyle x \in B} \!\!\!\! W \left( \Gradx{\varphi(x)} , \, \GRADx{\Gradx{\varphi(x)}} \right) \, \di x \,.
\end{align}
We first transform eq.\eqref{IntWGradGrad} to new coordinates $\xi \in \Euklid$ via introducing an orientation preserving diffeomorphism $\zeta: \, B \subset \Euklid \mapsto \Euklid$,
\begin{align}\label{ZetaXix}
x &= \zeta(\xi) \, , \qquad \xi = \zeta^{-1}(x) \, , \qquad x=\zeta(\zeta^{-1}(x)) \, ,
\end{align}
see also \cite{Neff_techmech07,Neff_Svendsen08}.
For further reference, we compute some relations connected to the transformation \eqref{ZetaXix}. Taking the first and second derivative with respect to $x$ in \eqref{ZetaXix}$_3$ and using \eqref{ZetaXix}$_2$ we obtain
\begin{align}\label{TranszetaGrad}
 \id= \Gradxi{\zeta(\zeta^{-1}(x))} \, \Gradx{\zeta^{-1}(x)} \quad \Leftrightarrow \quad \Gradx{\zeta^{-1}(x)} = (\Gradxi{\zeta(\xi)})^{-1} \, ,
\end{align}
and
\begin{align}\label{TranszetaGrad2}
& 0= \GRADx{\Gradxi{\zeta(\zeta^{-1}(x))} \, \Gradx{\zeta^{-1}(x)}} \notag \\
& 0= \GRADx{\Gradxi{\zeta(\zeta^{-1}(x))}} \, \Gradx{\zeta^{-1}(x)} + \Gradxi{\zeta(\zeta^{-1}(x))} \, \GRADx{\Gradx{\zeta^{-1}(x)}} \notag \\
\Leftrightarrow \,\, &  0= \GRADxi{\Gradxi{\zeta(\zeta^{-1}(x))}} \, \Gradx{\zeta^{-1}(x)} \, \Gradx{\zeta^{-1}(x)} + \Gradxi{\zeta(\zeta^{-1}(x))} \, \GRADx{\Gradx{\zeta^{-1}(x)}} \notag \\
\Leftrightarrow \,\, & \Gradxi{ \zeta(\zeta^{-1}(x))} \, \GRADx{\Gradx{\zeta^{-1}(x)}} = - \GRADxi{\Gradxi{\zeta(\xi)}} \, \Gradx{\zeta^{-1}(x)} \, \Gradx{\zeta^{-1}(x)} \,.
\end{align}
With help of eq.\eqref{TranszetaGrad} the latter implies
\begin{align}\label{TranszetaGrad3}
 \GRADx{\Gradx{\zeta^{-1}(x)}} = -(\Gradxi{ \zeta(\xi)})^{-1} \, \GRADxi{\Gradxi{\zeta(\xi)}} \, (\Gradxi{\zeta(\xi)})^{-1} \, (\Gradxi{\zeta(\xi)})^{-1} \,.
\end{align}
Connected to the coordinate transformation \eqref{ZetaXix} we consider the deformation expressed in these new coordinates via setting
\begin{align}\label{PhibZeta}
\varphi^{\flat}(\xi) \colonequals \varphi(\zeta(\xi)) \, , \qquad  \varphi^{\flat}(\zeta^{-1}(x)) = \varphi(x) \, .
\end{align}
The standard chain rule induces for the first and second derivative of $\varphi(x)$ with respect to $x$ the following relations
\begin{align}\label{GradNewVarphi1}
\Gradx{\varphi(x)} = \Gradxi{\varphi^{\flat}(\zeta^{-1}(x))} \, \Gradx{\zeta^{-1}(x)}
\, ,
\end{align}
and
\begin{align}\label{GradNewVarphi2}
& \GRADx{\Gradx{\varphi(x)}} \notag \\
 & \qquad = \GRADxi{\Gradxi{\varphi^{\flat}(\xi)}} \, \Gradx{\zeta^{-1}(x)} \, \Gradx{\zeta^{-1}(x)} + \Gradxi{\varphi^{\flat}(\xi)} \, \GRADx{\Gradx{\zeta^{-1}(x)}} \,.
\end{align}
The elastic energy \eqref{IntWGradGrad}, expressed in new coordinates via the transformation of integrals formula reads then, inserting \eqref{GradNewVarphi1} and \eqref{GradNewVarphi2}
\begin{align}\label{IntWGradGrad2}
\int \limits_{\displaystyle \xi \in \zeta^{-1}(B)} & \!\!\!\!\!\!\!\!\! W \bigg(  \Gradxi{\varphi^{\flat}(\xi)} \, \Gradx{\zeta^{-1}(x)} , \, \GRADxi{\Gradxi{\varphi^{\flat}(\xi)}} \, \Gradx{\zeta^{-1}(x)} \, \Gradx{\zeta^{-1}(x)} \notag \\
&\qquad \qquad \qquad \qquad + \Gradxi{\varphi^{\flat}(\xi)} \, \GRADx{\Gradx{\zeta^{-1}(x)}} \bigg) \, \det(\Gradxi{\zeta(\xi)}) \, \di \xi \,.
\end{align}
We rewrite \eqref{IntWGradGrad2} as follows, using again $\Gradx{\zeta^{-1}(x)} = (\Gradxi{\zeta(\xi)})^{-1}$ from \eqref{TranszetaGrad}
\begin{align}\label{IntWGradGrad3}
\int \limits_{\displaystyle \xi \in \zeta^{-1}(B)}\!\!\!\!\!\!\!\! W \bigg(  \Gradxi{\varphi^{\flat}(\xi)} \,  & (\Gradxi{\zeta(\xi)})^{-1} , \,
 \GRADxi{\Gradxi{\varphi^{\flat}(\xi)}} \, (\Gradxi{\zeta(\xi)})^{-1} \, (\Gradxi{\zeta(\xi)})^{-1} \notag \\
&\qquad \qquad \qquad  + \Gradxi{\varphi^{\flat}(\xi)} \, \GRADx{\Gradx{\zeta^{-1}(x)}} \bigg) \, \det(\Gradxi{\zeta(\xi)}) \, \di \xi \,.
\end{align}
Further, we use relation \eqref{TranszetaGrad3} to substitute $\GRADx{\Gradx{\zeta^{-1}(x)}}$ yielding
\begin{align}\label{IntWGradGrad4}
   & \int \limits_{\displaystyle \xi \in \zeta^{-1}(B)} \!\!\!\!\!\!\!\! W\bigg(  \Gradxi{\varphi^{\flat}(\xi)}  \, (\Gradxi{\zeta(\xi)})^{-1} , \,
 \GRADxi{\Gradxi{\varphi^{\flat}(\xi)}} \, (\Gradxi{\zeta(\xi)})^{-1} \, (\Gradxi{\zeta(\xi)})^{-1} \notag \\
 \quad
& - \Gradxi{\varphi^{\flat}(\xi)} \,  (\Gradxi{ \zeta(\xi)})^{-1}\, \GRADxi{\Gradxi{\zeta(\xi)}} \, (\Gradxi{\zeta(\xi)})^{-1} \, (\Gradxi{\zeta(\xi)})^{-1} \bigg) \, \det(\Gradxi{\zeta(\xi)}) \, \di \xi \,.
\end{align}
We say that the elastic energy \eqref{IntWGradGrad} is {\bf form-invariant} with respect to the (referential) coordinate transformation $\zeta$ if and only if
\begin{align}\label{IntWGradGrad5a}
   & \int \limits_{\displaystyle \xi \in \zeta^{-1}(B)} \!\!\!\!\!\!\!\! W\bigg(  \Gradxi{\varphi^{\flat}(\xi)}  \, (\Gradxi{\zeta(\xi)})^{-1} , \,
 \GRADxi{\Gradxi{\varphi^{\flat}(\xi)}} \, (\Gradxi{\zeta(\xi)})^{-1} \, (\Gradxi{\zeta(\xi)})^{-1} \notag \\
 \quad
& - \Gradxi{\varphi^{\flat}(\xi)} \,  (\Gradxi{ \zeta(\xi)})^{-1}\, \GRADxi{\Gradxi{\zeta(\xi)}} \, (\Gradxi{\zeta(\xi)})^{-1} \, (\Gradxi{\zeta(\xi)})^{-1} \bigg) \, \det(\Gradxi{\zeta(\xi)}) \, \di \xi \notag \\
& \qquad \qquad \qquad \qquad \qquad \qquad  =   \int \limits_{\displaystyle \xi \in \zeta^{-1}(B)} \!\!\!\!\!\!\!\! W \bigg(  \Gradxi{\varphi^{\flat}(\xi)} , \,
 \GRADxi{\Gradxi{\varphi^{\flat}(\xi)}} \bigg) \, \di \xi \,,
\end{align}
and we call then $\zeta$ a {\bf material symmetry transformation}, cf. \cite[eq.(11)]{ElzanowskiEpstein92}, \cite{Epstein2007, EpsteinLeon98} and \cite[p.~75, Def. 2.4]{Bertram2016}.  It is simple to see that for equality \eqref{IntWGradGrad5a} to be satisfied at all, we must have
\begin{align}\label{detGradxi1}
\det(\Gradxi{\zeta(\xi)}) = 1 \,.
\end{align}
It can be shown as well that we may restrict attention to the situation
\begin{align}\label{GradXizeta}
\Gradxi{\zeta(\xi)} \in {\mathcal G} \subseteq \SO(3) \qquad \forall \, \xi \in \zeta^{-1}(B) \,,
\end{align}
which is consistent with the local statement in the previous subsection \ref{KapEnerConsider}. Here, ${\mathcal G} \subseteq \SO(3)$ is a subgroup of the proper rotation group $\SO(3)$.\footnote{If ${\mathcal G} \subseteq {\rm O}(3)$ and $- \id \in {\mathcal G}$ then the material is called centro-symmetric. We do not discuss this possibility.}
The group ${\mathcal G}$ can be discrete\footnote{The first classification of discrete subgroups ${\mathcal G} \subseteq \SO(3)$ has been given in Hessel \cite{Hessel1831}.} (for e.g. orthotropy) or continuous (continuous for transverse isotropy and isotropy). We may therefore specify \eqref{IntWGradGrad5a} to the case of ${\mathcal G}\subseteq \SO(3)$ and writing
\begin{align}\label{GradxizetaQ}
(\Gradxi{\zeta(\xi)})^{-1} = Q(\xi) \,,
\end{align}
we obtain as first concise {\rm form-invariance} statement for material symmetry
\begin{align}\label{IntWGradGrad5}
  \int \limits_{\displaystyle \xi \in \zeta^{-1}(B)} \!\!\!\!\!\!\!\! W \bigg(  &  \Gradxi{\varphi^{\flat}(\xi)} \, Q(\xi) , \,
 \GRADxi{\Gradxi{\varphi^{\flat}(\xi)}} \, Q(\xi) \, Q(\xi) \notag \\
   & \qquad  - \Gradxi{\varphi^{\flat}(\xi)} \, Q(\xi) \, \GRADxi{Q^T(\xi)} \, Q(\xi) \, Q(\xi) \bigg) \, 1 \, \di \xi \notag \\
&  \qquad \qquad \qquad \qquad  =   \int \limits_{\displaystyle \xi \in \zeta^{-1}(B)} \!\!\!\!\!\!\!\! W \bigg(  \Gradxi{\varphi^{\flat}(\xi)} , \, \GRADxi{\Gradxi{\varphi^{\flat}(\xi)}} \bigg) \, \di \xi \quad \forall Q(\xi)\in {\mathcal G} \,.
\end{align}
The format \eqref{IntWGradGrad5} can be found e.g. in \cite{Epstein96, Murdoch77}. The homogeneous material is usually said to be {\bf isotropic} following the local reasoning,\footnote{In Bertram \cite[p.~75]{Bertram2016} we read "This group is used to define isotropy or anisotropy. If the symmetry group is a subgroup of [SO(3)] in the first entry and the zero in the second, $(Q, 0)$, these transformations can be interpreted as rigid rotations, and we call the respective reference placement an undistorted state. If a material allows for such undistorted states, it is a solid. If it contains all orthogonal dyadics in the first entry, then the material is called isotropic". This statement leaves open that for isotropy, we may choose inhomogeneous rotations $Q(x)$. In our setting we understand by isotropy using necessarily rigid rotations and therefore only $(Q,0)$.} whenever ${\mathcal G} \equiv \SO(3)$ and \eqref{IntWGradGrad5} holds for all non-constant rotation fields $Q(\xi) \in \SO(3)$. In the following, consistent with previous definitions, we will call \eqref{IntWGradGrad5} {\bf right-local SO(3)-invariance}.\\

\noindent
Let us proceed by showing that right-local SO(3)-invariance in \eqref{IntWGradGrad5} is, however, {\bf misconceived} as a general condition for isotropy. The misconception in \eqref{IntWGradGrad5} consists in suggesting a generality of transformation behaviour  followed for the rotation field $Q(\xi)$ that is, in reality, not present. To see this we refer to Figure \ref{ObjectivityPic}. Requiring in \eqref{GradxizetaQ} that
\begin{align}\label{GradxizetaQSO3}
(\Gradxi{\zeta(\xi)})^{-1} = Q^T(\xi) \in \SO(3) \quad \Leftrightarrow \quad \Gradxi{\zeta(\xi)} = Q(\xi) \, \in \SO(3)
\end{align}
means, by a {\bf classical geometric rigidity} result (see e.g. \cite{Neff_curl06}) that indeed
 \begin{align}\label{GradxiGradxizetaQ}
\fbox{
\parbox[][1.0cm][c]{14cm}{
$
\Gradxi{\zeta(\xi)} = Q(\xi) \, \in \SO(3) \quad \Rightarrow \quad Q(\xi) = \oQ={\rm const.} \quad {\rm and} \quad \zeta(\xi)= \oQ \, \xi + \overline{b} \\
\\
\text{and therefore} \quad \GRADxi{\Gradxi{\zeta(\xi)}} = 0 \qquad \text{\bf compatibility of rotations} \,.
$}}
 \end{align}

\noindent
Thus, the correct statement for isotropy, in our view is,\footnote{Consistent with condition \eqref{IntWGradGrad6} Huang \cite[eq.(3),(4)]{Huang67} arrives at the isotropy requirement
\begin{align}\notag
&W((F \, \oQ^T)^T \, F \, \oQ^T \!\! - \id, \, \Curl[(F \, \oQ^T)^T \, F \, \oQ^T \!\! - \id])
=W(\oQ \, (F^T F - \id) \, \oQ^T , \, \Curl[\oQ \, (F^T F - \id) \, \oQ^T]) \notag \\
&=W(\oQ \, (F^T \, F - \id) \, \oQ^T , \, \oQ \, \Curl[F^T \, F - \id] \, \oQ^T) = W(F^T \, F - \id , \, \Curl[F^T \, F - \id]) \, , \notag
\end{align}
 where $W$ is the free-energy of the nonlinear couple stress model of Toupin \cite{Toupin64}. Additionally, in \cite{Huang68} also the anisotropic case is given, c.f. \cite{Eremeyev2006, Eremeyev2012, Eremeyev2016}.} noting that for a global rotation of the coordinates we have $\zeta^{-1}(B)=B$
\begin{align}\label{IntWGradGrad6}
  \int \limits_{\displaystyle \xi \in B} \!\!\!\! W \bigg(  \Gradxi{\varphi^{\flat}(\xi)} & \, \overline{Q}^T , \,
 \GRADxi{\Gradxi{\varphi^{\flat}(\xi)}} \, \overline{Q}^T \, \overline{Q}^T \bigg) \, \di \xi \notag \\
   &=   \int \limits_{\displaystyle \xi \in B} \!\!\!\! W \bigg(  \Gradxi{\varphi^{\flat}(\xi)} , \, \GRADxi{\Gradxi{\varphi^{\flat}(\xi)}} \bigg) \, \di \xi \quad \forall \overline{Q} \in \SO(3) \,.
\end{align}
We will refer to this condition as {\bf right-global SO(3)-invariance}, which, for us, is {\bf isotropy}.\footnote{Toupin \cite[p.~108]{Toupin64} notes: "The material symmetry group of a homogeneous, isotropic material is the set of all distance-preserving transformations $[ \overline{Q} \, x +  \overline{b}]$."} We appreciate that the right-local SO(3)-invariance condition \eqref{IntWGradGrad5} is much to restrictive in that arbitrary, inhomogeneous rotation fields are allowed instead of only constant rotations $\overline{Q}$. The reader should carefully note that we started by using a {\rm coordinate transformation} $x=\zeta(\xi)$ and therefore we require in the end that $\zeta(\xi)=\overline{Q} \, \xi + \bar{b}$, in line with our understanding to really rotate the finite sized ball. There is no other coordinate transformation $\zeta$ such that $\Gradxi{\zeta(\xi)}=Q(\xi) \, \in \SO(3)$ everywhere, provided a minimum level of smoothness is assumed.
We note that \eqref{IntWGradGrad5} $\Rightarrow$ \eqref{IntWGradGrad6}.\footnote{More can be said in case of splitted energies $W(F, \, \GRAD{F}) = W_{\rm local}(F) + W_{\rm curv}(\GRAD{F})$. There, $W_{\rm local}$ must be isotropic and the condition on $W_{\rm curv}$ can be localized, yielding in the end
\begin{align}\notag
\int \limits_{\displaystyle \xi \in B} \!\!\!\! W_{\rm curv} (\GRADxi{\Gradxi{\varphi^{\flat}(\xi)}} \, Q^T\!\!(\xi) \, Q^T\!\!(\xi)) \, \di \xi = \int \limits_{\displaystyle \xi \in B} \!\!\!\! W_{\rm curv} (\GRADxi{\Gradxi{\varphi^{\flat}(\xi)}}) \, \di x \quad \forall Q(\xi) \in \SO(3) \, .
\end{align}
Thus, \eqref{IntWGradGrad6} can be expressed with inhomogeneous rotation fields since the space derivatives do not act on the rotations in this representation.}
 \begin{align}\notag
\fbox{
\parbox[][1.0cm][c]{14cm}{
$
\qquad \text {\bf right-local SO(3)-invariance}  \hspace{0.2cm} \Rightarrow \hspace{0.2cm} \text {\bf right-global SO(3)-invariance} \Leftrightarrow \text{\bf isotropy}\\
 \\
\text{ } \hspace{2.0cm} \eqref{IntWGradGrad5} \hspace{6cm} \eqref{IntWGradGrad6}
$}}
 \end{align}

\noindent
The condition \eqref{IntWGradGrad6} is consistent with \cite[eq.16]{Murdoch77} in that the term with additional first $\xi$-derivatives on $Q$ in \eqref{IntWGradGrad5} is absent. We remark again, that in the local theory the condition
\begin{align}\label{IntWGrad}
  \int \limits_{\displaystyle \xi \in B} \!\!\!\! W \left( \Gradxi{\varphi^{\flat}(\xi)} \, \overline{Q}^T \right) \di \xi =   \int \limits_{\displaystyle \xi \in B} \!\!\!\! W \left( \Gradxi{\varphi^{\flat}(\xi)} \right) \, \di \xi \quad \forall \overline{Q} \in \SO(3)
\end{align}
cannot distinguish between constant or non-constant rotations, as seen in \eqref{InvariantenvonFQ}. This means for the local theory that \eqref{IntWGrad} implies and is therefore equivalent to
\begin{align}\label{IntWGrad2}
  \int \limits_{\displaystyle \xi \in B} \!\!\!\! W \left( \Gradxi{\varphi^{\flat}(\xi)} \, \overline{Q}^T\!\!(\xi) \right) \di \xi =   \int \limits_{\displaystyle \xi \in B} \!\!\!\! W \left( \Gradxi{\varphi^{\flat}(\xi)} \right) \, \di \xi \, .
\end{align}
The latter might explain why one may be inclined to allow non-constant rotation fields in \eqref{IntWGradGrad5}, which is forbidden for higher gradient materials, as we saw.
\subsection{Objectivity and isotropy in nonlinear elasticity - the non-local case}\label{KapOINonlin-Non-Local-Case}
In higher gradient elasticity, however, not every (left-global SO(3)-invariant) objective energy is also left-local SO(3)-invariant.\footnote{Any $W(F , \GRAD{F}) = \Psi(C, \GRAD{C})$ is automatically already left-local SO(3)-invariant.} This can be seen by looking at
\begin{align}\label{WFPsiCPlus}
 W(F,\GRAD{F}) = \Psi(C) + \sum_{i=1}^3 \norm{F^T \, \frac{\partial F}{\partial x_i}}^2 \, .
\end{align}
The curvature expression $\norm{F^T \, \frac{\partial F}{\partial x_i}}^2$ is objective, but not left-local SO(3)-invariant, i.e.
\begin{align}\label{QFGradxQF}
 (Q (x)\, F)^T \, \frac{\partial (Q(x) \, F)}{\partial x_i} = (Q(x) \, F)^T \, \frac{\partial Q(x)}{\partial x_i} \, F + F^T \, Q(x)^T \, Q(x) \, \frac{\partial F}{\partial x_i} = \underbrace{F^T \, Q(x)^T \, \frac{\partial Q(x)}{\partial x_i} \, F}_{\displaystyle \in \so(3)} + F^T  \, \frac{\partial F}{\partial x_i}   \, ,
\end{align}
implying that
\begin{align}\label{notleftlocalSO3invariant}
\sum_{i=1}^3 \norm{(Q (x)\, F)^T \, \frac{\partial (Q(x) \, F)}{\partial x_i}}^2 = \sum_{i=1}^3 \norm{ F^T \, Q(x)^T \, \frac{\partial Q(x)}{\partial x_i} \, F + F^T  \, \frac{\partial F}{\partial x_i}}^2 \neq \sum_{i=1}^3 \norm{F^T \, \frac{\partial F}{\partial x_i}}^2 \, .
\end{align}
In order to obtain directly a left-local SO(3)-invariant expression, we could choose
\begin{align}\label{WSchlangeFPsiCPlus}
 \widetilde{W}(F) = \Psi(C) + \widetilde{\Psi}_{\rm curv}(\GRAD{C}) = \Psi(C) + \sum_{i=1}^3 \norm{\sym \left(F^T \, \frac{\partial F}{\partial x_i}\right)}^2 \, ,
\end{align}
since $\widetilde{\Psi}_{\rm curv}(\Grad{C})$ is depending on the full gradient of $F$ via
\begin{align}\label{GradCusingSym}
 \GRADx{C} &= \left[ \frac{\partial (F^T F)}{\partial x_1} \, , \, \frac{\partial (F^T F)}{\partial x_2} \, , \, \frac{\partial (F^T F)}{\partial x_3} \, \right] \notag \\
 &= \left[ \frac{\partial (F^T)}{\partial x_1} \, F + F^T \, \frac{\partial (F)}{\partial x_1} \, , \, \frac{\partial (F^T)}{\partial x_2} \, F+ F^T \, \frac{\partial (F)}{\partial x_2} \, , \, \frac{\partial (F^T)}{\partial x_3} \, F+ F^T \, \frac{\partial (F)}{\partial x_3} \, \right] \, .
\end{align}
Thus, using
\begin{align}\label{DFxiFT}
\left( \frac{\partial F^T}{\partial x_i} \, F \right)^T \!\!\! = F^T \! \left( \frac{\partial F^T}{\partial x_i} \right)^T = F^T \! \left( \frac{\partial F}{\partial x_i} \right)\, , \quad \text{for} \quad  i=1,2,3 \, ,
\end{align}
column-wise we get that
\begin{align}\label{PsicurvSchlange}
\widetilde{\Psi}_{\rm curv}(\GRAD{C}) &=  \norm{ \12 \, \left[ \left( F^T \, \frac{\partial (F)}{\partial x_1}\right)^T \!\!\! + F^T \, \frac{\partial (F)}{\partial x_1} \, , \, \left( F^T \, \frac{\partial (F)}{\partial x_2}\right)^T \!\!\! + F^T \, \frac{\partial (F)}{\partial x_2} \, , \, \left( F^T \, \frac{\partial (F)}{\partial x_3}\right)^T \!\!\! + F^T \, \frac{\partial (F)}{\partial x_3} \right] }^2 \notag \\
&=  \norm{ \12 \, \left[ 2 \, \sym \left( F^T \, \frac{\partial (F)}{\partial x_1}\right) \, , \, 2 \, \sym \left( F^T \, \frac{\partial (F)}{\partial x_2}\right) \, , \, 2 \, \sym \left( F^T \, \frac{\partial (F)}{\partial x_3}\right) \right]}^2 \notag\\
&=  \sum_{i=1}^3 \norm{ \sym \left( F^T \, \frac{\partial (F)}{\partial x_i}\right)}^2 = \widehat{\Psi}_{\rm curv}(F^T \, \GRADx{F}) \, .
\end{align}
Note that \eqref{WFPsiCPlus} is right-global SO(3)-invariant, but not right-local SO(3)-invariant, which can be seen similarly to the above discussion.\\

Next, we assume for simplicity a free energy density of the form $W=W(F) + W_{\rm curv}(\GRADx{F})$. Moreover, we already consider a further reduced objective format $W=\Psi(F^T \, F) + \Psi_{\rm curv}(\GRADx{F^T \, F})$, which is left-local and left-global SO(3) invariant as seen above. Now, for investigation of isotropy, it suffices to concentrate on the induced conditions on $\Psi_{\rm curv}(\GRADx{F^T \, F})$. We repeat the (referential) coordinate transformation and we already observed in \eqref{PsicurvSchlange} that we can always write
\begin{align}\label{PsiGradFTF}
 \Psi_{\rm curv}(\GRADx{F^T \, F})  = \widehat{\Psi}_{\rm curv}(F^T \, \GRADx{F})
 \, .
\end{align}
Moreover, from the results of the previous section right-global SO(3)-invariance under {\bf constant rotations} $\overline{Q} \in \SO(3)$ is equivalent to
\begin{align}\label{PsiCurv1}
\int \limits_{\displaystyle \xi \in B} \!\!\!\! \widehat{\Psi}_{\rm curv} \big( (\Gradxi{\varphi^{\flat}(\xi)} \,  & \overbrace{\Gradx{\zeta^{-1}(x)}}^{\overline{Q}})^T \, \GRADxi{\Gradxi{\varphi^{\flat}(\xi)}} \, \overline{Q} \, \, \overline{Q} \, \big) \, \di \xi \notag \\
& \qquad \qquad=   \int \limits_{\displaystyle \xi \in B} \!\!\!\! \widehat{\Psi}_{\rm curv} \big(  (\Gradxi{\varphi^{\flat}(\xi)})^T \, \GRADxi{\Gradxi{\varphi^{\flat}(\xi)}} \big) \, \di \xi \,.
\end{align}
And this is equivalent to
\begin{align}\label{PsiCurv2}
\int \limits_{\displaystyle \xi \in B} \!\!\!\! \widehat{\Psi}_{\rm curv} \big( \oQ^T \, & (\Gradxi{\varphi^{\flat}(\xi)})^T \, \GRADxi{\Gradxi{\varphi^{\flat}(\xi)}} \, \overline{Q} \, \, \overline{Q} \, \big) \, \di \xi \notag \\
&\qquad =   \int \limits_{\displaystyle \xi \in B} \!\!\!\! \widehat{\Psi}_{\rm curv} \big(  (\Gradxi{\varphi^{\flat}(\xi)})^T \, \GRADxi{\Gradxi{\varphi^{\flat}(\xi)}} \big) \, \di \xi \quad \forall \overline{Q} \in \SO(3) \,.
\end{align}
As an example, let us quickly check that the objective (left-global SO(3)-invariant) curvature expression
\begin{align}\label{NormGradGradPhi}
\norm{\GRADx{\Gradx{\varphi(x)}}}^2 \,,
\end{align}
is isotropic in the sense of right-global SO(3)-invariance \eqref{IntWGradGrad6}. Indeed, we have to check only that
\begin{align}\label{NormGradGradPhi2}
\norm{\GRADxi{\Gradxi{\varphi^{\flat}(\xi)}} \oQ \, \, \oQ }^2 = \norm{\GRADxi{\Gradxi{\varphi^{\flat}(\xi)}}}^2 \,.
\end{align}
The latter is verified since in indices we get
\begin{align}\label{NormGradGradPhi3}
 & \norm{\GRADxi{\Gradxi{\varphi^{\flat}(\xi)}} \oQ \, \oQ }^2 = \norm{(\GRADxi{\Gradxi{\varphi^{\flat}(\xi)}})_{ija} \oQ_{an} \, \oQ_{jk} \, e_i \otimes e_k \otimes e_n }^2 \notag\\
& = \scal{(\GRADxi{\Gradxi{\varphi^{\flat}(\xi)}})_{ija} \oQ_{an} \, \oQ_{jk} \, e_i \otimes e_k \otimes e_n \, , \, (\GRADxi{\Gradxi{\varphi^{\flat}(\xi)}})_{opq} \oQ_{qs} \, \oQ_{pt} \, e_o \otimes e_t \otimes e_s} \notag \\
& = (\GRADxi{\Gradxi{\varphi^{\flat}(\xi)}})_{ija} \, (\GRADxi{\Gradxi{\varphi^{\flat}(\xi)}})_{opq} \, \delta_{io} \, \oQ_{an} \, \delta_{ns} \, \oQ_{jk} \, \delta_{kt} \, \oQ_{qs} \, \oQ_{pt} \notag \\
& = (\GRADxi{\Gradxi{\varphi^{\flat}(\xi)}})_{ija} \, (\GRADxi{\Gradxi{\varphi^{\flat}(\xi)}})_{ipq} \, \oQ_{as} \, \oQ^T_{sq} \, \oQ_{jt} \,  \oQ^T_{tp} \notag \\
& = (\GRADxi{\Gradxi{\varphi^{\flat}(\xi)}})_{ija} \, (\GRADxi{\Gradxi{\varphi^{\flat}(\xi)}})_{ipq} \, \delta_{aq} \, \delta_{jp} \notag \\
&=   (\GRADxi{\Gradxi{\varphi^{\flat}(\xi)}})_{ija} \, (\GRADxi{\Gradxi{\varphi^{\flat}(\xi)}})_{ija} =  \norm{\GRADxi{\Gradxi{\varphi^{\flat}(\xi)}}}^2  \,.
\end{align}
Note, however, that \eqref{NormGradGradPhi} is a prototype expression which does not satisfy the right-local SO(3)-invariance condition \eqref{IntWGradGrad5}. This is the case, since for inhomogeneous rotation fields $Q(\xi) \in \SO(3)$ we have
\begin{align}\label{NormGradGradPhi4}
\norm{\GRADxi{\Gradxi{\varphi^{\flat}(\xi)}} Q(\xi) \, \, Q(\xi) - \Gradxi{\varphi^{\flat}(\xi)} \, Q(\xi) \, \GRADxi{Q(\xi)} & \, Q(\xi) \, \, Q(\xi) }^2 \notag\\
& \neq \norm{\GRADxi{\Gradxi{\varphi^{\flat}(\xi)}}}^2 \,.
\end{align}
It is instructive to see that the right-local SO(3)-invariance \eqref{IntWGradGrad5} can be satisfied nevertheless within a second gradient model by taking
\begin{align}\label{WPsiSchlangeGrad}
W(F, \, \GRADx{F}) = \Psi(C, \, \GRADx{C})
 = \widehat{\Psi}(I_1(C), \, I_2(C), \, I_3(C)) + \sum_{k=1}^3 \norm{\gradx{I_k(C)}}^2 \,,
\end{align}
with $\widehat{\Psi}: \Skalar^3 \mapsto \Skalar$ arbitrary. In this formulation, objectivity is clear, as the functional dependence is expressed only in $C=F^T \, F$ and isotropy will be secured by working only with the principal invariants $I_k$ and their derivatives. In order to see that right-local SO(3)-invariance \eqref{IntWGradGrad5} is satisfied we look exemplarily at
\begin{align}\label{NormGradInv1}
\norm{\gradx{I_1(C)}}^2_{\Skalar^3} = \norm{\gradx{\norm{\Gradx{\varphi(x)}}^2_{\Skalar^{3 \times 3}}}}^2_{\Skalar^3} \,.
\end{align}
We compute directly, following the steps leading to \eqref{IntWGradGrad5}
\begin{align}\label{NormGradInv1b}
\norm{\Gradx{\varphi(x)}}^2_{\Skalar^{3 \times 3}} &= \norm{\Gradxi{\varphi^{\flat}(\zeta^{-1}(x))} \, (\hspace{-0.8cm}\underbrace{\Gradx{\zeta^{-1}(x)}}_{\displaystyle (\Gradxi{\zeta(\xi)})^{-1} \in \SO(3)} \hspace{-0.8cm} ) }^2_{\Skalar^{3 \times 3}} =
\norm{\Gradxi{\varphi^{\flat}(\zeta^{-1}(x))} \, Q(\xi) }^2_{\Skalar^{3 \times 3}} \notag \\
& = \norm{\Gradxi{\varphi^{\flat}(\zeta^{-1}(x))} }^2_{\Skalar^{3 \times 3}}  \,,
\end{align}
where we have allowed that $Q(\xi) \in \SO(3)$ may be inhomogeneous.
Further, using indices
\begin{align}\label{NormGradInv1c}
&\gradx{ \norm{\Gradxi{ \varphi^{\flat}(\zeta^{-1}(x))}}^2 } \notag \\
&\qquad \qquad \qquad = 2 \, (\Gradxi{\varphi^{\flat}(\zeta^{-1}(x))})_{ij} \, (\GRADxi{\Gradxi{\varphi^{\flat}(\zeta^{-1}(x))}})_{ijb} \, (\hspace{-0.8cm}\underbrace{\Gradx{\zeta^{-1}(x)}}_{\displaystyle (\Gradxi{\zeta(\xi)})^{-1} \in \SO(3)} \hspace{-0.8cm} )_{bk} \, e_k  \,.
\end{align}
This implies finally
\begin{align}\label{NormGradInv1d}
\norm{\gradx{I_1(C)}}^2_{\Skalar^3} &= \norm{2 \, (\Gradxi{\varphi^{\flat}(\zeta^{-1}(x))})_{ij} \, (\GRADxi{\Gradxi{\varphi^{\flat}(\zeta^{-1}(x))}})_{ijb} ( \hspace{-0.6cm} \underbrace{\Gradxi{\zeta(\xi)}}_{\displaystyle \in \SO(3) \, \forall \xi \in \zeta^{-1}(B)} \hspace{-0.6cm} )^{-1}_{bk} \, e_k }^2_{\Skalar^3} \notag \\
 &=\norm{( 2 \, \Gradxi{\varphi^{\flat}(\zeta^{-1}(x))})_{ij} \, (\GRADxi{\Gradxi{\varphi^{\flat}(\zeta^{-1}(x))}})_{ijk} \, e_k }^2_{\Skalar^3} \notag \\
 &=\norm{ 2 \, \Gradxi{\varphi^{\flat}(\zeta^{-1}(x))} \, \GRADxi{\Gradxi{\varphi^{\flat}(\zeta^{-1}(x))}}}^2_{\Skalar^3} \notag \\
 &=\norm{ \gradxi{\norm{\Gradxi{\varphi^{\flat}(\zeta^{-1}(x))}}^2}}^2_{\Skalar^3} \notag \\
 &=\norm{ \gradxi{\tr\{(\Gradxi{\varphi^{\flat}(\zeta^{-1}(x))})^T \, \Gradxi{\varphi^{\flat}(\zeta^{-1}(x))} \} }}^2_{\Skalar^3} \notag \\
 &=\norm{ \gradxi{\tr C^{\flat}}}^2_{\Skalar^3}=\norm{ \gradxi{I_1(C^{\flat})}}^2_{\Skalar^3} \,,
\end{align}
which is form-invariance in the sense of condition \eqref{IntWGradGrad5} and implies therefore also condition \eqref{IntWGradGrad6}. Note again that we did not use constant rotations. Therefore, $\norm{\gradx{I_1(C)}}^2_{\Skalar^3}$ is already right-local SO(3)-invariant.
\subsubsection{Right-local SO(3)-invariant energies}
In the lack of having an encompassing representation theorem for right-local $\SO(3)$-invariant functions let us describe a large class of such expressions. For this we assume that
\begin{align}\label{PsiPsiCurv}
\Psi(C, \, \GRADx{C})  = \widehat{\Psi}(I_1(C), \, I_2(C), \, I_3(C)) + \widehat{\Psi}_{\rm curv}(\Gradx{\widehat{I}(C)}) \,,
\end{align}
where we defined the vector of principal invariants $\widehat{I}: \, \Sym^+(3) \mapsto \Euklid , \, \widehat{I}(C) \colonequals ((I_1(C), \, I_2(C), \, I_3(C))$
and we specify, furthermore, that\footnote{Note that the reduction in eq.\eqref{PsihatGradI} is not warranted by objectivity assumptions but just a simplification.}
\begin{align}\label{PsihatGradI}
\widehat{\Psi}_{\rm curv}(\Gradx{\widehat{I}(C)}) = \widehat{\widehat{\Psi}}_{\rm curv}((\Gradx{\widehat{I}(C)})^T \Gradx{\widehat{I}(C)}) \,.
\end{align}
Then right-local $\SO(3)$-invariance is entirely dependent upon properties of $\widehat{\widehat{\Psi}}$, namely
\begin{align}\label{PsihathatGradI}
\widehat{\widehat{\Psi}}_{\rm curv}( Q(\xi) \, \Gradx{\widehat{I}(C)}^T \Gradx{\widehat{I}(C)} \, Q(\xi)^T) = \widehat{\widehat{\Psi}}_{\rm curv}(\Gradx{\widehat{I}(C)}^T \Gradx{\widehat{I}(C)}) \quad \forall \, Q(\xi) \in \SO(3)  \,.
\end{align}
Thus, similar to the local case, $\widehat{\widehat{\Psi}}_{\rm curv}$ can be further reduced to a {\bf symmetric function} of the {\bf singular values} of $\Gradx{\widehat{I}(C)}$.\\

In fact, another representation, also reminiscent of the local theory suggests itself. Any $\widehat{\Psi}_{\rm curv}(\Gradx{\widehat{I}(C)})$ satisfies right-local SO(3)-invariance if an only if
\begin{align}\label{PsihatcurvInv2}
\widehat{\Psi}_{\rm curv}(\Gradx{\widehat{I}(C)} \, Q(\xi)) = \widehat{\Psi}_{\rm curv}(\Gradx{\widehat{I}(C)}) \, .
\end{align}
The latter may be used to obtain the reduced format
\begin{align}\label{PsihatcurvInv3}
\widehat{\Psi}_{\rm curv}(\Gradx{\widehat{I}(C)}) = \widetilde{\Psi}_{\rm curv}(\underbrace{\Gradx{\widehat{I}(C)} \, (\Gradx{\widehat{I}(C)})^T}_{\displaystyle \colonequals \widetilde{B}_{\nabla}} ) \,,
\end{align}
where
\begin{align}\label{PsihatcurvInv4}
\widetilde{B}_{\nabla} \colonequals \Gradx{\widehat{I}(C)} \, (\Gradx{\widehat{I}(C)})^T = \Gradx{\widehat{I}(B)} \,  (\Gradx{\widehat{I}(B)})^T \, , \qquad B=F \, F^T \,,
\end{align}
is a Finger-type stretch tensor based on $ \Gradx{\widehat{I}(B)}$.\footnote{In the local theory it would follow that $B \mapsto \sigma(B)$ is an isotropic tensor function and the strain energy is an isotropic function of $B$. However, this further reduction possibility is not viable here since we would need an additional objectivity requirement. This is absent in \eqref{PsihatcurvInv3}.}
\subsection{A problem with compatibility and left-local SO(3)-invariance}\label{KapProblemCompInv}
Let us approach the problem with an easy example. We have seen in Sect.\ref{KapOINonlin-Non-Local-Case}  that $\sum_{i=1}^3 \norm{F^T \, \frac{\partial F}{\partial x_i}}^2$ is frame-indifferent (left-global SO(3)-invariant) but not left-local SO(3)-invariant. However, using an involved compatibility argument \cite{Vallee92,LankeitNeffOsterbrink15}
one can show that
\begin{align}\label{FTpartialiF}
 \sum_{i=1}^3 \norm{F^T \, \frac{\partial F}{\partial x_i}}^2 = \widehat{\Psi}(\sqrt{F^T \, F} , \GRAD{\sqrt{F^T \, F}} ) = \widehat{\widehat{\Psi}}(C , \GRAD{C})  \, ,
\end{align}
and the right-hand side is clearly left-local SO(3)-invariant. In fact, it can be shown that every objective second gradient energy can be represented as\footnote{A further equivalent possibility is to use the "connection" $F^{-1} \Grad{F} = C^{-1} (F^T \Grad{F})$. This is done in Bertram \cite{Bertram2016}.}
\begin{align}\label{WFGradFPsi}
 W(F,\GRAD{F}) = \Psi(F^T F, F^T \GRAD{F}) =  \widehat{\widehat{\Psi}}(C , \GRAD{C})  \, ,
\end{align}
see Toupin \cite[eq.(10.24)+(10.26)]{Toupin64} and \cite[eq.(5.12)]{Toupin62}. This may lead one to assume that working solely with $\widehat{\widehat{\Psi}}(C , \GRAD{C})$ exhausts all possible fomulations for frame-indifferent second gradient models. Surprisingly, the situation is more complicated. Returning to \eqref{FTpartialiF} we observe
\begin{align}\label{FTpartialiF2}
 \sum_{i=1}^3 \norm{F^T \, \frac{\partial F}{\partial x_i}}^2 \qquad \quad &\neq \qquad \quad \sum_{i=1}^3 \norm{(Q(x) \, F)^T \, \frac{\partial (Q(x) \, F)}{\partial x_i}}^2 \notag \\
 \Big\| \,{\substack{\displaystyle F = \Grad{\varphi} \\ \displaystyle  \text{compatible}}}\quad & \hspace{2.3cm} \not{\! \Big\|}  \,\, {\substack{\displaystyle Q(x) \, F \neq \Grad{\vartheta} \\ \displaystyle  \text{incompatible}}} \notag \\
 \Psi \left(\sqrt{F^T F}, {\rm GRAD}\!\left[\!\sqrt{F^T F}\right]\right) \quad  &= \qquad  \widehat{\Psi}\left(\sqrt{(Q(x) \, F)^T (Q(x) \, F)}, {\rm GRAD}\!\left[\!\sqrt{(Q(x) \, F)^T (Q(x) \, F)}\right]\right) \,.
\end{align}
The situation is completely analogous to the linearized case in which e.g.
\begin{align}\label{axlskewCurlsym}
\frac{1}{4} \, \norm{\Grad{\curl(u)}}^2 = \norm{\Grad{\axl(\skew \Grad{u})}}^2 = \norm{\Curl \sym \Grad{u} }^2 \, ,
\end{align}
see eq.\eqref{Curvature_V2} and \eqref{Curvature_V3}, and both expressions are (left) global $\so(3)$-invariant, i.e. invariant w.r.t. the transformation $\Grad{u} \mapsto \overline{A} + \Grad{u} \, , \, \overline{A} \in \so(3)$. Equality in \eqref{axlskewCurlsym} holds only for gradients $\Grad{u}$ and uses compatibility. As equal energy expressions, both terms in \eqref{axlskewCurlsym} give rise to equivalent Euler-Lagrange equations, while the formulation based on $\Grad{\axl(\skew \Grad{u})}$ delivers {\bf non-symmetric force-stresses} and the formulation based on $\Curl \sym \Grad{u}$ gives {\bf symmetric} force stresses, see \eqref{strongform}. The difference resides in an extra stress contribution which has zero divergence, see \cite{GhiNeMaMue15}. Therefore, both formulas, although based on the same energy, may yield essentially different boundary value problems. The difference occurs for certain boundary conditions and load cases \cite{MaGhiNeMue15}. Thus it makes a fundamental difference starting with the left or right term in \eqref{axlskewCurlsym}, similarly with \eqref{FTpartialiF}. In summary,
\begin{align}\notag
\fbox{
\parbox[][2cm][c]{14.5cm}{
the complete model formulation based on left-global SO(3)-invariance is not necessarily captured by starting with a left-local SO(3)-invariant expression. Hence, not all frame-indifferent variational formulations can be obtained by assuming $W(F, \GRAD{F}) = \Psi(C, \GRAD{C})$ and working with the latter expression $\Psi(C, \GRAD{C})$.\footnote{The reader should be aware that the energy in a higher-gradient model does not completely determine the model. Decisive is the way in which partial integration steps are performed and this may change boundary conditions. Now, using $\norm{\Grad{\axl(\skew \Grad{u})}}^2$ suggests some partial integration steps which leads in a natural way to nonsymmetric total force stresses.}
}}
\end{align}
Moreover, while representation theorems for isotropic energies in the reduced format $\Psi(C, \GRAD{C})$ are well-known, cf. Sect.\ref{KapRayleighProduct} and \cite{dellIsola2009}, the more complete format is given by $\Psi(F^T F, \, F^T \GRAD{F})$ and representation theorems for isotropic models in this extended setting are not standard.
\section{A short look at objectivity and isotropy in linearized elasticity}\label{KapObjLinearElas}
\setcounter{equation}{0}
The form of the stress-strain law for linear elastic solids \cite{Love27} reads
\begin{align}\label{sigmaCeps}
\fg \sigma = \C \, \fg \varepsilon \, .
\end{align}
In order to discuss isotropy, the classical procedure is to consider $\fg \sigma$ and $\fg \varepsilon$ as linear mappings and refer these linear mappings to new axes of coordinates. We note that any linear mapping has a meaning independent of the used coordinate system. Then isotropy means that the formula connecting stress components with strain components is independent of the direction of the new axis. In this case the quadratic energy function $\Psi_{\rm lin}(\fg \varepsilon) = \12 \scal{\C \, \fg \varepsilon , \fg \varepsilon}$ is invariant for all transformations from one set of orthogonal axes to another:
\begin{align}\label{IsotropicDemand}
\fg \sigma^{\sharp} \colonequals \f Q \, \fg \sigma \, \f Q^T \, , \qquad \fg \varepsilon^{\sharp} \colonequals \f Q \, \fg \varepsilon \, \f Q^T \, , \qquad
\Psi_{\rm lin}(\fg \varepsilon^{\sharp}) \overbrace{=}^{\rm invariance} \Psi_{\rm lin}(\fg \varepsilon) \,.
\end{align}
Hence, isotropy is equivalent to
\begin{align}\label{IsotropicDemandAdd}
\Psi_{\rm lin}(\fg \varepsilon^{\sharp}) = \Psi_{\rm lin}(\f Q \, \fg \varepsilon \, \f Q^T) \overbrace{=}^{\rm isotropy}  \Psi_{\rm lin}(\fg \varepsilon) \quad \forall \, \f Q \in \SO(3) \,.
\end{align}
 In terms of the elasticity tensor, eq.\eqref{IsotropicDemandAdd} implies
\begin{align}\label{IsotropicDemandAdd2}
\C= \12 \, {\rm D}_{\fg \varepsilon}^2 \Psi_{\rm lin}(\fg \varepsilon) = \12 \, {\rm D}_{\fg \varepsilon^{\sharp}}^2 \Psi_{\rm lin}(\fg \varepsilon^{\sharp}) = \C^{\sharp} \,.
\end{align}
In an isotropic material, when we do not want to directly work with the elastic energy, the principal axes of stress $\fg \sigma$ and strain $\fg \varepsilon$ must coincide:
\begin{align}\label{IsotropicDemand2}
\C \,\fg \varepsilon^{\sharp} = \C \, (\f Q \, \fg \varepsilon \, \f Q^T ) = \f Q \, (\C \, \fg \varepsilon) \, \f Q^T \, , \qquad \fg \sigma^{\sharp} = \C^{\sharp} \, \fg \varepsilon^{\sharp} = \f Q \, \fg \sigma \, \f Q^T \,.
\end{align}
Up to this point, the isotropy condition in eq.\eqref{IsotropicDemand} is not directly related to our $\sharp$-transformation since $(\fg \sigma , \, \fg \varepsilon)$ are taken merely as linear mappings. Next we show that by properly introducing the $\sharp$-transformation, isotropy in linear elasticity can be completely characterized. Due to linear isotropic elasticity the simultaneous transformation of spatial and referential coordinates\footnote{Note that rotating the {\bf spatial} and {\bf referential} frame by the {\bf same} rotation for a linear elastic framework makes sense since linear elasticity does not distinguish between these two frames.} yields the Cauchy stress $\fg \sigma = 2 \, \mu \sym \Grad{\f u} + \lambda \, \tr (\Grad{\f u})\, \id $ in the transformed representation
\begin{align}
\fg \sigma^\sharp(\fg \xi) = \fg \sigma^\sharp(\Gradxi{\uxi}) &= 2 \, \mu \, \sym(\f Q \, \Gradxi{\uxi} \f Q^T) + \lambda \tr(\f Q \, \Gradxi{\uxi} \f Q^T) \id \notag \\
&= 2 \, \mu \, \sym(\f Q \, \Gradx{\f u(\f x)} \f Q^T) + \lambda \tr(\f Q \, \Gradx{\f u(\f x)} \f Q^T) \id \notag \\
&= \f Q \, \{ 2\, \mu \, \sym \Gradx{\ux} +\lambda \, \tr (\Gradx{\ux}) \id \} \f Q^T \notag \\
&= \f Q \, \fg \sigma (\Gradx{\ux}) \f Q^T = \f Q \, \sigmax \, \f Q^T \,.
\end{align}
But in any (rotated) coordinate system the elasticity tensor has the same format from assuming isotropy\footnote{Eringen \cite[p.~151]{Eringen62}: "An elastic material is said to be isotropic in its natural state if the constitutive equation [$\sigma^{\sharp}=Q \, \sigma \, Q^T$] is form-invariant with respect to a full orthogonal group of transformations of the natural state..."}
\begin{align}
\sigmax = \C  \, \varepsilonx \qquad {\rm and} \qquad \sigmaxi = \C^{\sharp} \, \varepsilonxi \, ,
\end{align}
where
\begin{align}
\sigmaxi = \f Q \, \sigmax \, \f Q^T  \qquad {\rm and} \qquad \varepsilonxi = \f Q \, \varepsilonx \, \f Q^T \, .
\end{align}
Thus, the requirement for isotropy is simply
\begin{align}\label{CCsharp}
\C = \C^{\sharp} \, ,
\end{align}
in terms of our $\sharp$-transformation. This yields
\begin{align}\label{QsigmaQ}
\f Q \, \sigmax \, \f Q^T &= \C^{\sharp}(\f Q \, \varepsilonx \, \f Q^T) \notag \\
\f Q \, (\C \, \varepsilonx) \, \f Q^T &= \C^{\sharp}(\f Q \, \varepsilonx \, \f Q^T) \quad \forall \, \f Q \in \SO(3) \quad \fg \varepsilon \in \Sym(3) \, .
\end{align}
Eq.\eqref{CCsharp} and \eqref{QsigmaQ} imply $\f Q \, (\C \, \fg \varepsilon) \, \f Q^T = \C \, (\f Q \, \fg \varepsilon \, \f Q^T)$. Since this is true for any $\fg \varepsilon \in \Sym(3)$ we have in components:
\begin{align}\label{CC}
\C_{ijkl} = Q_{ia} \, Q_{jb} \, Q_{kc} \, Q_{ld} \, \, \C_{abcd}   \, ,
\end{align}
the general solution, for symmetric tensors $\varepsilon$ is the well-known relation
\begin{align}\label{Ceps}
\C \, \fg \varepsilon = 2 \, \mu \, \fg \varepsilon + \lambda \, \tr(\fg \varepsilon) \, \id   \, ,
\end{align}
where $\mu$ is the shear modulus and $\lambda$ is the second Lam\'{e} parameter.

\noindent
In linear elasticity the quadratic energy $W_{\rm lin}$ is defined as a function of the displacement gradient
\begin{align}\label{Gradxvonu}
\Gradx{\ux} = \Grad{\phix} - \id \qquad {\rm for} \quad  \fg \varphi(\f x) = \f x + \f u(\f x)\,.
\end{align}
Linearized frame-indifference requires that
\begin{align}\label{Wlin}
W_{\rm lin}(\Grad{\f u} + \overline{\f W})= W_{\rm lin}(\Grad{\f u}) \qquad \forall \, \overline{\f W} \in \so(3) \,.
\end{align}
This implies the reduced representation
\begin{align}\label{Wlinreduced}
W_{\rm lin}(\Grad{\f u})= \Psi_{\rm lin}(\sym \Grad{\f u}) = \Psi_{\rm lin}(\fg \varepsilon) \, , \qquad \fg \varepsilon \colonequals \sym \Grad{\f u} \,.
\end{align}
Describing the deformation $\fg \varphi(\f x) = \f x + \f u(\f x) $ by the displacement field $\f u(\f x)$ one also gets
\begin{align}
\fg \varphi^\sharp(\fg \xi) &= \f  Q \, \fg \varphi (\f Q^T \, \fg \xi) = \f Q \, [\f Q^T \, \fg \xi + \f u(\f Q^T \, \fg \xi)] = \fg \xi + \f Q \, \f u(\f Q^T \, \fg \xi) = \fg \xi + \f u^\sharp(\fg \xi) \, ,
\end{align}
which shows that the tranformation in eq.\eqref{phisharpxi} is compatible with respect to displacements and deformations. Considering
\begin{align}
(\Gradxi{\fg \varphi^{\sharp}})^T \, \Gradxi{ \fg \varphi^{\sharp}} &= \f Q \, (\Gradx{\fg \varphi})^T \, \f Q^T \, \f Q \, \Gradx{\fg \varphi} \, \f Q^T= \f Q \, (\Gradx{\fg \varphi})^T \, \Gradx{\fg \varphi} \, \f Q^T \notag \\
&=\f Q \, [\id + 2 \, \sym \Gradx{\ux} + {\rm h.o.t.}]\, \f Q^T =\id + 2 \, \f Q \sym \Gradx{\ux} \,  \f Q^T  + {\rm h.o.t.} \notag \\
& =\id + 2 \,  \sym \f Q \, \Gradx{\ux} \,  \f Q^T  + {\rm h.o.t.} = \id + 2 \sym \Gradxi{\f u^\sharp(\fg \xi)} + {\rm h.o.t.} \, ,
\end{align}
we may use transformation \eqref{phisharpxi} for a geometrically linear displacement formulation and probe objectivity and isotropy in problems with linearized kinematics, which is investigated in the next section\footnote{Only transforming the referential coordinate system by defining $\fg \varphi^{\flat}(\fg \xi) \colonequals \fg \varphi(\f Q^T \, \fg \xi)$ is not compatible with the linear displacement formulation, since, $\Gradxi{\fg \varphi^{\flat}(\fg \xi)} = (\Gradx{\fg \varphi(\f x)}) \, \f Q$ and
\begin{align*}
(\Gradxi{\fg \varphi^{\flat}})^T \, \Gradxi{ \fg \varphi^{\flat}} &= \f Q \, (\Gradx{\fg \varphi})^T \, \Gradx{\fg \varphi} \, \f Q^T = \f Q \, [\id + 2 \, \sym ( \Gradx{\ux} + {\rm h.o.t.})\, \f Q^T  \notag \\
 &=\id + 2 \, \f Q \sym ( \Gradx{\ux} ) \,  \f Q^T  + {\rm h.o.t.} =\id + 2 \,  \sym ( \f Q \, \Gradx{\ux} \,  \f Q^T ) + {\rm h.o.t.}  \notag \\
 &\neq \id + 2 \, \sym \Gradxi{\f u^\flat(\fg \xi)} + {\rm h.o.t.} = \id + 2 \, \sym (\Gradx{\ux} \, \f Q) + {\rm h.o.t.}
\end{align*}
Note that $\fg \varphi^{\flat}(\fg \xi) = \fg \varphi^{\flat}(\f Q^T \, \fg \xi)$ is the correct transformation in nonlinear elasticity when we only want to discuss isotropy there as seen in Section \ref{KapOINonlin-Non-Local-Case}.}. As a preliminary, for the isotropic Cauchy stress tensor $\sigma$ in linear elasticity
\begin{align}
\fg \sigma = 2\, \mu \sym \Gradx{\ux} +\lambda \tr(\sym \Gradx{\ux}) \, \id \,,
\end{align}
we obtain
\begin{align}
\fg \sigma^\sharp(\fg \xi)= \f Q \, \fg \sigma(\f x) \, \f Q^T = 2\, \mu \sym \Gradxi{\f u^\sharp(\fg \xi)} +\lambda \tr(\sym \Gradxi{\f u^\sharp(\fg \xi)}) \, \id \,,
\end{align}
as expected.\\

\noindent
It is illuminating to see that, although in a linearized context, we may still use the full $\sharp$-transformation\footnote{c.f. \cite[eq.~3.18]{Bertram2016}.}
\begin{align}\label{usharpGradusharp}
\uxi \colonequals \f Q \, \f u(\f Q^T \, \fg \xi) \,, \qquad \Gradxi{\uxi} = \f Q \, \Gradx{\ux} \, \f Q^T
\end{align}
for the displacement $\f u$, since it is fully consistent with the corresponding transformation for the deformation:
\begin{align}\label{phisharpGradphisharp}
\phixi = \f Q \, \fg \varphi(\f Q^T \, \fg \xi) \,, \qquad \Gradxi{\phixi} = \f Q \, \Gradx{\phix} \, \f Q^T
\end{align}
and
\begin{align}\label{Gradphisharpminusone}
\Gradxi{\phixi} - \id &= \f Q \, \Gradx{\phix} \, \f Q^T  - \id = \f Q \, (\Gradx{\phix} - \id) \, \f Q^T
= \f Q \, \Gradx{\ux} \, \f Q^T \notag \\
&= \Gradxi{\uxi} \,.
\end{align}
With the help of the $\sharp$-transformation we may alternatively characterize frame-indifference and isotropy.

\begin{align}\label{Wlingradu}
\fbox{
\parbox[][3cm][c]{14.5cm}{
The linear hyperelastic formulation is (linearized) frame-indifferent and isotropic {\bf if and only if} there exists $\Psi_{\rm lin} : \Sym(3) \rightarrow \Skalar$ such that\\
\\
$
W_{\rm lin}(\Grad{\f u})  \overbrace{=}^{{\rm linearized} \atop {\rm frame-indifferent}} \Psi_{\rm lin}(\sym \Grad{\f u}) \,, \quad {\rm and}
$
\\
\\
$
W_{\rm lin}(\Gradxi{\uxi}) \overbrace{=}^{\rm form-invariance} W_{\rm lin}(\Gradx{\ux}) \quad \left[ {\rm i.e.} \quad \Psi_{\rm lin}(\f Q \, \fg \varepsilon \, \f Q^T) = \Psi_{\rm lin}(\fg \varepsilon) \, \forall \, \f Q \in \SO(3) \,. \right]
$
}}
\end{align}
The latter result is in complete correspondence to the finite strain setting in eq.\eqref{WvonGradxi}.
In conclusion: if we work already within the reduced format of
\begin{align}\label{Wlingradu}
W_{\rm lin}(\Grad{\f u}) = \Psi_{\rm lin}(\sym \Grad{\f u}) \, \rightarrow \, ({\rm linearized \, frame-indifference}) \,,
\end{align}
then, isotropy is nothing else but form-invariance of the energy $W_{\rm lin}$ under the $\sharp$-transformation. It is this simple observation that will be subsequently made operable for linear higher-gradient elasticity.

The reader should carefully note that linearized elasticity, as we treat it here, is not frame-indifferent in the precise sense that it is not left-invariant under rotations in $\SO(3)$. By this we mean that for $\fg \varphi(\f x) = \f x + \f u(\f x)$ one can write, formally, the linear elastic energy (wether isotropic or anisotropic) as
\begin{align}
W_{\rm lin}(\Grad{\fg \varphi}) = W_{\rm lin}\left(\sym ( \Grad{\fg \varphi} - \id) \right) \,,
\end{align}
and in this representation, $W_{\rm lin}$ is not left-invariant under the $\SO(3)$-action.\footnote{ The claim in \cite{Steigmann07b} on frame-invariance of linear elasticity seems to be due to a non-standard re-interpretation of objectivity, see also \cite{Yavari08} and \cite{FosdickSerrin79}.}
\section{Simultaneous rotation of spatial and referential coordinates}\label{KapRotInvariance}
\setcounter{equation}{0}
Let us consider the Lagrangian description of the displacement field
\begin{align}\label{uvonx}
\f u(\f x) = u_i(\f x) \, \f e_i \, ,
\end{align}
with components $u_i$ in the basis $\f e_i$ at the referential position $\f x$.
 \begin{align}\label{usharpdef}
\fbox{
\parbox[][2cm][c]{12cm}{
We investigate transformation rules for the {\bf simultaneous} transformation of spatial and referential coordinates by the {\bf same} rigid rotation $\f Q \in \SO(3)$
\begin{center}
$
\f u^{\sharp}(\f Q \, \f x) \colonequals \f Q \, \f u(\f x) \, , \quad \fg \xi := \f Q \, \f x \, , \quad \f x  := \f Q^T \, \fg \xi \,
$
\end{center}
\begin{center}
$
\uxi = \f Q \, \f u(\f Q^T \, \fg \xi)
\quad  \Leftrightarrow \quad \f u (\f x) = \f Q^T\,\uxi \, .
$
\end{center}
}}
 \end{align}

\begin{figure}[h]
\centering
 \begin{picture}(130,110)
  \put(0,0){\epsfig{file=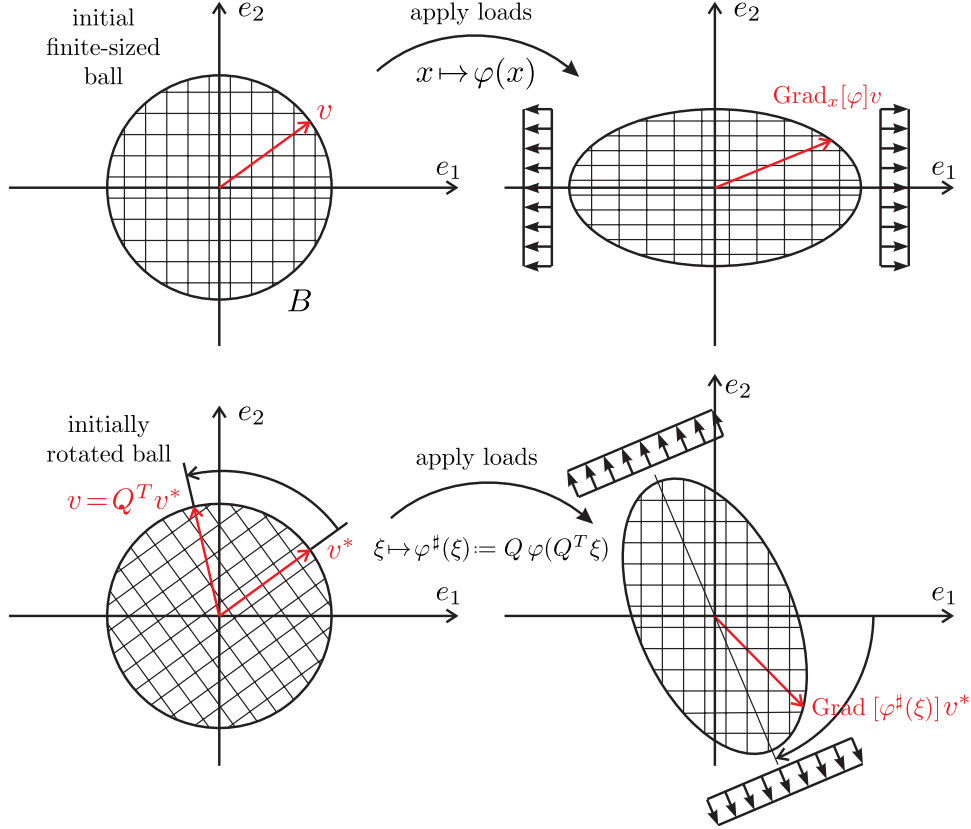, height=110mm, angle=0}}
 \end{picture}
\caption{Left: Consider a finite-sized ball $B$ in an undistorted initial configuration, which may be rotated by a constant rotation $Q \in SO(3)$. Right: Isotropy requires that deformation is rotated by $Q$ if we apply loading rotated by the same $Q$: the energy content again depends only on the shape of the ellipse.}
 \label{ObjectivityPic2}
\end{figure}

\noindent
Then the Jacobian matrix is given by the transposition of the rotation tensor:
\begin{align}\label{Jacobi}
\frac{\partial (x(\fg \xi))_i}{\partial \xi_j} =  \frac{\partial (Q^T_{ia} \, \xi_a)}{\partial \xi_j} = \f Q^T_{ia} \, \delta_{aj} =  Q^T_{ij} \, .
\end{align}
%
The simultaneous rotation of spatial and referential coordinates leaves the energy content stored due to deformation of a finite sized ball invariant. It depends only on the shape of the ellipse, not on its (irrelevant) inner structure nor its rotation in space, see Fig. \ref{ObjectivityPic2}.
\subsection{Standard transformation rules}\label{KapTrafoRules}
In this section the transformation \eqref{usharpdef} of standard vectorial and tensorial functions is under investigation. As already seen in eq.\eqref{Gradxiphi} the gradient of the displacement field $\f u$ transforms via
\begin{align}\label{Gradxiusharp}
\Gradxi{\uxi} = \Gradxi{ \f Q \, \f u (\f Q^T \, \fg \xi)} = \frac{\partial \, [ \f Q \, \f u (\f Q^T \, \fg \xi)]}{\partial \, \fg \xi}= \f Q \, \frac{\partial \, \f u (\f x)}{\partial \, \f x} \, \frac{\partial \, \f x}{\partial \, \fg \xi}= \f Q \, \Gradx{\f u (\f x)} \, \f Q^T  \,.
\end{align}
The divergence of the displacement field yields
\begin{align}\label{divxiuxi}
\divxi[\uxi] & = \tr(\Gradxi{\uxi}) = \tr(\Gradxi{\f Q \, \f u(\f Q^T \, \fg \xi)}) = \tr(\f Q \, \Gradx{\f u(\f x)} \, \f Q^T) \notag \\
&= \scal{\f Q \, \Gradx{\f u(\f x)} \, \f Q^T , \id} = \scal{ \Gradx{\f u(\f x)} , \f Q^T \, \f Q} = \scal{ \Gradx{\f u(\f x)} , \id}= \tr(\Gradx{\f u(\f x)}) \notag \\
& = \divx[\f u(\f x)] \,,
\end{align}
which is independent of $\f Q$. However, the divergence of the displacement gradient reads
\begin{align}\label{DivxiGradxiusharp}
\Divxi \{ \Gradxi{\uxi} \} &= (\GRADxi{ \Gradxi{\f Q \, \f u(\f Q^T \, \fg \xi)}})_{ijk} \, \delta_{jk} \, \f e_i = (\GRADxi{ \f Q \, \Gradx{\ux} \, \f Q^T })_{ijk} \, \delta_{jk} \, \f e_i \notag \\
&= \left( \frac{\partial \, (\f Q \, \Gradx{\ux} \, \f Q^T)_{ij}}{\partial \, x_a} \, \underbrace{\frac{\partial \, x_a}{\partial \, \xi_k}}_{\displaystyle \f Q^T} \right)_{ijk} \!\!\!\!\!\! \delta_{jk} \, \f e_i = (\GRADx{\f Q \, \Gradx{\ux} \, \f Q^T})_{ija} \, \underbrace{(\f Q^T)_{ak} \, \delta_{jk}}_{\displaystyle Q^T_{aj}} \, \f e_i  \notag \\
&= (Q_{im} ({\Gradx{\ux}})_{mb} \, Q^T_{bj})_{,a} \, Q^T_{aj} \, \f e_i = Q_{im} (\GRADx{\Gradx{\ux}})_{mba} \, \underbrace{Q^T_{bj} \, Q_{ja}}_{\displaystyle \delta_{ba}} \, \f e_i \notag \\
& = Q_{im} (\GRADx{\Gradx{\ux}})_{mbb} \, \f e_i = \f Q \, \Divx \{ \Gradx{\ux} \} \,.
\end{align}
The curl of the displacement gradient is given by definition \eqref{CurlDefTensor} for a second order tensor $\f X \in \Skalar^{3\times3}$. From eq.\eqref{CurlGradVect} it follows
\begin{align}\label{CurlGradu}
 \Curlx \{ \Gradx{\ux} \}& =  \Curlxi\{ \Gradxi{\uxi} \}   = 0 \, .
\end{align}
Of course, eq.\eqref{CurlGradu} is not a transformation rule since the Curl of the gradient of the displacement field is generally zero. However, the transformation rule for an arbitrary argument which is not the gradient of a vector field will be discussed in section \ref{KapTrafoEnhanced}. Before, let us see that the transformation \eqref{usharpdef} is compatible with the application of a sym or skew operator to the gradient of the displacement:
\begin{align}\label{symGraduxi}
\sym \{\Gradxi{\uxi} \} &=\sym \{\Gradxi{\f Q \, \f u(\f Q^T \, \fg \xi)} \} = \sym \{\f Q \, \Gradx{\f u(\f x)} \, \f Q^T \} \notag \\
& = \12 \{\f Q \, \Gradx{\f u(\f x)} \, \f Q^T + \f Q \, \Gradx{\f u(\f x)}^T \, \f Q^T \} = \f Q \, \sym \{ \Gradx{\f u(\f x)} \} \, \f Q^T \,,
\end{align}
\begin{align}\label{skewGradu}
\skew \{\Gradxi{\uxi} \} &=\skew \{\Gradxi{\f Q \, \f u(\f Q^T \, \fg \xi)} \} = \skew \{\f Q \, \Gradx{\f u(\f x)} \, \f Q^T \}  \notag \\
&= \12 \{\f Q \, \Gradx{\f u(\f x)} \, \f Q^T - \f Q \, \Gradx{\f u(\f x)}^T \, \f Q^T \} = \f Q \, \skew \{ \Gradx{\f u(\f x)} \} \, \f Q^T \,.
\end{align}
From eq.\eqref{divxiuxi} it follows that the trace of $\Grad{\f u}$  is invariant under the $\sharp$-transformation, yielding
\begin{align}\label{trsymskewgradu}
\tr \{ \sym ( \Gradxi{\uxi}  ) \} &= \tr \{ \Gradxi{\uxi} \} = \tr \{ \Gradx{\ux}  \} \,,\\
\tr \{ \skew ( \Gradxi{\uxi}  ) \} &= \tr \{ \skew ( \Gradx{\ux}  ) \} = 0  \,.
\end{align}
Thus, the deviator of the symmetric gradient of $\f u$ transforms accordingly as
\begin{align}\label{devsymgradu}
\dev \sym ( \Gradxi{\uxi}  ) & = \f Q \, \sym \{ \Gradx{\f u(\f x)} \} \, \f Q^T  - \frac{1}{3} \, \tr \{ \Gradx{\ux} \}  \, \id \notag \\
& = \f Q \, \{ \sym ( \Gradx{\f u(\f x)} ) - \frac{1}{3} \, \tr \{ \Gradx{\ux} \} \id \} \, \f Q^T \notag \\
& = \f Q \, \{ \dev \sym ( \Gradx{\f u(\f x)}) \} \, \f Q^T
 \,.
\end{align}
The divergence of the symmetric and skew-symmetric part of  $\Grad{\f u}$ then transforms according to
\begin{align}\label{Divsymgradu}
\Divxi \{ \sym ( \Gradxi{\uxi}  ) \} &= \Divxi \{ \sym ( \Gradxi{\f Q \, \f u(\f Q^T \, \fg \xi)}  ) \} = \Divxi  \{ \sym ( \f Q \, \Gradx{\ux} \, \f Q^T ) \} \notag \\
& = \Divxi  \{ \f Q \, \sym \Gradx{\ux} \, \f Q^T \}
= \f Q \, \Divx \{ \sym \Gradx{\ux} \} \, ,
\end{align}
\begin{align}\label{Divskewgradu}
\Divxi \{ \skew ( \Gradxi{\uxi}  ) \} &= \Divxi \{ \skew ( \Gradxi{\f Q \, \f u(\f Q^T \, \fg \xi)}  ) \}= \Divxi  \{ \skew ( \f Q \, \Gradx{\ux} \, \f Q^T ) \}  \notag \\
&=  \Divxi  \{ \f Q \, \skew \Gradx{\ux} \, \f Q^T ) \} = \f Q \, \Divx \{ \skew \Gradx{\ux} \} \,.
\end{align}
To investigate the axial vector $\f a \in\R^3$ of a skew symmetric tensor $\f A \in \so(3)$ we use eq.\eqref{axlanti} and  consider the components of the tensor $\f A$ in the basis $\f e_k$, reading
\begin{align}\label{aaxlA}
\f a &= \axl \f A = -\12 \, A_{ij}\,\epsilon_{ijk} \, \f e_k \, .
\end{align}
Since the functions axl and anti are free of a basis their results depend on the basis of the argument, which may be rotated:
\begin{align}\label{asharp}
\f a^\sharp &= \f Q \, \f a = \f Q \, \axl[\f A] = Q_{ak}(\axl[\f A])_k \, \f e_a = -\12 \, Q_{ak} \, A_{ij} \, \epsilon_{ijk} \, \f e_a = -\12 \, A_{ij} \, \underbrace{Q_{ak} \, \epsilon_{kij}}_{\mathclap{\displaystyle Q_{bi} \, Q_{cj} \, \epsilon_{bca}}} \, \f e_a
= -\12 \, Q_{bi} \, A_{ij} \, Q^T_{jc} \, \epsilon_{bca}  \, \f e_a \notag \\
 &= \axl[\underbrace{\f Q \, \f A \, \f Q^T}_{\displaystyle \f A^\sharp}] \qquad \Leftrightarrow \qquad \axl[\f A^\sharp] = \f Q \, \axl[\f A] \, .
\end{align}
Two transformation rules can be found in eq.\eqref{asharp}: the transformation for the axial-operator and how the simultaneous transformation of spatial and referential coordinates transform a tensor, reading
\begin{align}\label{TransTensor}
\f A^\sharp = \f Q \, \f A \, \f Q^T \, .
\end{align}
Additionally, the transformation of $\anti : \R^3 \rightarrow \so(3)$ is given by eq.\eqref{asharp}, reading
\begin{align}\label{antiasharp}
\anti[\f a^\sharp] = \anti  [\f Q \, \f a] = \anti \axl [\f A^\sharp] = \anti \axl [\f Q \, \f A \, \f Q^T] = \f Q \, \f A \, \f Q^T \, \qquad \Leftrightarrow \qquad \anti[\f a^\sharp] = \f Q \, \anti[\f a] \f Q^T \, .
\end{align}
With help of eq.\eqref{skewGradu} and \eqref{asharp} we obtain for the curl of the displacement field $\f u$:
\begin{align}\label{curlTrans}
\curlxi[\uxi] & = 2 \axl \skew [\Gradxi{\uxi}] = 2 \axl \skew [\f Q \, \Gradx{\f u(\f x)}\, \f Q^T] = 2 \axl  [\f Q \, \skew \Gradx{\f u(\f x)} \, \f Q^T] \notag \\
& = \f Q \, [ 2  \axl \skew \Gradx{\f u(\f x)} ] = \f Q \, \curlx [\f u(\f x)] \, .
\end{align}
The gradient of the curl of the displacement field transforms via
\begin{align}\label{GradcurlTrans}
\Gradxi{\curlxi[\uxi]} &= \Gradxi{\f Q \, \underbrace{\curlx \f u(\f Q^T \, \fg \xi)}_{\displaystyle \colonequals \f v(\f Q^T \, \fg \xi)}} = \Gradxi{\f Q \, \f v(\f Q^T \, \fg \xi)} = \f Q \, \Gradx{\f v(\f x)} \, \f Q^T \notag \\
&= \f Q \, \Gradx{\curlx \ux} \, \f Q^T \, .
\end{align}
Thus,
the transformation of spatial and referential coordinates by the same rigid rotation yields
the curvature $\ks = \12 \, \Grad{\curl(\f u)}$ from eq.\eqref{Curvature_V1} to transform via
\begin{align}\label{kssharpxi}
\ks^\sharp(\fg \xi) & = \12 \, \Gradxi{\curlxi(\uxi)}
= \f Q \,  \underbrace{\12 \, \Gradx{\curlx \f u(\f x)}}_{\displaystyle = \, \ks(\f x)} \, \f Q^T= \f Q \, \ks(\f x) \, \f Q^T \,  .
\end{align}
The symmetric and skew symmetric part of $\ks^\sharp(\fg \xi)$ transform as
\begin{align}\label{symkxi}
\sym[\ks^\sharp(\fg \xi)] = \12 \, [\ks^\sharp(\fg \xi) + (\ks^\sharp(\fg \xi))^T] = \f Q \, \12 \, [\ks(\f x) +(\ks(\f x))^T] \, \f Q^T = \f Q \, \sym[\ks(\f x)] \, \f Q ^T \,,
\end{align}
\begin{align}
\skew[\ks^\sharp(\fg \xi)] = \12 \, [\ks^\sharp(\fg \xi) - (\ks^\sharp(\fg \xi))^T] = \f Q \, \12 \, [\ks(\f x) -(\ks(\f x))^T] \, \f Q^T = \f Q \, \skew[\ks(\f x)] \, \f Q ^T \,.
\end{align}
Since the trace of $\ks$ is zero, it is independent of $\f Q$ with
$\tr \ksxi = \tr \ksx =  0 $.
%
%
\subsection{Advanced transformation rules}\label{KapTrafoEnhanced}
For the transformations in this section we often use the Cauchy stress tensor $\fg \sigma$ as the argument, which is a second order tensor. The divergence lowers a tensor of second order, transforming according to
\begin{align}\label{divxisigmaxi}
\Divxi[\sigmaxi] &= \Divxi[ \f Q \, \sigmax \, \f Q^T ] = \left( \frac{\partial \, ( \f Q \, \sigmax \, \f Q^T)_{ij}}{\partial \, x_a} \, \frac{\partial \, x_a}{\partial \, \xi_k}\right)_{ijk} \!\!\!\!\! \delta_{jk} \, \f e_i = Q_{ib} \, \sigma_{bc,a} \, Q^T_{cj} \, Q^T_{ak} \, \delta_{jk} \, \f e_i \notag \\
&= Q_{ib} \, \sigma_{bc,a} \, \underbrace{Q^T_{cj} \, Q_{ja}}_{\displaystyle \delta_{ca}} \, \f e_i  = Q_{ib} \, \sigma_{ba,a} \, \f e_i = \f Q \, \Divx[ \sigmax ] \,.
\end{align}
In contrast, the gradient lifts a second order tensor to a tensor of third order transforming via
\begin{align}\label{gradxisigmasharpxi}
\GRADxi{\sigmaxi} &= \GRADxi{\f Q \, \sigmax \, \f Q^T}
= \frac{\partial \, ( \f Q \, \sigmax \, \f Q^T)_{ij}}{\partial \, x_a} \, \frac{\partial \, x_a}{\partial \, \xi_k} \, \f e_i \otimes \f e_j \otimes \f e_k  \notag \\
&= (\GRADx{\f Q \, \sigmax \, \f Q^T})_{ija} \, Q^T_{ak} \, \f e_i \otimes \f e_j \otimes \f e_k \notag \\
&=(Q_{ib} \, \sigma_{bc} \, Q^T_{cj})_{,a} \, Q^T_{ak} \, \f e_i \otimes \f e_j \otimes \f e_k \notag \\
&=Q_{ib} \, \sigma_{bc,a} \, Q^T_{cj} \, Q^T_{ak} \, \f e_i \otimes \f e_j \otimes \f e_k \notag \\
&= \GRADx{\f Q \, \sigmax \f Q^T} \, \f Q^T \,.
\end{align}
Since the divergence lowers the order of $\Grad{\fg \sigma}$ to second order, one obtains
\begin{align}\label{divxigradisigmasharp}
\DIVxi[\GRADxi{\sigmaxi}] &= \DIVxi[ \GRADx{\f Q \, \fg \sigmax \, \f Q^T} \, \f Q^T] \notag \\
&= \left( \frac{\partial \, (\GRADx{\f Q \, \sigmax \, \f Q^T} \, \f Q^T)_{ijk}}{\partial \, x_n} \, \frac{\partial \, x_n}{\partial \, \xi_m} \right)_{ijkm} \!\!\! \delta_{km} \, \f e_i \otimes \f e_j
 \notag \\
&=  (\GRADx{\GRADx{\f Q \, \sigmax \, \f Q^T} \, \f Q^T})_{ijkn} \, Q^T_{nm} \, \delta_{km} \, \f e_i \otimes \f e_j \notag \\
 & =  [(\GRADx{\f Q \, \sigmax \, \f Q^T})_{ija} \, Q^T_{ak}]_{,n} \, Q_{kn} \, \f e_i \otimes \f e_j \notag\\
 &=  (\GRADx{\f Q \, \sigmax \, \f Q^T})_{ija,n} \, \underbrace{\f Q^T_{ak} \, Q_{kn}}_{\displaystyle \delta_{an}}  \, \f e_i \otimes \f e_j \notag \\
& =  (\GRADx{\f Q \, \sigmax \, \f Q^T})_{ijn,n} \, \f e_i \otimes \f e_j
\notag \\
& =  (\f Q \, \sigmax \, \f Q^T)_{ij,nn}  \, \f e_i \otimes \f e_j =   Q_{ia} \,  (\sigmax)_{ab,nn} \, Q^T_{bj}  \, \f e_i \otimes \f e_j \notag \\
&= \DIVx[\GRADx{\f Q \, \fg \sigmax \f Q^T}] = \f Q \, ( \DIVx \GRADx{\sigmax} ) \, \f Q^T \,.
\end{align}
Similarly, the gradient of the divergence of $\fg \sigma$ yields
\begin{align}\label{gradxidivxisigmasharp}
\Gradxi{\Divxi \sigmaxi} &= \Gradxi{\f Q \, \Divx \sigmax } = \frac{\partial \, (\f Q \, \Divx \sigmax)}{\partial \, \f x} \, \frac{\partial \, \f x}{\partial \, \fg \xi} = \Gradx{\f Q \, \Divx \sigmax} \, \f Q^T \notag \\
&= \f Q \, \Gradx{\Divx \sigmax} \, \f Q^T \,.
\end{align}
Next, let us discuss the transformation of the Curl operator for arbitrary second order tensors in $ \R^{3\times 3}$. Therefore, we consider a generally non-symmetric field $\f P: \Skalar^3 \mapsto \Skalar^{3 \times 3}$. The transformation to rotated coordinates is given, similarly to the gradient $\Grad{\f u}$, see eq.\eqref{Gradxiusharp}, by $\Pxi \colonequals \f Q \, \f P(\f Q^T \, \fg \xi) \, \f Q^T$. Then
\begin{align}\label{Curlxip}
\Curlxi[\Pxi] = \Curlxi[\f Q \, \f P(\f Q^T \, \fg \xi) \, \f Q^T] = \f Q \, \Curlx[\Px] \, \f Q^T  \,.
\end{align}
This last formula can be seen in several ways. First, with the help of eq.\eqref{LeviCivitaQEpsQQ} we obtain
\begin{align}\label{CurlSigma}
\Curlxi[\f P^{\sharp}(\fg \xi)] &= \Curlxi[\f Q \, \f P(\f Q^T \, \fg \xi) \, \f Q^T]
= - \left( \frac{\partial \, (\f Q \, \Px \, \f Q^T)_{ia}}{\partial \, x_c} \, \frac{\partial \, x_c}{\partial \, \xi_b} \right)_{iab} \, \epsilon_{abj} \, \f e_i \otimes \f e_j \notag \\
&=- (\Gradx{\f Q \, \Px \, \f Q^T})_{iac} \, Q^T_{cb} \, \epsilon_{abj} \, \f e_i \otimes \f e_j = - (\f Q \, \Px \, \f Q^T)_{ia,c} \, Q^T_{cb} \, \epsilon_{abj} \, \f e_i \otimes \f e_j \notag \\
&= - Q_{im} \, P_{mn,c} \, Q^T_{na} \, Q^T_{cb} \, \epsilon_{abj} \, \f e_i \otimes \f e_j = - Q_{im} \, P_{mn,c} \, \underbrace{Q_{an} \, Q_{bc} \, \epsilon_{abj}}_{\displaystyle Q_{jp} \, \epsilon_{pnc}} \, \f e_i \otimes \f e_j  \notag \\
& = Q_{im} \, \underbrace{(- P_{mn,c} \, \epsilon_{ncp})}_{\displaystyle (\Curlx[\Px])_{mp}} \, Q^T_{pj} \, \f e_i \otimes \f e_j =  \f Q \, \Curlx[\f P(\f x)] \, \f Q^T \,.
\end{align}
Alternatively, from \cite{Neff_Chelminski07_disloc}, p. 318, we have the transformation law
\begin{align}\label{NeffTrafoLaw1}
\Curlxi[\f P(\underbrace{\f Q^T \, \fg \xi}_{\displaystyle \f x}) \, \f Q^T] = ( \Curlx[\f P(\f x)] ) \, \f Q^T \,.
\end{align}
Combining eq.\eqref{NeffTrafoLaw1} with
\begin{align}\label{NeffTrafoLaw2}
\Curlxi[\f Q \, \f P(\f Q^T \, \fg \xi)] = \Lin_{\f P(\f Q^T \, \fg \xi)} \underbrace{(\nablaxi \, \f Q)}_{\displaystyle = \, 0} + \f Q \, \Curlxi[\f P(\f Q^T \, \fg \xi)] \,,
\end{align}
where $\Lin_{\f P(\f Q^T \, \fg \xi)}$ is a linear operator, we obtain formally
\begin{align}\label{NeffTrafoLaw3}
\Curlxi[\f Q \, \f P(\f Q^T \, \fg \xi) \, \f Q^T] =\f Q \, \Curlxi[\f P(\f x) \, \f Q^T] = \f Q \, \Curlx[\f P(\f x)] \, \f Q^T \,.
\end{align}
Obviously, repeating Curl operators yield the same structure of transformation as in eq.\eqref{Curlxip}, e.g.
\begin{align}\label{CurlCurlSigma}
\Curlxi\Curlxi[\sigmaxi] & = \Curlxi [ \f Q \, \underbrace{\Curlx [\sigmax]}_{\displaystyle \colonequals \Px}  \, \f Q^T ] =  \f Q \, \Curlx[\Px] \, \f Q^T  = \f Q \, \Curlx[\Curlx[\sigmax]] \, \f Q^T
 \,.
\end{align}
Note that the Curl of a second order tensor field has generally no divergence
\begin{align}\label{divxcurlxsigma}
\Div ( \Curl[\fg \sigma] ) & = \Div( - \sigma_{ia,b} \, \epsilon_{abj} \, \f e_i \otimes \f e_j )  = - \sigma_{ia,bk} \, \epsilon_{abj} \, \delta_{jk} \, \f e_i =  - \sigma_{ia,bk} \, \epsilon_{abk} \, \f e_i \notag \\
&= - (\sigma_{i1,23}-\sigma_{i1,32}+\sigma_{i2,31}-\sigma_{i2,13}+\sigma_{i3,12}-\sigma_{i3,21}) \, \f e_i = 0 \,,
\end{align}
which is independent of $\f Q$, reading $\Divxi(\Curlxi[\sigmaxi]) = \Divx(\Curlx[\fg \sigmax]) = 0 $.\\
%
%
Let us reformulate the result in eq.\eqref{gradxisigmasharpxi} by using $\Grad{\f u}$ instead of $\fg \sigma$ as second order argument:
\begin{align}\label{gradxigradxiusharp}
\m^{\sharp}(\fg \xi) \cong \GRADxi{\Gradxi{\uxi}} &= \GRADxi{\f Q \, \Gradx{\ux} \, \f Q^T} \notag \\
&=(\GRADx{\f Q \, \Gradx{\ux} \, \f Q^T})_{ija} \, Q^T_{ak} \, \f e_i \otimes \f e_j \otimes \f e_k \notag \\
&=(Q_{ib} \, u_{b,c} \, Q^T_{cj})_{,a} \, Q^T_{ak} \, \f e_i \otimes \f e_j \otimes \f e_k \notag \\
&=Q_{ib} \, u_{b,ca} \, Q^T_{cj} \, Q^T_{ak} \, \f e_i \otimes \f e_j \otimes \f e_k \notag \\
&= (\GRADx{\f Q \, \Gradx{\ux} \, \f Q^T}) \, \f Q^T \,.
\end{align}
Then, lowering the tensorial order of $\GRAD{\Grad{\f u}} $ by the divergence yields
\begin{align}\label{divxigradxigradxiuxi}
\DIVxi\{\GRADxi{\Gradxi{\uxi}}\} &=\DIVxi\{\GRADx{\underbrace{\f Q \, \Gradx{\ux}\, \f Q^T}_{\displaystyle \colonequals \Px}}\,\f Q^T \} \notag \\
&=\left( \frac{\partial \, (\GRADx{\Px} \, \f Q^T)_{ijk}}{\partial \, x_a} \, \frac{\partial \, x_a}{\partial \, \xi_l} \right)_{ijkl} \, \delta_{kl} \, \f e_i \otimes \f e_j \notag \\
 & = (\GRADx{\Px})_{ijb,a} \, \underbrace{Q^T_{bk} \, Q^T_{al} \, \delta_{kl}}_{\displaystyle \delta_{ba}} \, \f e_i \otimes \f e_j = (\GRADx{\Px})_{ija,a} \, \f e_i \otimes \f e_j \notag \\
& = [\Px]_{ij,aa} \, \f e_i \otimes \f e_j =  \DIVx\{\GRADx{\f Q \, \Gradx{\ux}\, \f Q^T} \} \,.
\end{align}
Since $\f Q$ is a constant rotation  tensor it is possible to reformulate the result of eq.\eqref{divxigradxigradxiuxi} into
\begin{align}\label{divxigradxigradxiuxi2}
\DIVxi\{\GRADxi{\Gradxi{\uxi}}\} &= \DIVx\{\GRADx{\f Q \, \Gradx{\ux}\, \f Q^T} \} = [\f Q \, \Gradx{\ux}\, \f Q^T]_{ij,aa} \notag \\
&= ( Q_{ib}  \, (\Gradx{\ux})_{bc} \, Q_{cj}^T )_{,aa} \, \f e_i \otimes \f e_j = Q_{ib}  \, (\Gradx{\ux})_{bc,aa} \, Q_{cj}^T \notag \\
 &= \f Q \, \DIVx \{\GRADx{\Gradx{\ux}} \} \f Q^T = \f Q \, \Gradx{\Divx\{ \Gradx{\ux} \} } \f Q^T \,.
\end{align}
Finally, with help of eq.\eqref{divxigradxigradxiuxi2} let us lower the tensorial order of $\DIV \Grad{\Grad{\f u}} $ by the divergence again to find
\begin{align}\label{divxidivxigradxigradxiuxi}
\Divxi ( \DIVxi \{\GRADxi{\Gradxi{\uxi}}\} ) &= \Divxi ( \f Q \, \underbrace{\Gradx{\Divx\{ \Gradx{\ux} \} }}_{\displaystyle \colonequals \Px} \f Q^T ) \notag \\
&=\left( \frac{\partial \, (\f Q \, \Px \, \f Q^T)_{ij}}{\partial \, x_a} \, \frac{\partial \, x_a}{\partial \, \xi_k} \right)_{ijk} \, \delta_{jk} \, \f e_i =  (\f Q \, \Px \, \f Q^T)_{ij,a} \, Q^T_{ak} \, \delta_{jk} \, \f e_i \notag \\
&=  Q_{ib} \, (\Px)_{bc,a} \, Q^T_{cj} \, Q_{ja} \, \f e_i =  Q_{ib} \, (\Px)_{bc,a} \, \delta_{ca} \, \f e_i   =  Q_{ib} \, (\Px)_{ba,a} \, \f e_i\notag \\
&=\f Q \, \Divx(\Px) =  \f Q \, \Divx (\Gradx{\Divx\{ \Gradx{\ux} \} }) \,.
\end{align}
\subsection{Form-invariance of the linear momentum balance equation}\label{KapTrafoLinMom}
In this subsection we consider the question of form-invariance of the balance equations. To start these considerations, we consider the Laplace equation. Note that the Laplacian $\Laplace$ of a scalar field $h(\f x) = h(\f Q^T \fg \xi) = h^\sharp(\fg \xi)$ is {\bf form-invariant} concerning our simultaneous transformation of spatial and referential coordinates by the same rigid rotation:
\begin{align}
\Laplace_{\fg \xi}h^\sharp(\fg \xi) & \colonequals \divxi \, \gradxi{h^\sharp(\fg \xi)} =  \divxi \, \gradxi{h(\f Q^T \, \fg \xi)} = \divxi \left[ \frac{\partial \, h(\f Q^T \, \fg \xi)}{\partial \fg \xi} \right] = \divxi \left[ \frac{\partial \, h(\f x)}{\partial x_i} \, \frac{\partial \, x_i}{\partial \, \xi_j} \, \f e_j \right] \notag \\
&= \divxi[h_{,i} \, Q^T_{ij} \, \f e_j] = \divxi[\f Q \, \gradx{h(\f x)}] = \scal{\frac{\partial \, (\f Q \, \gradx{h(\f x)})}{\partial \, \f x} \, \frac{\partial \, \f x}{\partial \, \fg \xi} , \id} \notag \\
& = \scal{\f Q \, \Gradx{\gradx{h(\f x)}} \, \f Q^T, \id} = \scal{\Gradx{\gradx{h(\f x)}} , \f Q^T \, \f Q} = \scal{\Gradx{\gradx{h(\f x)}} , \id} \notag \\
& = \tr\Gradx{\gradx{h(\f x)}} = \divx \, \gradx{h(\f x)} = \Laplace_{\f x}h(\f x) \, .
\end{align}
Therefore, the Laplace-equation is form-invariant with respect to the $\sharp$-transformation for scalars\footnote{The Dirichlet-integral for scalar-valued functions is easily seen to be {\bf form-invariant} under coordinate rotations:
\begin{align*}
\int_B \frac{1}{2} \, \norm{\gradx{h(x)}}^2_{\Euklid} \, \di x \, , \quad h: \, \Euklid \rightarrow \Skalar \, .
\end{align*}
Consider $h^{\flat}(\xi) \colonequals h(Q^T \, \xi)$. Then $\gradxi{h^{\flat}(\xi)}$ is a row-vector and it holds
\begin{align*}
\gradxi{h^{\flat}(\xi)} = \gradx{h(Q^T \, \xi)} \, Q^T \, .
\end{align*}
Moreover,
\begin{align*}
 \frac{1}{2} \, \norm{\gradxi{h^{\flat}(\xi)}}^2_{\Euklid} &= \frac{1}{2} \, \norm{(\gradxi{h^{\flat}(\xi)})^T}^2_{\Euklid} = \frac{1}{2} \, \norm{(\gradx{h(Q^T \, \xi)} \, Q^T)^T}^2_{\Euklid} = \frac{1}{2} \, \norm{Q \, (\gradx{h(Q^T \, \xi)})^T}^2_{\Euklid} \notag \\
 & = \frac{1}{2} \, \norm{ (\gradx{h(x)})^T}^2_{\Euklid} = \frac{1}{2} \, \norm{ \gradx{h(x)}}^2_{\Euklid} \, .
\end{align*}
Therefore
\begin{align*}
\int_{x \in B} \frac{1}{2} \, \norm{\gradx{h(x)}}^2_{\Euklid} \, \di x  = \int_{\xi \in B} \frac{1}{2} \, \norm{\gradxi{h^{\flat}(\xi)}}^2_{\Euklid} \, \det Q \di \xi = \int_{\xi \in B} \frac{1}{2} \, \norm{\gradxi{h^{\flat}(\xi)}}^2_{\Euklid} \, \di \xi \,,
\end{align*}
which is form-invariance and in terms of our previous terminology right-global SO(3)-invariance.}, as is well-known.\\

\noindent
The classical balance of linear momentum in linear elasticity includes the divergence of the Cauchy stress $\fg \sigma$ and the net force $\f f$ which both transform via
\begin{equation}\label{BLM1Add}
 {\begin{array}{lrcl}
& \Divxi[ \sigmaxi ] + \f f^{\sharp}(\fg \xi) &=& 0 \, , \qquad \qquad \f f^{\sharp}(\, \fg \xi) \colonequals \f Q \, \f f(\f Q^T \, \fg \xi) \vspace{0.2cm} \\
\Leftrightarrow \quad & \Divxi[\f Q \, \fg \sigma(\f Q^T \, \fg \xi) \, \f Q^T] + \f Q \, \f f(\f Q^T \, \fg \xi) &=& 0 \vspace{0.2cm} \\
\Leftrightarrow \quad & \Divxi[\f Q \, \sigmax \, \f Q^T] + \f Q \,  \f f(\f x) &= & 0 \vspace{0.2cm} \\
\Leftrightarrow \quad & \f Q \, \Divx[\sigmax] + \f Q \,  \f f(\f x) &=& 0 \vspace{0.2cm} \\
\Leftrightarrow \quad & \f Q \, \{ \Divx[\sigmax] + \f f(\f x) \} &=& 0 \vspace{0.2cm} \\
\Leftrightarrow \quad & \Divx[\sigmax] + \f f(\f x) &=& 0 \,.
 \end{array}}
\end{equation}
Therefore, the linear elasticity equations are form-invariant with respect to the $\sharp$-transformation.

The transformation with $\f Q$ from the left also holds for the divergence of second order couple stress $\f m \sim \ks^\sharp(\fg \xi)$ concerning its symmetric part
\begin{align}
\Divxi(\sym[\ks ^\sharp(\fg \xi)]) = \Divxi[\f Q \, \sym[\ks(\f x)] \, \f Q^T] = \f Q \, \Divx(\sym[\ks(\f x)]) \,,
\end{align}
and also concerning its skew-symmetric part
\begin{align}
\Divxi(\skew[\ks ^\sharp(\fg \xi)]) = \Divxi[\f Q \, \skew[\ks(\f x)] \, \f Q^T] = \f Q \, \Divx(\skew[\ks(\f x)]) \,.
\end{align}
In the isotropic indeterminate couple stress theory the nonlocal force stress is given by $\widetilde{\fg \tau}= - \frac{1}{2} \anti \Div \f m$ and
adds to the local force stress yielding the total stress $\widetilde{\fg \sigma} = \fg \sigma + \widetilde{\fg \tau}$. Since $\f m \sim \ks$ we find with help of eq.\eqref{antiasharp} that its transformation reads for the symmetric part
\begin{align}
\anti\{\Divxi(\sym[\ks ^\sharp(\fg \xi)])\} = \anti \{\f Q \,\Divx \, ( \sym[\ks(\f x)]) \} = \f Q \, \anti\{ \Divx(\sym[\ks(\f x)]) \} \, \f Q^T \,,
\end{align}
and similarly for its skew-symmetric part
\begin{align}
\anti\{\Divxi(\skew[\ks ^\sharp(\fg \xi)])\} = \anti \{\f Q \,\Divx \, ( \skew[\ks(\f x)]) \} = \f Q \, \anti\{ \Divx(\skew[\ks(\f x)]) \} \, \f Q^T \,.
\end{align}
Thus, the divergence of the nonlocal stress contribution transforms like the local stress via
\begin{align}
\Divxi[\anti\{\Divxi(\sym[\ks ^\sharp(\fg \xi)])\}] &= \Divxi [\f Q \,\anti \{\Divx \, ( \sym[\ks(\f x)]) \}\,\f Q^T ] \notag \\
&= \f Q \,\Divx[\anti\{ \Divx(\sym[\ks(\f x)]) \}] \,,
\end{align}
and
\begin{align}
\Divxi[\anti\{\Divxi(\skew[\ks ^\sharp(\fg \xi)])\}] &= \Divxi [\f Q \,\anti \{\Divx \, ( \skew[\ks(\f x)]) \}\,\f Q^T ] \notag \\
& = \f Q \,\Divx[\anti\{ \Divx(\skew[\ks(\f x)]) \}] \,.
\end{align}
This can also be seen if we directly write the (extended) balance of linear momentum of the indeterminate couple stress model via
\begin{equation}\label{BLM_ICST}
 {\begin{array}{lrcl}
& \Divxi \{ \sigmaxi - \12 \, \anti \Divxi [ \f m^{\sharp}(\xi)] \} + \f f^{\sharp}(\fg \xi) &=& 0 \vspace{0.2cm} \\
\Leftrightarrow \quad & \Divxi \{ \f Q \, \fg \sigma(\f Q^T \, \fg \xi) \, \f Q^T - \12 \, \anti \Divxi [\f Q \,  \f m(\f x) \, \f Q^T] \} + \f Q \, \f f(\f Q^T \, \fg \xi) &=& 0 \vspace{0.2cm} \\
\Leftrightarrow \quad & \Divxi \{ \f Q \, \sigmax \, \f Q^T - \12 \, \anti (\f Q \,  \Divx [ \f m(\f x) ] ) \} + \f Q \,  \f f(\f x) &= & 0 \vspace{0.2cm} \\
\Leftrightarrow \quad & \Divxi \{ \f Q \, \sigmax \, \f Q^T - \12 \, \f Q \, ( \anti \Divx [\f m(\f x)] ) \, \f Q^T \} + \f Q \,  \f f(\f x) &= & 0 \vspace{0.2cm} \\
\Leftrightarrow \quad & \f Q \, \Divx \{ \sigmax - \12 \, \anti \Divx [\f m(\f x)] \} + \f Q \,  \f f(\f x) &=& 0 \vspace{0.2cm} \\
\Leftrightarrow \quad & \Divx \{ \sigmax - \12 \, \anti \Divx [\f m(\f x)] \} + \f f(\f x) &=& 0 \,.
 \end{array}}
\end{equation}
Therefore, the balance equations of the indeterminate couple stress model are form-invariant with respect to the $\sharp$-transformation.\\

Alternatively, the hyperstress tensor $\m$ of {\bf third order} may appear in higher gradient theories  via \newline
${\m^\sharp(\fg \xi) = \L.\GRADxi{\Gradxi{\uxi}} = \L.{\rm D}_{\fg \xi}^2 \uxi}$, where $\L$ is a 6th order tensor. We only consider $\L$ to map like the identity $\L. \n = \n$ for any $\n \in \Skalar^{3 \times 3 \times 3}$ for simplicity. Thus, we obtain from eq.\eqref{divxidivxigradxigradxiuxi} the transformation
\begin{align}
\Divxi \DIVxi[\m^\sharp(\fg \xi)]  & = \Divxi \DIVxi \{\GRADxi{\Gradxi{\uxi}}\} = \f Q \, \Divx (\Gradx{\Divx\{\Gradx{\ux}\}} ) \notag \\
& =\f Q \, \Divx \DIVx \underbrace{\{\GRADx{\Gradx{\ux}}\}}_{\displaystyle = \m(\f x)} \notag \\
& = \f Q \, \Divx\{\DIVx[\m(\f x)]\} \,.
\end{align}
Hence the balance of linear momentum in a full higher gradient theory transforms according to
\begin{equation}\label{BLM2}
 {\begin{array}{lrcl}
& \Divxi[ \sigmaxi + \DIVxi [ \m^{\sharp}(\fg \xi)] ] + \f f^{\sharp}(\fg \xi) &=& 0 \vspace{0.2cm} \\
\Leftrightarrow \quad & \Divxi[\sigmaxi] + \f f^{\sharp}(\fg \xi) +  \Divxi \, \DIVxi[\m^\sharp(\fg \xi)] &=& 0 \vspace{0.2cm} \\
\Leftrightarrow \quad & \Divxi[\f Q \, \fg \sigma(\f Q^T \, \fg \xi) \, \f Q^T] + \f Q \, \f f(\f Q^T \, \fg \xi) +  \Divxi \, \DIVxi[\m^\sharp(\fg \xi)] &=& 0 \vspace{0.2cm} \\
\Leftrightarrow \quad & \f Q \, \Divx[\fg \sigma(\f x)] + \f Q \, \f f(\f x) + \f Q \, \{ \Divx \, \DIVx[\m(\f x)]  \} & = & 0 \vspace{0.2cm} \\
\Leftrightarrow \quad & \Divx[\fg \sigma(\f x) +\DIVx \m(\f x) ] + \f f(\f x)  & = & 0 \,.
 \end{array}}
\end{equation}
Therefore, the balance equations of the higher gradient model are form-invariant with respect to the $\sharp$-transformation.\\

Finally, we show that the form invariance condition for the system of eq.\eqref{strongform} in a variant of the indeterminate couple stress theory \cite{GhiNeMaMue15} always holds true:\footnote{In the physics oriented literature, {\bf form-invariance} of objects under certain coordinate transformations is called {\bf covariance}. Obviously, we discuss here covariance of the balance equations under simultaneous rotations of the referential and spatial coordinates, c.f. \cite[p.~156]{Marsden83}.}
\begin{equation}\label{BLM1}
 {\begin{array}{lrcl}
& \Divxi[ \sigmaxi + \sym \Curlxi(\f m^{\sharp}(\fg \xi)) ] + \f f^{\sharp}(\fg \xi) &=& 0 \vspace{0.2cm} \\
\Leftrightarrow \quad & \Divxi[\f Q \, \fg \sigma(\f Q^T \, \fg \xi) \, \f Q^T + \sym \Curlxi(\f Q \, \f m(\f Q^T \, \fg \xi) \, \f Q^T) ] + \f Q \, \f f(\f Q^T \, \fg \xi) &=& 0 \vspace{0.2cm} \\
\Leftrightarrow \quad & \Divxi[\f Q \, \sigmax \, \f Q^T + \sym \{ \f Q \, \Curlx( \f m(\f x)) \, \f Q^T \} ] + \f Q \,  \f f(\f x) &= & 0 \vspace{0.2cm} \\
\Leftrightarrow \quad & \Divxi[\f Q \, \sigmax \, \f Q^T +  \f Q \, \sym \Curlx( \f m(\f x))  \, \f Q^T ] + \f Q \,  \f f(\f x) &=& 0 \vspace{0.2cm} \\
\Leftrightarrow \quad & \Divxi[\f Q \, \{ \sigmax + \sym \Curlx( \f m(\f x)) \} \, \f Q^T ] + \f Q \,  \f f(\f x) &=& 0\vspace{0.2cm} \\
\Leftrightarrow \quad & \f Q \, \Divx[\sigmax + \sym \Curlx( \f m(\f x))] + \f Q \,  \f f(\f x) &=& 0 \vspace{0.2cm} \\
\Leftrightarrow \quad & \f Q \, \{ \Divx[\sigmax + \sym \Curlx( \f m(\f x))] + \f f(\f x) \} &=& 0 \vspace{0.2cm} \\
\Leftrightarrow \quad & \Divx[\sigmax + \sym \Curlx( \f m(\f x))] + \f f(\f x) &=& 0 \,.
 \end{array}}
\end{equation}
\section{Objectivity and isotropy in linear higher gradient elasticity}\label{KapObjectIso}
\setcounter{equation}{0}
Now we turn back to the invariance considerations in Sections \ref{KapObjNonlinearElas} and \ref{KapObjLinearElas}. Motivated by the insights obtained there, we require {\bf form-invariance} of the {\bf energy expression} also in higher gradient or couple-stress models. As a primer, consider the higher order curvature expression
\begin{align}\label{wcurv}
W_{\rm curv} = \mu \, L_c^2 \, \| \sym \Curl (\sym \Grad{\f u}) \|^2 \, .
\end{align}
Linearized objectivity in eq.\eqref{wcurv} is clear, since the  curvature is defined in terms of linear combinations of partial derivatives of $\fg \varepsilon = \sym \Grad{\f u}$. Moreover, since
\begin{align}\label{CurlxisymGradxi}
\Curlxi(\sym\{\Gradxi{\uxi}\}) &= \Gradxi{\axl \skew \Gradxi{\uxi}} = \Gradxi{\axl \skew (\f Q \, \Gradx{\ux} \, \f Q^T)} \notag \\ & =  \Gradxi{\axl (\f Q \, ( \skew \Gradx{\ux}) \, \f Q^T) }=  \Gradxi{\f Q \, \underbrace{( \axl \skew \Gradx{\ux})}_{\displaystyle \colonequals \f v(\f x)}} \notag \\
&= \Gradxi{\f Q \, \f v(\f Q^T \, \fg \xi)} = \f Q \, \Gradx{\f v(\f x)} \, \f Q^T = \f Q \, (\Gradx{\axl \skew \Gradx{\ux}}) \notag \\
&= \f Q \, \Curlx(\sym \Gradx{\ux}) \, \f Q^T \,,
\end{align}
%
%
%
eq.\eqref{kssharpxi} yields
\begin{align}
\Curlxi(\sym\{\Gradxi{\uxi}\}) =(\ks^{\sharp} (\fg \xi))^T = \f Q \, (\ks (\f x))^T \, \f Q^T\,.
\end{align}
Thus, the curvature energy in eq.\eqref{wcurv} is form-invariant since
\begin{align}
\|\sym\Curlxi(\sym\{\Gradxi{\uxi}\}) \|^2 &=\| \sym(\ks^{\sharp} (\fg \xi) \|^2 = \| \sym  \f Q \, \ks (\f x)  \, \f Q^T \|^2 = \| \sym \ks(\f x) \|^2 \notag \\
&=  \| \sym \Curlx \sym \Gradx{\f u(\f x)} \|^2 \,.
\end{align}
\subsection{Energetic considerations}\label{KapEnerConsider}
Here, we consider (linearized) frame-indifferent second gradient elasticity written in the reduced format\footnote{Here, it must be noted that $W(\sym \Grad{u} , \, \Grad{\Grad{u}})$ is also (linearized) frame-indifferent. The reduction to $W(\sym \Grad{u} , \, \Grad{\sym \Grad{u}})$ is possible since any second derivative $\Grad{\Grad{u}}$ is a linear combination of gradients of strain $\sym \Grad{u}$. However, there remains a subtle difference between $W(\sym \Grad{u} , \, \Grad{\Grad{u}})$ and $W(\sym \Grad{u} , \, \Grad{\sym \Grad{u}})$ in that formulations based on either expression may differ in the higher order traction boundary conditions while the energies coincide, see subsection \ref{KapProblemCompInv}.}
\begin{align}
W(\Grad{u} , \, \GRAD{\Grad{u}}) = W(\sym \Grad{u} , \, \GRAD{\sym \Grad{u}}) \,,
\end{align}
where $W$ is a quadratic form. Assuming the natural additive split
\begin{align}\label{WSplit}
W(\sym \Grad{\f u} , \, \GRAD{\sym \Grad{\f u}}) = W_{\rm lin}(\sym \Grad{\f u}) + W_{\rm curv}(\GRAD{\sym \Grad{\f u}}) \,,
\end{align}
we are interested in further reduced forms for isotropy. The reduced format for $W_{\rm lin}$ is nothing but isotropic linear elasticity, in which isotropy is in effect the requirement of form-invariance of
\begin{align}\label{WlinGradxi}
W_{\rm lin}(\sym \Gradxi{\uxi}) \overbrace{=}^{\rm form-invariant} W_{\rm lin}(\sym \Gradx{\ux}) \,,
\end{align}
under the transformation \eqref{usharpdef}. From eq.\eqref{WlinGradxi} we obtain the standard condition
\begin{align}\label{WlinQT}
W_{\rm lin}(\f Q \, \sym \Grad{\f u} \, \f Q^T) = W_{\rm lin}(\sym \Grad{\f u}) \qquad \forall \, \f Q \in \SO(3) \,,
\end{align}
which implies, as is well-known, the representation
\begin{align}\label{Wlinsymgradu}
W_{\rm lin}(\sym \Grad{\f u})  &=  \mu \, \| \sym \Grad{\f u} \|^2 + \frac{\lambda}{2} \, [\tr \sym \Grad{\f u} ]^2 \notag\\
&=  \mu \, \|\dev \sym \Grad{\f u} \|^2 + \frac{\kappa}{2} \, [\tr \sym \Grad{\f u} ]^2 \,,
\end{align}
where $\mu$, $\lambda$ are the classical Lam\'{e}-constants, and $\kappa$ is the bulk modulus.\\

\noindent
For the general curvature term in eq.\eqref{WSplit} we require for isotropy the same formal form-invariance as in eq.\eqref{WlinGradxi}, namely
\begin{align}\label{WcurvGradxi}
W_{\rm curv}(\GRADxi{\sym \Gradxi{\uxi}}) \overbrace{=}^{\rm form-invariant} W_{\rm curv}(\GRADx{\sym \Gradx{\ux}}) \,.
\end{align}
From eq.\eqref{symGraduxi} and \eqref{gradxigradxiusharp} we obtain
\begin{align}\label{GradsymGradxi}
\GRADxi{\sym \Gradxi{\uxi}} = \GRADx{\sym \f Q \, \Gradx{\ux} \, \f Q^T} \, \f Q^T \, ,
\end{align}
yielding
\begin{align}\label{WcurvGradxi2}
 W_{\rm curv}(\GRADx{\sym \f Q \, \Gradx{\ux} \, \f Q^T} \, \f Q^T) = W_{\rm curv}(\GRADx{\sym \Gradx{\ux}})
\end{align}
in this general case of a strain gradient model.\\

\noindent
In the case of the indeterminate couple stress model with symmetric total stress tensor we have that
\begin{align}\label{WcurvCurl}
W_{\rm curv}(\GRADx{\sym \Gradx{\ux}}) = \widehat W_{\rm curv}\left( (\Curlx{\sym \Gradx{\ux}})^T \right)  \,,
\end{align}
and we obtain similarly the form-invariance condition
\begin{align}\label{WcurvCurl2}
\widehat W_{\rm curv}((\Curlxi{\, \sym \Gradxi{\uxi} })^T ) \overbrace{=}^{\rm form-invariant} \widehat W_{\rm curv}((\Curlx{\sym \Gradx{\ux}})^T) \,,
\end{align}
which implies
\begin{align}\label{WcurvCurl3}
\widehat W_{\rm curv}(\f Q \,  (\Curlx{\, \sym \Gradx{\ux} })^T \, \f Q^T ) = \widehat W_{\rm curv}((\Curlx{\, \sym \Gradx{\ux}})^T) \,.
\end{align}
If this is true for all $\f Q \in \SO(3)$ as in \eqref{Wlinsymgradu} it can be concluded that
\begin{align}\label{WcurvCurl4}
\widehat W_{\rm curv}(\underbrace{(\Curlx{\, \sym \Gradx{\ux} })^T}_{\displaystyle \ks(\f x)}) = \mu \, L_c^2 \, \{ \alpha_1 \norm{\dev \sym \ks(\f x) }^2 + \alpha_2 \, \norm{\skew \ks(\f x)}^2 + \alpha_3 \, [\tr \ks(\f x)]^2 \}
\end{align}
with nondimensional constants $\alpha_1$, $\alpha_2$, $\alpha_3$ and we note that $\tr \ks(\f x) = 0$.

\begin{align}\label{LinHGFormuBox}
\fbox{
\parbox[][5cm][c]{14.5cm}{
The linear second-gradient hyperelastic formulation is (linearized) frame-indifferent and isotropic {\bf if and only if} there exists a function $\widehat{\Psi}_{\rm lin} : \Sym(3) \times \Skalar^{18} \rightarrow \Skalar$ such that\footnotemark\\
\\
$
W_{\rm lin}(\Grad{u} , \, \GRAD{\Grad{u}}) \!\!\!\!\!\!\!\! \overbrace{=}^{\text{linearized frame-indifferent} \atop \text{ reduced format}} \!\!\!\!\!\!\!\! \widehat{\Psi}_{\rm lin}(\sym \Grad{\f u} \, , \, \GRAD{\sym \Grad{\f u}}) \,, \quad {\rm and}
$
\\
\\
$
W_{\rm lin}(\Gradxi{\uxi} \, , \, \GRADxi{\Gradxi{\uxi}}) \overbrace{=}^\text{ form-invariance} W_{\rm lin}(\Gradx{\ux} \, , \, \GRADx{\Gradx{\ux}})
$
\\
\\
$
[ {\rm i.e.} \quad \widehat{\Psi}_{\rm lin}(\sym \Gradxi{\uxi} \, , \, \GRADxi{\sym \Gradxi{\uxi}})
$
\center{
$
= \widehat{\Psi}_{\rm lin}(\sym \Gradx{\ux} \, , \, \GRADx{\sym \Gradx{\ux}}) \quad \forall \, \f Q \in \SO(3)] \,.
$}
}}
\end{align}\footnotetext{A remark on the {\bf only if}: also energies like $\norm{\Grad{\axl(\skew \Grad{u})}}^2$ can be subsumed, since from compatibility the former can be written as $\norm{\Curl \sym \Grad{u}}^2 = \widehat{\Psi}_{\rm curv}(\Grad{\sym \Grad{u}})$. As seen in Sect.\ref{KapProblemCompInv} this does not mean that all modelling variants of isotropic and (linearized) frame-indifferent formulations are exhausted using the representation with $\widehat{\Psi}_{\rm lin}(\sym \Gradx{\ux} \, , \, \GRADx{\sym \Gradx{\ux}})$.}
In terms of $\widehat{\Psi}_{\rm lin}$ a second gradient material is said to be isotropic if the stored energy satisfies for all $Q \in \SO(3)$:
\begin{align}\label{PsiStoredEnergy}
\widehat{\Psi}_{\rm lin}&(\sym \Gradx{\ux} \, , \, \GRADx{\sym \Gradx{\ux}}) \notag \\
&= \widehat{\Psi}_{\rm lin} \bigg( Q \, (\sym \Gradx{\ux}) \, Q^T \, , \, \underbrace{\GRADx{ Q \, (\sym \Gradx{\ux}) \, Q^T} \, Q^T }_{\displaystyle (Q_{ib} \, \12 [u_{b,c}+u_{c,b} ] \, Q^T_{cj})_{,a} \, Q^T_{ak} \, \f e_i \otimes \f e_j \otimes \f e_k\footnotemark}\bigg)
\end{align}\footnotetext{Following Truesdell \cite[p.~60]{Truesdell66} the rotation $Q \in \SO(3)$ can only be constant and therefore\\
${(Q_{ib} \, \12 [u_{b,c}+u_{c,b} ] \, Q^T_{cj})_{,a} \, Q^T_{ak} \, \f e_i \otimes \f e_j \otimes \f e_k = (Q_{ib} \, \12 [u_{b,ca}+u_{c,ba} ] \, Q^T_{cj}) \, Q^T_{ak} \, \f e_i \otimes \f e_j \otimes \f e_k }$. Note again that any coordinate transformation $\zeta: \, \Euklid \mapsto \Euklid$ whose gradient is a rotation every where, i.e. $\Gradx{\zeta(x)} = Q(x)$ is necessarily of the affine form $\zeta(x)=\overline{Q} \, x + \bar{b}$, with constants $\bar{b} \in \Euklid$ and $\overline{Q} \in \SO(3)$, see \cite{Neff_curl06}.}and this is an immediate consequence of \eqref{LinHGFormuBox}. Considering condition \eqref{PsiStoredEnergy} in the case of the theory without higher gradients we have
\begin{align}\label{PsiLinWithoutHigh}
{\Psi}_{\rm lin}&(\sym \Gradx{\ux}) = {\Psi}_{\rm lin} \left( Q \, (\sym \Gradx{\ux}) \, Q^T \right) \quad \forall \f Q \in \SO(3)
\end{align}
and we may probe \eqref{PsiLinWithoutHigh} also with non-constant rotation fields $\f Q = \f Q(\f x) \, \in \SO(3)$ without altering the imposed invariance condition. However, it is important to realize that the adopted invariance requirement for isotropy in higher-gradient materials via the $\sharp$-transformation \eqref{LinHGFormuBox} does not allow for space dependent rotation fields and the invariance condition for isotropy {\bf does not} read
\begin{align}\label{PsiNOT}
\widehat{\Psi}_{\rm lin}&(\sym \Gradx{\ux} \, , \, \GRADx{\sym \Gradx{\ux}}) \notag \\
&= \widehat{\Psi}_{\rm lin} \left( Q(\f x) \, (\sym \Gradx{\ux}) \, Q^T(\f x) \, , \, \GRADx{ Q(\f x) \, (\sym \Gradx{\ux}) \, Q^T(\f x)} \, Q^T(\f x) \right) \,.
\end{align}
\section{The isotropic linear indeterminate couple stress theory}\label{KapIndeterminate}
\setcounter{equation}{0}
This section does not contain new results but is included to connect our investigation of isotropy conditions to specific models proposed in the literature. The linear indeterminate couple stress model is a specific second gradient elastic model, in which the higher order interaction is restricted to the continuum rotation $\curl \, \f u$,
where $\f u : \Omega \mapsto \R^3$ is the displacement of the body.
It is therefore interpreted  to be sensitive to rotations of material points and it is
possible to prescribe boundary conditions of rotational type.

Superficially, this is the simplest possible generalization of linear elasticity in order
to include the gradient of the local continuum rotation as a source of additional strains
and stress with an associated energy. In this paper, therefore, we limit our analysis to isotropic
materials and only to the second gradient of the displacement:
\begin{align}\label{D2u}
\GRAD{\Grad{u}} \sim {\rm D}^2_{\f x} \f u=\underbrace{\frac{\partial^2 u_i}{\partial x_j \, \partial
x_k}}_{\displaystyle u_{i,jk}} \, \f e_i \otimes \f e_j \otimes \f e_k
=(\varepsilon_{ji,k}+\varepsilon_{ki,j}-\varepsilon_{jk,i}) \,\f e_i \otimes \f e_j \otimes
\f e_k \, ,
 \end{align}
where
\begin{align}\label{strain}
\fg \varepsilon= \varepsilon_{ij} \, \f e_i \otimes \f e_j = \12 (u_{i,j} +
u_{j,i})\, \f e_i \otimes \f e_j = \12 (\text{Grad} [\f u ]+ (\text{Grad} [\f u])^T ) = \sym
\text{Grad}[ \f u]
 \end{align}
is the symmetric linear strain tensor. Thus, from eq.\eqref{D2u} all second derivatives of the displacement field $\f u$ can be obtained from linear combinations of derivatives of $\fg \varepsilon$. In general, strain gradient models do not introduce additional independent degrees of
freedoms\footnote{In contrast to Cosserat models where an independent rotation field is
under consideration.} aside the displacement field $\f u$. Thus, the higher derivatives
introduce a ``latent-microstructure" (constraint microstructure \cite{Grioli03}).
However, this apparent simplicity has to be payed with more complicated and intransparent boundary
conditions, as treated in a series of papers \cite{MaGhiNeMue15,NeGhiMaMue15}.\\

\noindent
We assume the elastic energy to be given by
\begin{align}\label{Welastic}
I(\f u)&= \int_\Omega W_{\rm lin}(\fg \varepsilon) + W_{\rm curv}(\ks) \, \di V \notag\\
&= \int_\Omega \mu \norm{\fg \varepsilon}^2 +
\frac{\lambda}{2} [\tr(\fg \varepsilon)]^2 + \mu\,L_c^2\,(\alpha_1\, \|\dev\sym
[\widetilde{\f k}]\|^2+ \alpha_2 \, [\tr(\widetilde{\f k})]^2 + \alpha_3\,
\|\skew[\widetilde{\f k}]\|^2) \, \di V \, ,
\end{align}
where $\mu$ and $\lambda$ are the constitutive Lam\'{e} constants and the curvature energy is
based on the second order curvature tensor $\ks$
\begin{align}\label{kcurlu}
\widetilde{\f k}=\Grad{ \axl(\skew \Grad{ \f u})}=\12 \, \Grad{ \curl \, \f u} \, ,
\end{align}
with additional dimensionless constitutive parameters $\alpha_1$, $\alpha_2$, $\alpha_3$, and the characteristic length $L_c > 0$. Taking free variations $\delta \f u\in C_0^\infty(\Omega,\Gamma)$ of the elastic energy $I(\f u)$ yields the virtual work principle
\begin{align}\label{gradeq211}
\frac{\rm d}{\rm dt}\mid_{t=0} \, I(\f u + t \, \delta \f u)=&\int_\Omega 2\mu\,\langle \fg \varepsilon, \Grad{ \delta \f u} \rangle+\lambda \tr(\fg \varepsilon)\,\tr( \Grad{ \delta \f u}) \notag\\
&+ 2 \, \mu \, L_c^2 \, \alpha_1\, \langle \dev \sym (\ks),\dev\sym \Grad{\axl \, \skew \Grad{ \delta \f u}} \rangle \notag \\
& + 2 \, \mu \, L_c^2 \, \alpha_2\, \tr(\ks) \, \tr(\Grad{\axl \, \skew \Grad{ \delta \f u}}) \notag \\
& + 2 \, \mu \, L_c^2 \, \alpha_3\, \langle \skew(\ks) ,
\skew(\Grad{ \axl\,\skew \, \Grad{ \delta \f u}})\rangle + \langle \f f , \delta \f u\rangle \di
V =0 \, .
\end{align}
Using the classical divergence theorem for the curvature term in eq.\eqref{gradeq211} it
follows after some simple algebra that
\begin{align}\label{germaneq311}
\int_\Omega \langle \Div ({\fg \sigma}+\widetilde{\fg \tau})+\f f, \delta \f u \rangle \, \di
V-\int_{\partial \Omega}\langle ({\fg \sigma}+\widetilde{\fg \tau}).\,\f  n, \delta \f
u\rangle \, \di A + \int_{\partial\Omega}&\langle  {\f m}.\, \f n, \axl\skew \Grad{ \delta \f
u} \rangle \di A=0,
\end{align}
where  $\fg \sigma$ is the symmetric local force-stress tensor
\begin{align}
\fg \sigma=2\, \mu \, \fg \varepsilon+\lambda \, \tr( \fg \varepsilon)\, \id \quad \in
{\rm Sym}(3) \,,
\end{align}
and $\widetilde{\fg \tau}$ represents  the nonlocal force-stress  tensor
\begin{align}\label{nonlocalForcestress}
 \widetilde{\fg \tau}&= - \frac{1}{2}\anti {\rm Div}[ {\f m}] \quad \in \so(3) \,,
\end{align}
which here is automatically skew-symmetric\footnote{Lederer and Khatibi refer in  \cite{Lederer15b} to  equation \eqref{nonlocalForcestress}, which is properly invariant under the $\sharp$-transformation: Since $\widetilde{\fg \tau}$ and $\f m$ are tensors of second order it holds $\widetilde{\fg \tau}^{\sharp}(\fg \xi) = \f Q \, \widetilde{\fg \tau}(\f x) \, \f Q^T$ and $\f m^{\sharp}(\fg \xi) = \f Q \, \f m(\f x) \, \f Q^T$, respectively. Further, eq.\eqref{divxisigmaxi} yields $\Divxi [\f m^{\sharp}(\fg \xi)] = \f Q \, \Divx [\f m(\f x)]$ and with help of eq.\eqref{antiasharp} we obtain $\anti(\Divxi[\f m^{\sharp}(\fg \xi)]) = \anti (\f Q \, \Divx [\f m(\f x)]) = \f Q \, \anti(\Divx [\f m(\f x)]) \, \f Q^T$.}.\\

\noindent
The second order couple stress\footnote{Also denominated hyperstress or moment stress.} tensor $\f m$  in
eq.\eqref{nonlocalForcestress} reads
\begin{align}\label{DefCoupleStress}
{\f m}&=\mu \, L_c^2 \, [{\alpha_1}\dev \sym (\Grad{ \curl \, \f u}) + \alpha_2 \,
\underbrace{\tr(\Grad{ \curl \, \f u})}_{\displaystyle = \,0} \, \id
+ {\alpha_3}\,\skew(\Grad{ \curl \, \f u})] \notag \\
&=2 \, \mu \, L_c^2\,[\alpha_1 \, \dev \sym (\ks) + \alpha_3
\, \skew(\ks)] \,,
\end{align}
which may or may not be symmetric, depending on the material parameters $\alpha_1$, $\alpha_3$. However, $\f m$ in eq.\eqref{DefCoupleStress} is automatically trace free since both the deviator and the skew
operator yield trace free tensors.

However, the asymmetry of the nonlocal force stress $\widetilde{\fg \tau}$ appears as a constitutive assumption.
Thus, if the test function $\delta \f u\in C_0^\infty(\B,\Gamma)$ also satisfies
$\axl(\skew \Grad{ \delta \f u})=0$ on $\Gamma$ (equivalently $\curl \, \delta \f u=0$), then
we obtain the balance of momentum
\begin{align}\label{Div2muvarepsilon}
\Div \bigg\{&\underbrace{2\, \mu \, \fg \varepsilon+\lambda \, \tr(\fg \varepsilon)\,
\id}_{\displaystyle \text{local force stress}\  \sigma \, \in \, {\rm Sym}(3)}
\underbrace{- \anti {\rm Div}\{\underbrace{\mu\,L_c^2\,\alpha_1\, \dev\sym
(\ks)+ \mu\,L_c^2\,\alpha_2\,\skew (\ks)}_{\displaystyle \text{\rm hyperstress } \tfrac{1}{2}\f m \, \in \, \gl(3)}\}}_{\displaystyle \text{\rm
completely skew-symmetric non-local force stress}\ \widetilde{\fg \tau} \, \in \,
\so(3)}\bigg\}+ \, \f f=0 \, .
\end{align}
Since the balance of angular momentum is given by eq.\eqref{nonlocalForcestress}, we can combine them to see
compact equilibrium equation

\begin{align}\label{CoupleStressBalanceMomentum}
\fbox{
\parbox[][2cm][c]{12cm}{$
\Div \widetilde{\fg \sigma} + \f f = 0 \, , \quad \quad \widetilde{\fg
\sigma} = \fg \sigma + \widetilde{\fg \tau} \, \notin \Sym(3) \, , \quad \text{total force stress} \\
\\
\fg \sigma = 2 \, \mu \, \sym \Grad{ \f u} + \lambda \tr(\Grad{ \f u}) \id \, ,
\quad \quad \widetilde{\fg \tau} = - \frac{1}{2} \anti \Div \f m \in \so(3) \\
\\
 \Div{\f m} + 2 \axl(\widetilde{\fg \tau}) = 0  \quad \Leftrightarrow \quad
 \Div{\f m} + 2 \axl(\skew \widetilde{\fg \sigma}) = 0 \, .
 $}}
\end{align}
\subsection{Related models in isotropic second gradient elasticity}
Let us consider the strain and curvature energy as minimization problem
\begin{align}
W(\f u)=\int_\B \left[\mu\, \|{\rm sym} \, \Grad{ \f u}\|^2+\frac{\lambda}{2}\, [\tr({\rm
sym} \, \Grad{ \f u})]^2+W_{\rm curv}({\rm D}^2_{\f x} \f u)\right] \di V \quad \mapsto \quad
\text{min. w.r.t.} \, \f u,
\end{align}
admitting unique minimizers under some appropriate boundary condition. Here $\lambda,\mu$
are the Lam\'{e} constitutive coefficients of isotropic linear elasticity, which is
fundamental to small deformation gradient elasticity. If the curvature energy has the
form $W_{\rm curv}({\rm D}^2_{\f x} \f u)=W_{\rm curv}({\rm D}_{\f x} \sym \Gradx{\f u})$, the model
is called {\bf a strain gradient model}. We define the hyperstress tensor of third order
as $\m={\rm D}_{{\rm D}^2_{\f x} \f u}W_{\rm curv}({\rm D}^2_{\f x} \f u)$. Note that $\Div \m$ is in general not symmetric.

In the following we recall  some curvature energies proposed in different isotropic second gradient
elasticity models:
\begin{itemize}
\item {\bf the indeterminate couple stress model} (Grioli-Koiter-Mindlin-Toupin model)
\cite{Grioli60,Aero61,Koiter64,Mindlin62,Toupin64,Sokolowski72,Grioli03}
in which the higher derivatives (apparently)  appear only through derivatives of the infinitesimal
continuum rotation $\curl \, \f u$.  Hence, the curvature energy  has  the equivalent forms
\begin{align}\label{KMTe}\notag
W_{\rm curv}(\ks)&=\frac{\mu \,L_c^2}{4}(\alpha_1 \, \|\sym \Grad{ \curl\, \f u}\|^2+\alpha_3\,\| \skew \Grad{ \curl\, \f u}\|^2\notag\\
&=\mu\,L_c^2(\alpha_1\, \|\sym\underbrace{\Grad{\axl(\skew \Grad{ \f u})}}_{\displaystyle \ks}\|^2+\,\alpha_3\,\| \skew \underbrace{ \Grad{\axl(\skew
\Grad{ \f u})}}_{\displaystyle \ks}\|^2 \, , \\
\f m &= 2\,\mu\,L_c^2(\alpha_1\,\sym\ks+\alpha_3\,\skew\ks)  \, .
\end{align}
Although  this energy admits  the equivalent forms \eqref{KMTe}$_1$ and \eqref{KMTe}$_2$,
the equations and the boundary value problem of the indeterminate couple stress model is
usually formulated only using the form \eqref{KMTe} of the energy. We remark that the
spherical part of the couple stress tensor is zero since $\tr(2\,\ks)=\tr(\nabla \curl \, \f u)={\rm
div} (\curl \, \f u)=0$ as seen before. In order to ensure the pointwise uniform positive definiteness it is
assumed that $\alpha_1>0, \alpha_3>0$. Note that pointwise uniform positivity is often assumed
\cite{Koiter64} when deriving analytical solutions for simple boundary value problems
because it allows to invert the couple stress-curvature relation.
\item
  {\bf the modified symmetric couple stress model - the conformal model}.  On the other hand, in the  conformal case  \cite{Neff_Jeong_IJSS09,Neff_Paris_Maugin09} one may consider $\alpha_3=0$, which makes  the couple stress tensor ${\f m}$ symmetric and trace free \cite{dahler1963theory}.  This conformal  curvature case has been derived by Neff in  \cite{Neff_Jeong_IJSS09}, the curvature energy having the form
\begin{align}
W_{\rm curv}(\ks)&=\mu \, L_c^2\, \alpha_1\, \|\dev \sym \ks\|^2 , \qquad \f m=2\,\mu\,L_c^2\,\alpha_1\,\dev \sym\,\ks\,.
\end{align}
Indeed, there are two major reasons uncovered in \cite{Neff_Jeong_IJSS09} for using the
modified couple stress model. First, in order to avoid non physical singular stiffening behaviour for
smaller and smaller samples in bending \cite{Neff_Jeong_bounded_stiffness09} one has to
take $\alpha_3=0$. Second, based on a homogenization procedure invoking a natural
``micro-randomness" assumption (a strong statement of microstructural isotropy) implies
conformal invariance, which is again $\alpha_3=0$. Such a modell is still well-posed
\cite{Neff_JeongMMS08} leading to existence and uniqueness results with only one
additional material length scale parameter, while it is {\bf not} pointwise uniformly
positive definite.
\item {\bf the skew-symmetric couple stress model}.
  { Hadjesfandiari and Dargush} strongly advocate
  \cite{hadjesfandiari2011couple,hadjesfandiari2013fundamental,hadjesfandiari2013skew} the opposite
  extreme case, $\alpha_1=0$ and $\alpha_3>0$, i.e. they used the  curvature  energy
\begin{align}
W_{\rm curv}(\ks)&=\mu\,L_c^2\,\frac{\alpha_3}{4}\, \|\skew \Grad{{\rm curl}\, \f u}\|^2=\mu\,L_c^2\,\alpha_3\, \|\skew \ks\|^2\, , \quad \quad \f m=2\,\mu\,L_c^2\,\alpha_3\,\skew\,\ks\,.
\end{align}
In that model the nonlocal force stress tesnor $\widetilde{\fg \tau}$ is skew-symmetric and the couple stress tensor $\f m$ is assumed to be completely skew-symmetric. Their reasoning, based on boundary conditions, is critically discussed in Neff et al.\cite{NeMueGhiMa15}.
 \end{itemize}
\subsection{A variant of the indeterminate couple stress model with symmetric total force stress}
In Ghiba et al.\cite{GhiNeMaMue15} the isotropic, linear indeterminate couple stress
model has been modified so as to have symmetric total force stress $\widehat{\fg \sigma}$, while retaining the
same weak form of the Euler-Lagrange equations. This is possible since the force stresses
appearing in the balance of forces is only determined up to a self-equilibrated
stress-field $\bar{\fg \sigma}$, i.e.
\begin{align}\label{selfstress}
\Div \widetilde{\fg \sigma} + \f f = 0 \, \Leftrightarrow \, \Div (\widetilde{\fg \sigma}
+ \bar{\fg \sigma}) + \f f = 0 \, , \quad {\text{for any}} \, \bar{\fg \sigma} \,
{\text{with}} \quad \Div \bar{\fg \sigma}=0 \,.
\end{align}
The curvature energy expression of this new model is
\begin{align}
W_{\rm curv}({\rm D}^2_{\f x} \f u) = \mu \, L_c^2 ( \alpha_1 \, \norm{\dev \sym \Curl(\sym \Grad{ \f
u})}^2 + \alpha_2 \, \norm{\skew \Curl(\sym \Grad{ \f u})}^2 ) \, .
\end{align}
The strong form of the new model reads

\begin{align}\label{strongform}
\fbox{
\parbox[][2cm][c]{12cm}{$
\Div \widehat{\fg \sigma} + \f f = 0 \, , \quad \quad \widehat{\fg
\sigma} = \fg \sigma + \widehat{\fg \tau} \, \in
\Sym(3) \, ,  \quad \text{total force stress} \\
\\
\fg \sigma = 2 \, \mu \, \sym \Grad{ \f u} + \lambda \tr(\Grad{ \f u}) \id \, ,
\quad \quad \widehat{\fg \tau} = \sym \Curl(\widehat{\f m}) \, \in \Sym(3) \\
\\
\widehat{\f m} = 2 \, \mu \, L_c^2 ( \alpha_1 \, \dev \sym \Curl(\sym \Grad{ \f u})+
\alpha_2 \, \skew \Curl(\sym \Grad{ \f u})) \,. $}}
\end{align}
\\
The total force stress tensor is now $\widehat{\fg \sigma}=\fg \sigma + \widehat{\fg
\tau}$ and the second order couple stress tensor is $\widehat{\f m}$. Note that similarly as in the
indeterminate couple stress theory we have $\tr(\widehat{\f m})=0$. Compared to the
indeterminate couple stress theory on can show that $\Div(\widetilde{\fg
\sigma}-\widehat{\fg \sigma})=0$, as claimed. Thus, eq.\eqref{strongform} is a {\sl
couple stress model with symmetric total force stress $\widehat{\fg \sigma}$ and trace
free couple stress tensor $\widehat{\f m}$.} Moreover, the couple stress tensor $\f m$ is
symmetric itself for the possible choice $\alpha_2=0$.
\subsection{The general strain gradient model}
The strong form of the general strain gradient model reads

\begin{align}\label{strongformGSGM}
\fbox{
\parbox[][3cm][c]{14cm}{$
\Div(\fg \sigma + \DIV \m) + \f f = 0 \, , \qquad  \m \, \,  \text{third order hyperstress tensor} \\
\\
\fg \sigma = 2 \, \mu \, \sym \Grad{ \f u} + \lambda \tr(\Grad{ \f u}) \id \, , \quad \breve{\tau}=\DIV \m \, \notin \Sym(3)\\
\\
\breve{\sigma} = \sigma + \breve{\tau} \, \notin \Sym(3) \, , \qquad \Div \breve{\sigma}  + \f f = 0\\
\\
\m = \L.\GRADx{\Gradx{u(x)}} \, \in \Skalar^{3 \times 3 \times 3} , \quad \L: \Skalar^{3 \times 3 \times 3} \rightarrow \Skalar^{3 \times 3 \times 3} , \quad \text{linear 6th order tensor.} $}}
\end{align}
For isotropy, it has been shown in \cite{dellIsola2009} that $\L$ has 5 independent coefficients.
\section{Conclusions and outlook}\label{KapConclusion}
\setcounter{equation}{0}
We have investigated the isotropy requirements in elasticity theory by going back to the basic invariance requirements under transformation of coordinates. In nonlinear elasticity it is important to observe the difference between right-local $\SO(3)$-invariance and right-global $\SO(3)$-invariance. We have given examples that both notions do not coincide in general and have shown another example where even right-local $\SO(3)$-invariance is satisfied. For us, isotropy is right-global $\SO(3)$-invariance while right-local $\SO(3)$-invariance is a stronger condition imposing more restrictive invariance requirements. Next, in linear theories (whether first or second gradient or indeed higher order gradient theories) defined in the linearized frame-indifferent object $\varepsilon = \sym \Grad{\f u}$ and gradients of it, isotropy can be probed by simultaneous rigid rotation of spatial and referential coordinates. We have called this transformation the $\sharp$-mapping. Isotropy then is equivalent to {\bf form-invariance} of all energy terms with respect to the $\sharp$-transformation. Similarly, isotropy is equivalent to form-invariance of the balance equations with respect to the $\sharp$-transformation. We have investigated, in this respect, several variants of higher-gradient models.\\

\noindent
Our approach can be immediately extended to discuss anisotropy for higher gradient theories simply by requiring form-invariance of the energy in the $\sharp$-transformation not for all $\f Q \in \SO(3)$ (isotropy) but only for $\f Q \in \mathcal{G}$ (material symmetry group). This will lead to very succinct statements of anisotropy for higher gradient models.\\

\noindent
%
%
\addcontentsline{toc}{section}{Bibliography}
\footnotesize{
\bibliographystyle{plain}
\bibliography{literaturNeff,literaturGhiba,literaturMuench}
}

\addcontentsline{toc}{section}{Appendix}
\appendix
\setcounter{section}{1} \setcounter{equation}{0}
\normalsize{
%
%
%
%
\subsection{Overview}
For the convenience of the reader we gather some further calculations in this appendix.
We define the simultaneous change of spatial and referential coordinates by a rigid rotation
\begin{itemize}
  \item for a scalar function
\begin{align}\label{trafoscal}
h^{\sharp} (\fg \xi) \colonequals h(\f Q^T \, \fg \xi) \, ,
\end{align}
  \item for a vector
\begin{align}\label{trafovect}
\fg \varphi^{\sharp} (\fg \xi) \colonequals \f Q \, \fg \varphi(\f Q^T \, \fg \xi) \, ,
\end{align}
  \item for a tensor of second order
\begin{align}\label{trafovect}
\fg \sigma^{\sharp} (\fg \xi) \colonequals \f Q \, \fg \sigma(\f Q^T \, \fg \xi) \, \f Q^T \, .
\end{align}
\end{itemize}
Note that the latter is motivated by
\begin{align}\label{phiphi}
\fg \varphi^{\sharp}(\fg \xi) \otimes  \fg \varphi^{\sharp}(\fg \xi) = \f Q \, \fg \varphi(\f Q^T \, \fg \xi) \otimes \f Q \, \fg \varphi(\f Q^T \, \fg \xi) =  \f Q \, [ \fg \varphi(\f Q^T \, \fg \xi) \otimes \fg \varphi(\f Q^T \, \fg \xi) ] \, \f Q^T \, .
\end{align}
For a tensor of third order the standard symbolic notation is not sufficient to express the simultaneous change of spatial and referential coordinates, which can be seen from
\begin{align}\label{phiphiphi}
\fg \varphi^{\sharp}(\fg \xi) \otimes  \fg \varphi^{\sharp}(\fg \xi)  \otimes  \fg \varphi^{\sharp}(\fg \xi) = \f Q \, \fg \varphi(\f Q^T \, \fg \xi) \otimes \f Q \, \fg \varphi(\f Q^T \, \fg \xi)  \otimes \f Q \, \fg \varphi(\f Q^T \, \fg \xi) =  \f Q \, [ \fg \varphi(\f Q^T \, \fg \xi) \otimes \fg \varphi(\f Q^T \, \fg \xi) ] \, \f Q^T \otimes \fg \varphi(\f Q^T \, \fg \xi) \, \f Q^T \, .
\end{align}
It is not possible to separate the three rotations $\f Q$ from a central dyadic product, such that e.g. $\m^{\sharp}(\fg \xi) = \f Q \, \m(\f Q^T \, \fg \xi) \, \f Q^T \, \f Q^T $ does not hold true. Therefore, we need to use index notation in case of third order tensors to be precise.
\subsection{Quadratic norms}
The quadratic norm of the second gradient of the displacement field reads
\begin{align}\label{NormGradGradu}
\norm{\GRADx{\Gradx{\ux}}}^2 & = \langle u_{a,bc} \, \f e_a \otimes \f e_b \otimes \f e_c \, , \, u_{i,jk} \, \f e_i \otimes \f e_j \otimes \f e_k \rangle = u_{a,bc} \, u_{i,jk} \, \delta_{ai} \, \delta_{bj} \, \delta_{ck}
 = \, u_{a,bc} \, u_{a,bc} \, ,
\end{align}
which is form-invariant under our transformation since
\begin{align}\label{NormGradGraduxi}
\norm{\GRADxi{\Gradxi{\uxi}}}^2  & = \norm{\GRADx{\f Q \, \Gradx{\uxi} \, \f Q^T} \, \f Q^T}^2 \notag \\
& = \langle Q_{ia} \, Q_{jb} \, Q_{kc} \,  u_{a,bc} \, \f e_i \otimes \f e_j \otimes \f e_k \, , \,  Q_{md} \, Q_{ne} \, Q_{pf} \,  u_{d,ef} \, \f e_m \otimes \f e_n \otimes \f e_p \rangle \notag \\
& = Q_{ia} \, Q_{jb} \, Q_{kc} \,  u_{a,bc} \,  Q_{md} \, Q_{ne} \, Q_{pf} \, u_{d,ef} \, \delta_{im} \, \delta_{jn} \, \delta_{kp} \notag \\
& =  \underbrace{Q^T_{ai} \, Q_{id}}_{\displaystyle \delta_{ad}} \, \underbrace{Q^T_{bj} \, Q_{je}}_{\displaystyle \delta_{be}} \, \underbrace{Q^T_{ck} \, Q_{kf}}_{\displaystyle \delta_{cf}}  \, u_{a,bc} \, u_{d,ef} \notag \\
& = u_{a,bc} \, u_{a,bc} = \norm{\GRADx{\Gradx{\ux}}}^2 \, .
\end{align}
 This invariance also holds for the gradient of the symmetric part
\begin{align}\label{NormGradsymGradu}
\norm{ & \GRADxi{\sym \Gradxi{\uxi}}}^2
= \norm{\GRADxi{\f Q \, \sym \Gradx{\ux} \, \f Q^T}}^2  \notag \\
& = \norm{\GRADx{\f Q \, \sym \Gradx{\uxi} \, \f Q^T} \, \f Q^T}^2 \notag \\
& = \langle Q_{ia} \, Q_{jb} \, Q_{kc} \, (\GRADx{\sym \Gradx{\ux}})_{abc} \, \f e_i \otimes \f e_j \otimes  \f e_k  \, , \, \notag \\
 & \hspace{6cm} Q_{md} \, Q_{ne} \, Q_{pf} \,  (\GRADx{\sym \Gradx{\ux}})_{def} \, \f e_m \otimes \f e_n \otimes \f e_p \rangle \notag \\
& = Q_{ia} \, Q_{jb} \, Q_{kc}  \, (\GRADx{\sym \Gradx{\ux}})_{abc} \, Q_{md} \, Q_{ne} \, Q_{pf}  \, (\GRADx{\sym \Gradx{\ux}})_{def} \, \, \delta_{im} \, \delta_{jn} \, \delta_{kp} \notag \\
 & =  \underbrace{Q^T_{ai} \, Q_{id}}_{\displaystyle \delta_{ad}} \, \underbrace{Q^T_{bj} \, Q_{je}}_{\displaystyle \delta_{be}} \, \underbrace{Q^T_{ck} \, Q_{kf}}_{\displaystyle \delta_{cf}}  \, (\GRADx{\sym \Gradx{\ux}})_{abc} \, (\GRADx{\sym \Gradx{\ux}})_{def} \notag \\
& = (\GRADx{\sym \Gradx{\ux}})_{abc} \, (\GRADx{\sym \Gradx{\ux}})_{abc} = \norm{\GRADx{\sym \Gradx{\ux}}}^2 \, .
\end{align}
Preparing
\begin{align}\label{NormGradtrsymGradu}
\norm{\Gradx{\tr ( \sym \Gradx{\ux}) \, \id }}^2 &= \norm{\Gradx{\tr (\Gradx{\ux}) \, \id }}^2 = \norm{\Gradx{ u_{a,a} \, \delta_{ij} \, \f e_i \otimes \f e_j } }^2 \notag \\
& = \langle u_{a,ab} \, \delta_{ij} \, \f e_i \otimes \f e_j \otimes \f e_b \, , \, u_{c,cd} \, \delta_{kl} \, \f e_k \otimes \f e_l \otimes \f e_d \rangle \notag \\
& = u_{a,ab} \, \delta_{ij} \, u_{c,cd} \, \delta_{kl} \, \delta_{ik} \, \delta_{jl} \, \delta_{bd} = u_{a,ab} \, \delta_{ij} \, u_{c,cb} \, \delta_{ij}
\end{align}
and in accordance with eq.\eqref{trsymskewgradu} we obtain
\begin{align}\label{NormGradtrsymGradu2}
\norm{\Gradxi{\tr ( \sym \Gradxi{\uxi}) \, \id }}^2 &= \norm{\Gradxi{\tr (\Gradxi{\uxi}) \, \id }}^2
\notag \\
& = \langle u_{a,ab} \, \delta_{ij} \, Q_{kb} \, \f e_i \otimes \f e_j \otimes \f e_k \, ,\, u_{c,cd} \, \delta_{mn} \, Q_{pd} \, \f e_m \otimes \f e_n \otimes \f e_p \rangle \notag \\
& = u_{a,ab} \, u_{c,cd} \, \delta_{ij}  \, \delta_{mn} \, Q_{kb} \, Q_{pd} \, \delta_{im} \, \delta_{jn} \, \delta_{kp} = u_{a,ab} \, u_{c,cd} \, \delta_{ij} \, \delta_{ij}  \, \underbrace{Q^T_{bk} \, Q_{kd}}_{\displaystyle \delta_{bd}} \notag \\
& = u_{a,ab} \, \delta_{ij} \, u_{c,cb} \, \delta_{ij} = \norm{\Gradx{\tr ( \sym \Gradx{\ux}) \, \id }}^2 \, ,
\end{align}
which is also form-invariant under our transformation. Therefore, from eq.\eqref{NormGradsymGradu} and eq.\eqref{NormGradtrsymGradu2} it follows that $\norm{\GRAD{\dev ( \sym \Grad{\f u(\f x)})}}^2$ is also form-invariant meaning that
\begin{align}
 \norm{\GRADxi{\dev ( \sym \Gradxi{\uxi})}}^2 = \norm{\GRADx{\dev ( \sym \Gradx{\ux})}}^2 \,.
\end{align}
Further, the norm of the curl of the right Cauchy-Green tensor $C = F^T \, F$ is form-invariant:
\begin{align}
 \norm{\Curlxi[C^{\sharp}(\xi)]}^2 &= \norm{\Curlxi[\Fxi^T \, \Fxi]}^2 =  \norm{\Curlxi[Q \, \Fx^T \, Q^T \, Q \, \Fx \, Q^T]}^2 =  \norm{\Curlxi[Q \, \Fx^T \, \Fx \, Q^T]}^2 \notag \\
&= \norm{Q \, \Curlx[\Fx^T \, \Fx] \, Q^T}^2 = \scal{Q \, \Curlx[\Fx^T \, \Fx] \, Q^T \, , \, Q \, \Curlx[\Fx^T \, \Fx] \, Q^T} \notag \\
&= \scal{\Curlx[\Fx^T \, \Fx] \, , \, Q^T Q \, \Curlx[\Fx^T \, \Fx] \, Q^T Q}
\notag \\
&=\norm{\Curlx[\Fx^T \, \Fx]}^2 = \norm{\Curlx[C(x)]}^2 \,.
\end{align}
\subsection{Transformation law for Kröner's incompatibility tensor}
We consider Kröners's incompatibility tensor $\inc(\sym \f P) \colonequals \Curl\left( (\Curl \sym \f P)^T \right) \in \Sym(3) \subseteq \Skalar^{3 \times 3}$.
By substituting $\sigmaxi = \sym \Pxi$ in eq.\eqref{CurlCurlSigma} we find
\begin{align}\label{CurlCurlSymP}
\Curlxi\Curlxi[\sym \Pxi] = \Curlxi [ \f Q \, \Curlx [\sym \Px]  \, \f Q^T ] = \f Q \, \Curlx[\Curlx[\sym \Px]] \, \f Q^T \,,
\end{align}
where eq.\eqref{symkxi} yielding $\sym \Pxi = \f Q \, \sym \Px \, \f Q^T$ has been used. Finally we obtain
\begin{align}\label{IncSymP}
\incxi(\sym \Pxi) & \colonequals \Curlxi \left( (\Curlxi[\sym \Pxi])^T \right) = \Curlxi \left( \f Q \, (\Curlx [\sym \Px])^T \, \f Q^T \right) \notag \\
& = \f Q \, \Curlx \left( (\Curlx[\sym \Px])^T \right) \, \f Q^T = \f Q \,  \incx(\sym \Px) \, \f Q^T \,.
\end{align}
Therefore, $\inc$ transforms similarly as any other second order tensor. Note again that inhomogeneous rotation fields are not allowed in \eqref{IncSymP}.
\subsection{Transformation law for the dislocation density tensor}
We also consider the dislocation density tensor $\Curl P \in \Skalar^{3 \times 3}$. By looking back at eq.\eqref{NeffTrafoLaw3} we note
\begin{align}
\Curlxi[P^{\sharp}(\xi)]=\Curlxi[\f Q \, \f P(\f Q^T \, \fg \xi) \, \f Q^T] =\f Q \, \Curlxi[\f P(\f x) \, \f Q^T] = \f Q \, \Curlx[\f P(\f x)] \, \f Q^T \,,
\end{align}
where again it is clear that the rotation $Q$ must be homogeneous.
\subsection{The Rayleigh product}\label{KapRayleighProduct}
The Rayleigh product \cite[p.~46]{bertram2012elasticity} for a (symmetric) second order tensor $C$ is defined by the algebraic operation
\begin{align}\label{RayleighTensor2}
Q * C \colonequals Q_{ia} \, Q_{jb} \, C_{ab} = Q \, C \, Q^T \,.
\end{align}
Following \cite{Auffray2014}, where the Rayleigh product is defined for arbitrary tensorial rank, it reads for a third order tensor $\m$
\begin{align}\label{RayleighTensor3}
Q * \m \colonequals Q_{ia} \, Q_{jb} \, Q_{kc} \,  \m_{abc} \,.
\end{align}
In contrast to \cite{bertram2012elasticity,Bertram2015} we abstain from introducing symbolic notation in eq.\eqref{RayleighTensor3}. Note that the Rayleigh product is the simple result of rotating the basis of any tensor component by $Q$ in accordance with eq.\eqref{QmapsEuklid}. A triad $d_i \otimes d_j \otimes d_k$ from rotated basis vectors $d_i = Q \, e_i$, $d_j = Q \, e_j$, and $d_k = Q \, e_k$ reads
\begin{align}\label{Rayleightriad}
 d_i \otimes d_j \otimes d_k = Q *  e_i \otimes e_j \otimes e_k = Q_{ia} \, e_a \otimes Q_{jb} \, e_b \otimes Q_{kc} \, e_c  =  Q_{ia} \, Q_{jb} \, Q_{kc} \, e_a \otimes e_b \otimes e_c \,.
\end{align}
Note that the $\sharp$-transformation generates the action of the Rayleigh-product when applied to objective quantities like $C$. By now it is clear that we are not directly working on the level of (isotropic) scalar or tensor functions. An algebraic characterization of isotropy for a {\bf reduced format of strain energies} can be obtained with help of the Rayleigh-product. Indeed, a free energy depending on a number of tensorial variables
\begin{align}\label{FreeEnerNTensors}
\Psi = \Psi(\underbrace{C}_{\rm 2.order}, \, \underbrace{\nabla C}_{\rm 3.order}, \, \nabla^2 C, \, \ldots , \, \underbrace{\nabla^n C}_{\rm 2+n. \, order}) = \Psi(C, \, \GRADx{C}, \, \GRADx{\GRADx{C}}, \, \ldots) \,,
\end{align}
induces an isotropic formulation if and only if
\begin{align}\label{PsiRayleigh}
\Psi(Q^{(2)}\!\! * C, \, Q^{(3)}\!\! * \nabla C, \, \ldots, \, Q^{(2+n)}\!\! * \nabla^n C) =  \Psi(C, \, \nabla C, \, \nabla^2 C, \, \ldots , \, \nabla^n C) \quad \forall Q \in \SO(3) \, ,
\end{align}
where $Q^{(n)} * \textsl{M}$ denotes the Rayleigh product acting on n-th. order tensors.\footnote{Since in this {\bf special representation} the rotations are appearing already derivative free, the condition could be formulated with inhomogeneous rotation fields without altering the result.} Similarly, for the linearized kinematics, isotropy, and second strain gradient material we need to have
\begin{align}\label{PsiRayleigh2}
\Psi_{\rm lin}(Q^{(2)}\!\! * \sym \nabla u, \, Q^{(3)}\!\! * \nabla \sym \nabla u) =  \Psi_{\rm lin}(\sym \nabla u, \, \nabla \sym \nabla u) \quad \forall Q \in \SO(3) \, .
\end{align}
As we have seen in the main body of this work, special gradient elasticity formulations employ certain reduced representations. For these we note the isotropy conditions as well. For the indeterminate couple stress model we obtain
\begin{align}\label{PsiRayleigh3}
\Psi_{\rm lin}(Q^{(2)}\!\! * \sym \nabla u, \, Q^{(2)}\!\! * \nabla \axl(\skew \nabla u)) =  \Psi_{\rm lin}(\sym \nabla u, \, \nabla \axl(\skew \nabla u)) \quad \forall Q \in \SO(3) \, .
\end{align}
And for the version with symmetric stresses we get
\begin{align}\label{PsiRayleigh4}
\Psi_{\rm lin}(Q^{(2)}\!\! * \sym \nabla u, \, Q^{(2)}\!\! * \Curl \sym \nabla u ) =  \Psi_{\rm lin}(\sym \nabla u, \, \Curl \sym \nabla u ) \quad \forall Q \in \SO(3) \, .
\end{align}
The following remark is in order. It is true that the second order tensor $\Curl \sym \nabla u $ is a linear combination of components of the third order tensor $\nabla \sym \nabla u$. The isotropy condition based on $\nabla \sym \nabla u$ needs the Rayleigh product with $Q^{(3)}$ (i.e. 3 rotations) in \eqref{PsiRayleigh2}, while $\Curl \sym \nabla u$ necessitates to take $Q^{(2)}$ in \eqref{PsiRayleigh4}. This apparent inconsistency is resolved by observing that $\Curl \sym \nabla u $ is not just any linear combination of $\nabla \sym \nabla u$, but such that it holds
\begin{align}\label{RayleighCurl}
Q^{(2)}\!\! * \Curl \sym \nabla u = Q^{(2)}\!\! * \widehat{\L}_{\rm Curl} \nabla \sym \nabla u = \widehat{\L}_{\rm Curl} (Q^{(3)}\!\! *\nabla \sym \nabla u) \, ,
\end{align}
where $\widehat{\L}_{\rm Curl} : \Skalar^{3 \times 3 \times 3} \mapsto   \Skalar^{3 \times 3}$ is such that
\begin{align}\label{CurlAbbild}
\Curl \sym \nabla u = \widehat{\L}_{\rm Curl} \nabla \sym \nabla u \, .
\end{align}
In general, for arbitrary  $\L : \Skalar^{3 \times 3 \times 3} \mapsto   \Skalar^{3 \times 3}$
\begin{align}\label{nablasymnablauAbbild}
\L (Q^{(3)}\!\! * \nabla \sym \nabla u) \neq Q^{(2)}\!\! * \L \nabla \sym \nabla u \, .
\end{align}
The same observation applies to $ \nabla \axl(\skew \nabla u)$ in \eqref{PsiRayleigh3}.
}
\end{document}